\documentclass[a4paper,10pt]{article}

\usepackage{latexsym}

\usepackage{amssymb}

\def \zgran{\displaystyle}
\def \zpizq{\left(}
\def \zcizq{\left[}
\def \zpder{\right)}
\def \zcder{\right]}

\def \za{\alpha}
\def \zb{\beta}
\def \zg{\gamma}

\def \ze{\varepsilon}

\def \zh{\theta}

\def \zl{\lambda}
\def \zm{\mu}
\def \zn{\nu}
\def \zx{\xi}

\def \zp{\pi}

\def \zr{\rho}

\def \zs{\sigma}

\def \zt{\tau}

\def \zf{\varphi}
\def \zfi{\phi}

\def \zq{\psi}
\def \zw{\omega}

\def \zL{\Lambda}

\def \zF{\Phi}
\def \zQ{\Psi}
\def \zW{\Omega}

\def \zima{\imath}

\def \zlma{\ell}

\def \zsu{\sum}
\def \zpr{\prod}

\def \zin{\cap}

\def \zun{\cup}

\def \zex{\wedge}

\def \zdi{\oplus}

\def \zte{\otimes}

\def \zpu{\cdot}
\def \zpor{\times}
\def \zci{\circ}

\def \zmei{\leq}
\def \zmai{\geq}
\def \zco{\subset}

\def \zcco{\supset}

\def \zpe{\in}

\def \znoi{\neq}

\def \znope{\not\in}

\def \zpar{\partial}
\def \zinf{\infty}

\def \zfl{\rightarrow}

\def \zbv{\mid}
\def \zdbv{\parallel}
\def \z/{\over}

\begin{document}

{\bf THE LOCAL PRODUCT THEOREM FOR BIHAMILTONIAN STRUCTURES}
\bigskip

{\it \noindent Francisco-Javier Turiel

\noindent Geometr{\'\i}a y Topolog{\'\i}a, Facultad de Ciencias,
Campus de Teatinos, 29071 M\'alaga, Spain

\noindent e-mail: turiel@agt.cie.uma.es}
\vskip 1truecm

{\bf Abstract}

In this work one proves that, around each point of a dense open
set (regular points), a real analytic or holomorphic bihamiltonian
structure decomposes into a product of a Kronecker bihamiltonian
structure and a symplectic one if a necessary condition on the
characteristic polynomial of the symplectic factor holds. Moreover
we give an example of bihamiltonian structure for showing that
this result does not extend to the $C^{\zinf}$- category.

Thus a classical problem on the geometric theory of bihamiltonian
structures is solved at almost every point.
\bigskip

{\bf Introduction}

Given two Poisson structures $\zL,\zL_{1}$ on a real (at least $C^{\zinf}$)
or complex (holomorphic) manifold $M$, following Magri \cite{MAG}
one will say that $(\zL,\zL_{1})$ is a bihamiltonian structure (or
that $\zL,\zL_{1}$ are compatible) if $\zL+\zL_{1}$ is a Poisson
structure as well. Bihamiltonian structures are a useful tool for
dealing with some differential equations many of them with a
physical meaning; besides they are interesting from the
geometric viewpoint too that will be the case here.

The algebraic classification of the pairs of bivectors on a finite
dimensional real or complex vector space was given by Gelfand and
Zakharevich in \cite{GEB}. Essentially each pair decomposes into
the product of a Kronecker pair and a symplectic one (see
\cite{GEB,TUD}). Therefore  it is natural to ask whether this
decomposition into a product Kronecker-symplectic holds, at least
locally, for bihamiltonian structures as well, which would be an
important steep toward their classification. One recalls that
Kronecker bihamiltonian structures are intimately related to
Veronese webs (see \cite{TUD} for an exposition of the local
theory of Veronese webs and its relationship with Kronecker
bihamiltonian structures), whereas the local classification of
symplectic bihamiltonian structures, that is of pairs of
compatible symplectic forms, is known at almost every point (see
\cite{T89,TU}).

The chief goal of this work is to show that, around every point of
a dense open set (regular points), a real analytic or holomorphic
bihamiltonian structure decomposes into a product
Kronecker-symplectic if a necessary condition on the
characteristic polynomial of the symplectic factor holds (theorem
7.1). Moreover we exhibit an example of $C^{\zinf}$- bihamiltonian
structure for which this result fails. These results have been
annoced in \cite{TU11A,TU11B}.

As main tool for this purpose, to any bihamiltonian structure we
associate a new object called a Veronese flag, which generalizes the
notion of Veronese web introduced by Gelfand and Zakharevich in
\cite{GEB} (codimension one) and later on by others authors \cite{PAA,TUC}
(higher codimension). Roughly speaking the crucial point is to
show that, about each regular point, a Veronese flag is the
product of a Veronese web and a pair of compatible symplectic
forms. For that one has to prove, in a indirect way, the existence
of solutions of some differential equations not explicitly
formulated. In the complex case they are always ordinary whereas
in some real cases we have to deal with systems of partial
derivative equations, which may contain the Lewy's example \cite{LW}
as sub-system; thus the result fails in the $C^{\zinf}$ category. In
the real analytic case a method of complexification  transforms
the real problem on a complex one.
Once the decomposition of Veronese flags established, that of
bihamiltonian structures follows from it with a little extra-work.

The study of bihamiltonian structures at not regular points rather
belongs to the theory of singularities and, in spite of its great
interest, will be not considered here.

The present text consists of eight sections plus an appendix. In the
first one the Veronese flag, as quotient of a bihamiltonian structure,
and the bihamiltonian structure over a $(1,1)$-tensor field and
a foliation, which gives a simple method for constructing
bihamiltonian structures, are introduced.

Sections 2, 3 and 4, this last one rather technical, are devoted to prove
a local decomposition theorem for Veronese flags with only one eigenvalue.
In section 5 real Veronese flags, without real eigenvalue, are dealt with
by reducing them to holomorphic ones through the analyticity.

In section 6 one shows that, locally, Veronese flags are the fibered
product of those with just one eigenvalue, and in section 7 we prove
the theorem of local decomposition for real analytic or holomorphic
bihamiltonian structures. A $C^{\zinf}$ counter-example to this last
result is given in section 8.

Finally, in the appendix one proves a well known result on $(1,1)$-tensor
fields belonging to the folklore but without many accessible proofs.
\bigskip

{\bf 1. The quotient of a bihamiltonian structure}

As it well known to any Kronecker bihamiltonian structure one may associate a
Veronese web on the local quotient of the support manifold (see \cite{GEA,PAA,TUC}).
Here we will
associate a new structure, called a Veronese flag and defined on a local
quotient of the support manifold also, to a very large class of bihamiltonian structures.

{\it From now on all structured considered will be real $C^{\zinf}$ or complex holomorphic
unless another thing is stated.}
\bigskip

{\bf 1.1. The main construction.}

On a manifold $P$ consider a foliation ${\mathcal F}$ (that is an
involutive distribution) of positive codimension and a morphism of
vector bundles $\zlma:{\mathcal F}\zfl TP$. If $\za$ is a $s$-form
on an open set $B$ of $P$, then $\zlma^{*}\za$ ( we will write
$\za\zci \zlma$ as well) can be regarded as
$s$-form with domain $B$ on the leaves of ${\mathcal F}$. Let
$G:TP\zfl TP$ be a prolongation of $\zlma$; then
$(G^{*}\za)_{\zbv{\mathcal F}}$ equals $\zlma^{*}\za$. On the
other hand if $\zlma^{*}\za$ is closed on  ${\mathcal F}$ for
every closed $1$-form $\za$ such that $Ker\za\zcco{\mathcal F}$,
then the restriction of the Nijenhuis torsion $N_G$ of $G$ to
${\mathcal F}$ does not depend on the prolongation $G$ of $\zlma$
(see lemma 2.2 of \cite{TUD}) and it will called the {\it Nijenhuis
torsion $N_{\zlma}$ of $\zlma$}.

Let ${\mathcal A}(p)$, $p\zpe P$, be the largest $\zlma$-invariant
vector subspace of ${\mathcal F}(p)$.

We will say that the pair $({\mathcal F}, \zlma)$ is a {\it weak
Veronese flag} if the following three conditions hold:

\noindent 1) $\zlma^{*}\za$ is closed on ${\mathcal F}$ for every
closed $1$-form $\za$ such that $Ker\za\zcco{\mathcal F}$,

\noindent 2) $N_{\zlma}=0$,

\noindent 3) $dim{\mathcal A}(p)$ does not depend on $p$.

First of all let us see that the distribution ${\mathcal
A}=\zun_{p\zpe P}{\mathcal A}(p)$ is a foliation  when $({\mathcal
F}, \zlma)$ is a weak Veronese flag. Given any point $q\zpe P$ the
morphism $\zlma:{\mathcal F}(q)\zfl T_{q}P$ projects in a morphism
$\zf_{q}:{\mathcal F}(q)/{\mathcal A}(q)\zfl T_{q}P/{\mathcal
A}(q)$ without non-zero $\zf_{q}$-invariant vector subspace. Thus
$({\mathcal F}(q)/{\mathcal A}(q),\zf_{q})$ defines an algebraic
Veronese web $w_{q}$ on $T_{q}P/{\mathcal A}(q)$ of codimension
$\zmai  1$ by setting $w_{q}(t)=(\zf_{q}+tI)({\mathcal
F}(q)/{\mathcal A}(q))$, such that $w_{q}(\zinf)={\mathcal
F}(q)/{\mathcal A}(q)$. Moreover, as it is well known, if
$a_{1},...,a_{k}$, $k=dim(T_{q}P/{\mathcal A}(q))$, are non-equal
scalars then $w_{q}(a_{1})\zin...\zin w_{q}(a_{k})=\{ 0\}$, so
${\mathcal A}(q)=\zin_{j=1}^{k}((\zlma+a_{j}I){\mathcal F}(q))$.

On the other hand, if $Ker(\zlma+aI)=\{0\}$ on an open set then
$(\zlma+aI){\mathcal F}$ is involutive on this set. For showing
this last assertion we need the following result transcription of
lemma 2.1 of \cite{TUD}.
\bigskip

{\bf Lemma 1.1.} {\it Consider a $1$-form $\zr$, a $(1,1)$- tensor
field $H$ and two vector fields $X,Y$ on a manifold, then
$(d(\zr\zci H))(HX,Y)+(d(\zr\zci H))(X,HY) = d\zr (HX, HY) +
d(\zr\zci H^{2})(X,Y) + \zr (N_{H}(X,Y))$.}
\bigskip

Let $G$ be a (local) prolongation of $\zlma$. If $\zm$ is a closed
$1$-form and $Ker\zm\zcco {\mathcal F}$ then
$Ker(\zm\zci(G+aI)^{-1})\zcco (G+aI){\mathcal F}$ and by lemma 1.1
applied to $\zm\zci(G+aI)^{-1}$ and $(G+aI)$ one has
$d(\zm\zci(G+aI)^{-1})((G+aI){\mathcal F},(G+aI){\mathcal F})=
-d(\zm\zci(G+aI))({\mathcal F},{\mathcal F}) -\zm(N_{G}({\mathcal
F},{\mathcal F}))=0$ therefore $d(\zm\zci(G+aI)^{-1})_{\zbv
(G+aI){\mathcal F}}=0$, whence the involutivity of
$(G+aI){\mathcal F}$.

Since whichever $p\zpe P$ always there exist non-equal scalars
$a_{1},...,a_{k}$ such that $Ker(\zlma+a_{j}I)(p)=\{0\}$,
$j=1,...,k$, around this point ${\mathcal
A}=\zin_{j=1}^{k}((\zlma+a_{j}I){\mathcal F})$; so ${\mathcal A}$
is a foliation, called {\it the axis of the flag $({\mathcal
F},\zlma)$} from now on.

Let $\zp:P\zfl N$ be a local quotient of $P$ by ${\mathcal A}$;
then ${\bar w}(t)=\zp_{*}((\zlma+tI){\mathcal F})$ is a foliation
whose codimension equals that of ${\mathcal F}$ and ${\bar
w}=\{{\bar w}(t)\zbv t\zpe{\mathbb K}\}$ is a Veronese web on $N$.
Indeed, given $q\zpe P$ and ${\bar q}\zpe N$ such that
$\zp(q)={\bar q}$ then ${\bar w}({\bar q})=\{{\bar w}({\bar
q})(t)\zbv t\zpe{\mathbb K}\}$ is the algebraic Veronese web
defined by $({\mathcal F}(q)/{\mathcal A}(q),\zf_{q})$ when
$T_{q}P/{\mathcal A}(q)$ is identified to $T_{\bar q}N$; so ${\bar
w}$ is an algebraic Veronese web at each point of $N$. Moreover
${\bar w}(\zinf)=\zp_{*}({\mathcal F})$, which is a foliation;
therefore by proposition 2.1 of \cite{TUD} the family ${\bar w}$ is a
Veronese web.

Thus if ${\bar\zlma}:{\bar\mathcal F}\zfl TN$, where
${\bar\mathcal F}=\zp_{*}({\mathcal F})$, is the morphism
canonically associated to ${\bar w}$ then $({\bar\mathcal
F},{\bar\zlma})$ is the projection of $({\mathcal F},\zlma)$.
\bigskip

{\bf Lemma 1.2.} {\it Consider a weak Veronese flag $({\mathcal
F},\zlma)$ and for every integer $k\zmai 0$ set
$g_{k}=trace(({\zlma}_{\zbv\mathcal A})^{k})$. Then
$kdg_{k+1}=(k+1)dg_{k}\zci\zlma$ on ${\mathcal F}$.}
\bigskip

{\bf Proof.} As the problem is local one may extend $\bar\zlma$ to
a flat diagonalizable tensor field $J$ and consider an extension
$G$ of $\zlma$ projecting in $J$. Then $ImN_{G}\zco{\mathcal A}$
and $trace(G^{k})=g_{k}+c_{k}$ where $c_{k}\zpe{\mathbb K}$. Since
$trace(H_{1}\zci H_{2})=trace(H_{2}\zci H_{1})$ for every vector
field $X$ one has: $kd(trace(G^{k+1}))(X)=k(k+1)trace(G^{k}\zci
L_{X}G)=k(k+1)trace(G^{k-1}\zci
L_{GX}G)-k(k+1)trace(N_{G}(X,\quad))= (k+1)d(trace(G^{k}))(GX)
-k(k+1)trace(N_{G}(X,\quad))$.

But $trace(N_{G}(X,\quad))=0$ when $X\zpe{\mathcal F}$ because
$N_{G}({\mathcal F},{\mathcal F})=0$ and
$Im(N_{G}(X,\quad))\zco{\mathcal A}\zco{\mathcal F}$. $\square$

Now let $\zw,\zw_{1}$ be a couple of $2$-forms defined on
$\mathcal A$. One will say that $({\mathcal F}, \zlma,
\zw,\zw_{1})$ is a {\it Veronese flag} on $P$ if:

\noindent 1) $({\mathcal F}, \zlma)$ is a weak Veronese flag.

\noindent 2) $\zw$ is symplectic on $\mathcal A$, $\zw_{1}$ closed and
$\zw_{1}=\zw(\zlma,\quad)$ [that is $\zw_{1}(X,Y)=\zw(\zlma X,Y)$].

\noindent 3) Whenever $f$ is a function on an open set of $P$ such that
$\zlma^{*}df$ is closed on $\mathcal F$, then $L_{X_{f}}\zlma=0$ where
$X_{f}$ is the $\zw$-hamiltonian of $f$ along $\mathcal A$.
\bigskip

{\bf Remark.} Given, on a manifold,  a foliation ${\mathcal  G}$, a tensor field
${\mathcal  T}$ defined along ${\mathcal  G}$ and a ${\mathcal  G}$-foliate vector field
$X$, then the Lie derivative $L_{X}{\mathcal  T}$ is defined as a tensor field along
${\mathcal  G}$; moreover the flow of $X$preserves ${\mathcal  T}$ if and only if
$L_{X}{\mathcal  T}=0$. In condition 3) above $X_{f}$ is tangent to ${\mathcal  A}\zco {\mathcal  F}$
so ${\mathcal  F}$-foliate. Obviously this condition  implies $L_{X_{f}}\zw_{1}=0$.

When ${\mathcal A}=0$, Veronese web and Veronese flag are equivalent notions.

By technical reasons we need the following definition. Given $p_{0}\zpe P$ we
will say that $({\mathcal F}, \zlma,\zw,\zw_{1})$ is a {{\it Veronese flag at point
$p_{0}$} when 1) and 2) hold but 3) is replaced by:

\noindent 3') for any function $f$ defined on an open set $p_{0}\zpe B\zco P$ such that
$\zlma^{*}df$ is closed on $\mathcal F$, then $L_{X_{f}}\zlma=0$ on an open
set $p_{0}\zpe B'\zco B$.

Let us recall some facts about pairs of bivectors on real or
complex vector spaces (see \cite{GEB} and section 1.2 of \cite{TUD}).
Consider a pair of
bivectors $(\zl,\zl_{1})$ on a finite  dimensional vector space
$W$. By definition {\it the rank of $(\zl,\zl_{1})$} is the
maximum of the ranks of $(1-t)\zl+t\zl_{1}$, $t\zpe \mathbb K$,
and one has $rank(\zl,\zl_{1})= rank((1-t)\zl+t\zl_{1})$ except
for a finite number of scalars $t$, which is $\zmei \frac{dim W}
{2}$. A pair $(\zl,\zl_{1})$ is called {\it maximal} when
$rank(\zl)=rank(\zl_{1})=rank(\zl,\zl_{1})$. Given an odd
dimensional vector space $U$, the action of $GL(U)$ on
$(\zL^{2}U)\zpor(\zL^{2}U)$ possesses one dense open orbit, whose
elements are named {\it Kronecker elementary pairs}; they are
maximal and their rank equals $dimU-1$. According to the
classification by Gelfand and Zakharevich (see \cite{GEB} and propositions 1.4
and 1.5 of \cite{TUD}), every maximal pair decomposes  into a product of
Kronecker elementary pairs $(U_{j},\zm_{j},\zm_{1j})$, $j=1,...r$,
where $r=corank(\zl,\zl_{1})$, and a symplectic pair
$(U',\zm',\zm_{1}')$; moreover these factors are unique up to
isomorphism or change of order.

A bihamiltonian structure on a manifold is called {\it Kronecker}
when at each point its algebraic model is a product of Kronecker
elementary pairs only, and {\it symplectic} if at every point its
algebraic model only has the symplectic factor.

On a real or complex $m$-manifold $M$ consider a bihamiltonian structure $(\zL,\zL_{1})$
such that:

\noindent 1) $(\zL,\zL_{1})$ is maximal, that is every $(\zL(p),\zL_{1}(p))$, $p\zpe M$,
is maximal,

\noindent 2) the rank of $(\zL,\zL_{1})$ and the dimension of the the symplectic factor at
each point are constant.

As before set $r=corank(\zL,\zL_{1})$ and let $2m'$ be the dimension of the symplectic factor.
Since $r$ is the number of Kronecker elementary factors, $m+r$ is even and one may set
$m=2m'+2n-r$. Note that, at every point, $2n-r$ equals the sum of the dimensions of the
Kronecker elementary factors (warning these last dimensions could depend on the point).

Our next aim is locally to associate a Veronese flag in dimension
$2m'+n$ to $(\zL,\zL_{1})$. For each $p\zpe M$ let ${\mathcal
A}_{1}(p)$ be the intersection of all vector subspaces $Im(\zL+t\zL_{1})(p)$,
$t\zpe \mathbb K$, such that $rank(\zL+t\zL_{1})(p)=m-r$. From the
algebraic model follows that $dim{\mathcal A}_{1}(p)=m-n=2m'+n-r$,
which defines a foliation ${\mathcal A}_{1}$ called {\it the
(primary) axis of $(\zL,\zL_{1})$}. Indeed, given $p\zpe M$ one can
chose non-equal scalars $t_{1},...,t_{n}$ such that
$rank(\zL+t_{j}\zL_{1})(p)=m-r$, $j=1,...,n$; in particular
$\zin_{j=1}^{n}Im(\zL+t_{j}\zL_{1})(p)={\mathcal A}_{1}(p)$. By
continuity $rank(\zL+t_{j}\zL_{1})=m-r$, $j=1,...,n$ and
$\zin_{j=1}^{n}Im(\zL+t_{j}\zL_{1})={\mathcal A}_{1}$ around $p$.

It is not hard to see that ${\mathcal A}_{1}\zco Im\zL_{1}$ and
$dim(Im(\zL+t\zL_{1})+{\mathcal A}_{1})=m-r$, $t\zpe\mathbb K$.
Set ${\tilde w}(t)=Im(\zL+t\zL_{1})+{\mathcal A}_{1}$, $t\zpe
\mathbb K$; then ${\tilde w}=\{ {\tilde w}(t)\zbv t\zpe \mathbb K
\}$ is a family of foliations of codimension $r$ whose limit at
each point, when $t\zfl\zinf$, is $Im\zL_{1}$. Indeed, given
$p\zpe M$ and $t_{0}\zpe \mathbb K$ consider functions
$f_{1},...,f_{k}$ and vector fields $X_{1},...,X_{m-r-k}$ tangent
to ${\mathcal A}_{1}$, all of them defined around $p$, such that
$\{
(\zL+t_{0}\zL_{1})(df_{1},\quad),...,(\zL+t_{0}\zL_{1})(df_{k},\quad),X_{1},...,X_{m-r-k}\}$
at $p$ is a basis of ${\tilde w}(t_{0})(p)$. By continuity $\{
(\zL+t\zL_{1})(df_{1},\quad),...,(\zL+t\zL_{1})(df_{k},\quad),X_{1},...,X_{m-r-k}\}$
is a basis of ${\tilde w}(t)(q)$ when $(q,t)$ is close to
$(p,t_{0})$ on $M\zpor{\mathbb K}$, so ${\tilde w}$ is a family of
distributions. But the set $D=\{ (q,t)\zpe M\zpor{\mathbb K}\zbv
rank(\zL+t\zL_{1})(q)=m-r\}$ is dense and open and, obviously,
${\tilde w}$ is a family of foliations on $D$, therefore ${\tilde
w}$ is a family of foliations on $M\zpor{\mathbb K}$. Finally note
that ${\mathcal A}_{1}\zco Im\zL_{1}$ and
$Im(\zL+t\zL_{1})=Im(s\zL+\zL_{1})$ when $s=t^{-1}$.

Let $N$ be the local quotient of $M$ by ${\mathcal A}_{1}$, which is a $n$-dimensional manifold,
and $\zp_{N}:M\zfl N$ the canonical projection. Then
${\bar w}=\{{\bar w}(t)=(\zp_{N})_{*}{\tilde w}(t)\zbv t\zpe {\mathbb K}\}$ is a Veronese web on $N$
of codimension $r$. Indeed, the algebraic model shows that, at each point of $N$, the family of
foliations ${\bar w}$ is an algebraic Veronese web. On the other hand its limit when $t\zfl\zinf$
equals $(\zp_{N})_{*}(Im\zL_{1})$, which is a foliation too; so ${\bar w}$ is a Veronese web
(see proposition 2.1 of \cite{TUD}).

The Poisson structure $\zL$ is given by a symplectic form
${\tilde\zw}$ defined on $Im\zL$ while $\zL_{1}$ is given by a
symplectic form form ${\tilde\zw}_{1}$ on $Im\zL_{1}$. Therefore
the restricted $2$-forms ${\tilde\zw}_{\zbv{\mathcal A}_{1}}$ and
${\tilde\zw}_{1\zbv{\mathcal A}_{1}}$ are closed; besides (see
proposition 1.4 of \cite{TUD}) $Ker({\tilde\zw}_{\zbv{\mathcal
A}_{1}})=Ker({\tilde\zw}_{1\zbv{\mathcal A}_{1}}) =\zL({\mathcal
A}'_{1},\quad)=\zL_{1}({\mathcal A}'_{1},\quad)$ where ${\mathcal
A}'_{1}$ is the annihilator of ${\mathcal A}_{1}$ and
$dim(Ker({\tilde\zw}_{\zbv{\mathcal A}_{1}}))=n-r$. Thus
${\mathcal A}_{2}=Ker({\tilde\zw}_{\zbv{\mathcal A}_{1}})$ is a
foliation of dimension $n-r$, which will be called {\it the
secondary axis of $(\zL,\zL_{1})$}, and ${\mathcal
A}_{2}\zco{\mathcal A}_{1}$.

Let $P$ be the local quotient of $M$ by ${\mathcal A}_{2}$ and
$\zp_{P}:M\zfl P$ the canonical projection; then $dimP=2m'+n$,
${\mathcal A}_{1}$ projects into a $2m'$-dimensional foliation
${\mathcal A}$ and ${\tilde\zw}_{\zbv{\mathcal A}_{1}}$,
${\tilde\zw}_{1\zbv{\mathcal A}_{1}}$ in two symplectic forms
$\zw$, $\zw_{1}$ on ${\mathcal A}$.
Moreover $\zL$ projects in the Poisson structure defined by
$\mathcal{A}$ and $\zw$, whereas $\zL_{1}$ does in the Poisson
structure defined by $\mathcal{A}$ and $\zw_{1}$.
Let ${\mathcal F}$ be the
$r$-codimensional foliation on $ P$ projection of $Im\zL_{1}$.
Obviously the local quotient of $P$ by ${\mathcal A}$ is
identified in a natural way to $N$ and $\zp\zci\zp_{P}=\zp_{N}$
where $\zp:P\zfl N$ is the canonical projection. In short we have
three of the four elements of a Veronese flag on $P$. Let us
construct the fourth one.

As $\zL({\mathcal A}'_{1},\quad)=\zL_{1}({\mathcal
A}'_{1},\quad)={\mathcal A}_{2}$ and ${\mathcal A}'_{1}$ contains
$Ker\zL$ and $Ker\zL_{1}$, the Poisson structures $\zL$, $\zL_{1}$
give rise to two isomorphisms ${\tilde\zl}$, ${\tilde\zl}_{1}$
from $\frac{T^{*}M} {{\mathcal A}'_{1}}$ to $\frac{Im\zL}
{{\mathcal A}_{2}}$ and $\frac{Im\zL_{1}} {{\mathcal A}_{2}}$
respectively, by setting ${\tilde\zl}([\za])=[\zL(\za,\quad)]$ and
${\tilde\zl}_{1}([\za])=[\zL_{1}(\za,\quad)]$. Thus
${\tilde\zlma}={\tilde\zl}\zci{\tilde\zl}_{1}^{-1}$ is a
monomorphism from $\frac{Im\zL_{1}} {{\mathcal A}_{2}}$ to
$\frac{TM} {{\mathcal A}_{2}}$ whose image equals $\frac{Im\zL}
{{\mathcal A}_{2}}$. By construction ${\tilde\zlma}$ is an
invariant of $(\zL,\zL_{1})$ and, for every $q\zpe M$, there
exists a monomorphism $\zf:(\zp_{P})_{*}(Im\zL_{1}(q))\zfl
T_{\zp_{P}(q)}P$ with $Im\zf=(\zp_{P})_{*}(Im\zL(q))$ that is the
projection of ${\tilde\zlma}(q)$; moreover from the algebraic
model follows that $\zw_{1}(u,v)=\zw(\zf u,v)$, $u,v\zpe{\mathcal
A}(\zp_{P}(q))$, and
$\zL_{1}(\zp^{*}_{P}\zf^{*}\zb,\quad)=\zL(\zp^{*}_{P}\zb,\quad)$,
$\zb\zpe T^{*}_{\zp_{P}(q)}P$ [note that $\zL_{1}$ can be regarded
as a linear map from $(Im\zL_{1})^{*}$ to $TM$ since
$\zL_{1}({\tilde\zb},\quad)=0$ whenever
${\tilde\zb}(Im\zL_{1})=0$].

For proving that, in fact, ${\tilde\zlma}$ projects into a
suitable morphism ${\zlma}:{\mathcal F}\zfl TP$ we will need some
extra-work, essentially local, which allows us to do it around the
points of $M$. Therefore, given non-equal and non-vanishing
scalars $a_{1},...,a_{n-r},a$ we may assume the existence on $N$
of coordinates $(x_{1},...,x_{n})$, closed $1$-forms
$\za_{1},...,\za_{r}$ and a $(1,1)$-tensor field $J$ such that
$dx_{j}\zci J=a_{j}dx_{j}$, $j=1,...,n-r$, $dx_{j}\zci J=adx_{j}$,
$j=n-r+1,...,n$, $Ker(\za_{1}\zex...\zex\za_{r})={\bar w}(\zinf)$,
$d(\za_{k}\zci J)\zex\za_{1}\zex...\zex\za_{r}=0$, $k=1,...,r$,
and that

\centerline{$\zg(t)=(\zpr_{j=1}^{n-k}(t+a_{j}))(t+a)^{k}(\za_{1}\zci(J+tI)^{-1})
\zex...\zex(\za_{r}\zci(J+tI)^{-1})$}

\noindent represents $\bar w$ (see theorem 2.1 of \cite{TUD}).

By identifying $\zt$ and $\zp^{*}_{N}\zt$, any $k$-form $\zt$, defined on open set of $N$,
can be regarded as an $\mathcal A_{1}$-basic $k$-form on an open set of $M$. Thus
$dx_{1},...,dx_{n}$ span ${\mathcal A}'_{1}$,
$Kerdx_{j}\zcco Im(\zL-a_{j}\zL_{1})$ whence
$\zL(dx_{j},\quad)=a_{j}\zL_{1}(dx_{j},\quad)$, $j=1,...,n-r$, and
$Kerdx_{j}\zcco Im(\zL-a\zL_{1})$ whence
$\zL(dx_{j},\quad)=a\zL_{1}(dx_{j},\quad)$, $j=n-r+1,...,n$.

On the other hand from the algebraic model at each point follows that two functions of
$(x_{1},...,x_{n})$ are always in involution for both $\zL$ and $\zL_{1}$, and the
families $\{dx_{1},...,dx_{n-r}, \za_{1}\zci J^{-1},...,\za_{r}\zci J^{-1}\}$ and
$\{dx_{1},...,dx_{n-r}, \za_{1},...,\za_{r}\}$ are linearly independent everywhere.
Consequently around each $p\zpe M$ one may chose functions
$y_{1},...,y_{n-r},z_{1},...,z_{2m'}$ such that
$(x_{1},...,x_{n},y_{1},...,y_{n-r},z_{1},...,z_{2m'})$ is a system of coordinates
and $\zL$ is given by $\za_{1}\zci J^{-1},...,\za_{r}\zci J^{-1}$ and the closed $2$-form
${\tilde\zW}=\zsu_{j=1}^{n-r}dx_{j}\zex dy_{j}+\zsu_{k=1}^{2m'}dz_{2k-1}\zex dz_{2m'}$.

Therefore ${\mathcal A}_{2}$ is spanned by $\zpar /\zpar y_{1},...,\zpar /\zpar y_{n-r}$
since $\zL(dx_{j},\quad)=\zpar /\zpar y_{j}$, $j=1,...,n-r$, and
$\zL(\za_{k}\zci J^{-1},\quad)=0$, $k=1,...,r$.

By the same reason around each point $p\zpe M$ there exist functions
$y'_{1},...,y'_{n-r}$, $z'_{1},...,z'_{2m'}$ such that
$(x_{1},...,x_{n},y'_{1},...,y'_{n-r},z'_{1},...,z'_{2m'})$ is a system of coordinates
while $\zL_{1}$ is given by $\za_{1},...,\za_{r}$ and the closed $2$-form
${\tilde\zW}_{1}=\zsu_{j=1}^{n-r}dx_{j}\zex dy'_{j}+\zsu_{k=1}^{2m'}dz'_{2k-1}\zex dz'_{2m'}$.

But $\zL(dx_{j},\quad)=\zpar/\zpar y_{j}$ and $\zL_{1}(dx_{j},\quad)=\zpar/\zpar y'_{j}$
whence $\zpar/\zpar y_{j}=a_{j}\zpar/\zpar y'_{j}$, $j=1,...,n-r$. So expressing
$dy'_{1},...,dy'_{n-r}$, $dz'_{1},...,dz'_{2m'}$ in terms of
$dx_{1},...,dx_{n}$, $dy_{1},...,dy_{n-r}$, $\, dz_{1},...,dz_{2m'}$ yields
${\tilde\zW}_{1}=\zsu_{j=1}^{n-r}a_{j}dx_{j}\zex dy_{j}+\zW'_{1}$ where $\zW'_{1}$ does
not contain any term involving $dy_{1},...,dy_{n-r}$ and its coefficient functions do not depend
on $(y_{1},...,y_{n-r})$.

Now it is clear that ${\tilde\zlma}$ projects in a partial tensor
field $\zlma:{\mathcal F}\zfl TP$ since the flow of each $\zpar
/\zpar y_{j}$, $j=1,...,n-r$, preserves ${\tilde\zlma}$. For
proving the remainder properties of $\zlma$ consider the product
manifold $M\zpor{\mathbb K}^{r}$ endowed with coordinates
$(x_{1},...,x_{n},y_{1},...,y_{n},z_{1},...,z_{2m'})$, where
$(y_{n-r+1},...,y_{n})$ are the canonical coordinates of ${\mathbb
K}^{r}$, and identify $M$ to $M\zpor \{0\}$. As before forms on
$N$, or on $M$, will be regarded, in the obvious way, as form on
$M\zpor{\mathbb K}^{r}$ when necessary. On this last manifold set
$\zW=\sum_{j=1}^{n}dx_{j}\zex dy_{j}+\sum_{k=1}^{m'}dz_{2k-1}\zex
dz_{2k}$ and $\zW_{1}=\sum_{j=1}^{n-r}a_{j}dx_{j}\zex
dy_{j}+\sum_{j=n-r+1}^{n}adx_{j}\zex dy_{j}+\zW'_{1}$. Then
$\zgran L_{{\frac{\zpar} {\zpar y_{j}}}}\zW= L_{{\frac{\zpar}
{\zpar y_{j}}}}\zW_{1}=0$, $j=1,...,n$.

Let $H$ and $\zW_{k}$ be the $(1,1)$-tensor field and the $2$-form
defined by $\zW_{1}(X,Y)=\zW(HX,Y)$ and
$\zW_{k}(X,Y)=\zW(H^{k}X,Y)$, $k\zpe \mathbb Z$, respectively
(obviously $\zW_{0}=\zW$). Then $H\zpar /\zpar y_{j}=a_{j}\zpar
/\zpar y_{j}$, $j=1,...,n-r$, $H\zpar /\zpar y_{j}=a\zpar /\zpar
y_{j}$, $j=n-r+1,...,n$, and
$\zW_{k}=\sum_{j=1}^{n-r}a_{j}^{k}dx_{j}\zex
dy_{j}+\sum_{j=n-r+1}^{n}a^{k}dx_{j}\zex dy_{j}+\zW'_{k}$ where
$\zW'_{k}$ does not contain any term involving $dy_{1},...,dy_{n}$
and its coefficient functions do not depend on
$y=(y_{1},...,y_{n})$.

As $(i_{X}\zW)\zci
H=\zW(X,H\quad)=\zW(HX,\quad)=\zW_{1}(X,\quad)$, one has
$dx_{j}\zci H=a_{j}dx_{j}$ and $dy_{j}\zci H=a_{j}dy_{j}+\zl_{j}$,
$j=1,...,n-r$; $dx_{j}\zci H=adx_{j}$ and $dy_{j}\zci
H=ady_{j}+\zl_{j}$, $j=n-r+1,...,n$, where $\zl_{1},...,\zl_{n}$
are functional combinations of
$dx_{1},...,dx_{n},dz_{1},...,dz_{2m'}$ whose coefficients do not
depend on $y$. By the same reason each $dz_{k}\zci H$,
$k=1,...,2m'$, is a functional combination of
$dx_{1},...,dx_{n},dz_{1},...,dz_{2m'}$ and its coefficients do
not depend on $y$.

Thus, if $\zp_{1}:M\zpor{\mathbb K}^{r}\zfl M$ is the first
projection, the tensor field $H$ projects in $J$ on $N$ through
$\zp_{N}\zci\zp_{1}$ and in a tensor field $G$ on $P$ through
$\zp_{P}\zci\zp_{1}$. In turns $G$ projects in $J$ via $\zp:P\zfl
N$.

 On the other hand if $\zt= \sum_{j=1}^{n}f_{j}dx_{j}$ then its
$\zW_{1}$-hamiltonian $\sum_{j=1}^{n-r}a_{j}^{-1}f_{j}\zpar/\zpar
y_{j}+\sum_{j=n-r+1}^{n}a^{-1}f_{j}\zpar/\zpar y_{j}$ equals the
$\zW$-hamiltonian of $\zt\zci J^{-1}=
\sum_{j=1}^{n-r}a_{j}^{-1}f_{j}dx_{j}+
\sum_{j=n-r+1}^{n}a^{-1}f_{j}dx_{j}$. Therefore the
$\zW_{1}$-hamiltonians of $\za_{1},...,\za_{r}$, or the
$\zW$-hamiltonians of $\za_{1}\zci J^{-1},...,\za_{r}\zci J^{-1}$,
define a $r$-dimensional foliation ${\mathcal A}_{0}$ on
$M\zpor{\mathbb K}^{r}$ transverse to the first factor, which is
$\zW_{1}$-symplecticly complete since $\za_{1},...,\za_{r}$ are
closed and $\zW$-symplecticly complete because $\za_{1}\zci
J^{-1},...,\za_{r}\zci J^{-1}$ regarded on $N$ define the
foliation ${\bar w}(0)$
(recall that a foliation is called symplecticly complete if
its symplectic orthogonal is a foliation too; see \cite{LP}).

Let $\zp':M\zpor{\mathbb K}^{r}\zfl \zgran {\frac {M\zpor{\mathbb
K}^{r}} {{\mathcal A}_{0}}}$ the canonical projection in the local
quotient. Then $\zp':M\zfl \zgran {\frac {M\zpor{\mathbb K}^{r}}
{{\mathcal A}_{0}}}$ is a diffeomorphism. Moreover as ${\mathcal
A}_{0}$ is bi-symplecticly complete, the Poisson structures
$\zL_{\zW}$ and $\zL_{\zW_{1}}$, associated to $\zW$ and $\zW_1$
respectively, project in two Poisson structures $\zL'$ and
$\zL'_{1}$ on $\zgran {\frac {M\zpor{\mathbb K}^{r}} {{\mathcal
A}_{0}}}$. But the restriction of of $\za_{1}\zci
J^{-1},...,\za_{r}\zci J^{-1}$ and $\zW$ to $M$ defines $\zL$ and
that of $\za_{1},...,\za_{r}$ and $\zW_1$ defines $\zL_1$, so
$\zp':M\zfl\zgran {\frac {M\zpor{\mathbb K}^{r}} {{\mathcal
A}_{0}}}$ transforms $(\zL,\zL_{1})$ in $(\zL',\zL'_{1})$ as a
straightforward algebraic calculation at each point $(p,0)$ shows
[or apply lemma 1.4 of \cite{TUD} to $T_{(p,0)}(M\zpor{\mathbb K}^{r})$,
${\mathcal A}_{0}(p,0)$, $T_{(p,0)}(M\zpor\{0\})$, $\zW(p,0)$ and
$\zW_{1}(p,0)$].

This last construction only needs the properties of $\za_{1},...,\za_{r}$ and $J$ on $N$ but not
the compatibility of $\zL,\zL_{1}$, which may be expressed by means of $\zW,\zW_{1}$. More
exactly:
\bigskip

{\bf Proposition 1.1.} {\it The Poisson structures $\zL,\zL_{1}$ are compatible if and only if
$\za_{1}\zex...\zex\za_{r}\zex d\zW_{2}=0$.}
\bigskip

For proving this proposition we need some auxiliary results.
\bigskip

{\bf Lemma 1.3.} {\it On an even dimensional manifold $\tilde M$
consider a couple of $2$-forms $\zb,\zb_{1}$ such that
$rank\zb=dim\tilde M$ everywhere. Let $K$ the $(1,1)$-tensor field
defined by $\zb_{1}=\zb(K,\quad)$ and set
$\zb_{2}=\zb(K^{2},\quad)$. Then for any vector fields
$X_{1},X_{2},X_{3}$ one has:

\noindent $\zb(N_{K}(X_{1},X_{2}),X_{3})+d\zb(KX_{1},KX_{2},X_{3})=-d\zb_{2}(X_{1},X_{2},X_{3})
+d\zb_{1}(KX_{1},X_{2},X_{3})+d\zb_{1}(X_{1},KX_{2},X_{3})$.}
\bigskip

{\bf Proof.} As the foregoing formula is tensorial, one may assume
$[X_{1},X_{2}]=[X_{1},X_{3}]=[X_{2},X_{3}]=0$ without loss of generality.
Then [recall that $\zb(K,\quad)=\zb(\quad,K)$]:

\noindent $d\zb_{2}(X_{1},X_{2},X_{3})=X_{1}\zb(KX_{2},KX_{3})-X_{2}\zb(KX_{1},KX_{3})
+X_{3}\zb(KX_{1},KX_{2})$

\noindent $d\zb_{1}(KX_{1},X_{2},X_{3})=(KX_{1})\zb(KX_{2},X_{3})-X_{2}\zb(KX_{1},KX_{3})
+X_{3}\zb(KX_{1},KX_{2})-
\zb(K[KX_{1},X_{2}],X_{3}])+\zb([KX_{1},X_{3}],KX_{2})$

\noindent $d\zb_{1}(X_{1},KX_{2},X_{3})=X_{1}\zb(KX_{2},KX_{3})-(KX_{2})\zb(KX_{1},X_{3})
+X_{3}\zb(KX_{1},KX_{2})-
\zb(K[X_{1},KX_{2}],X_{3}])-\zb([KX_{2},X_{3}],KX_{1})$.

Therefore the right side of the formula becomes:

\noindent
$(KX_{1})\zb(KX_{2},X_{3})-(KX_{2})\zb(KX_{1},X_{3})+X_{3}\zb(KX_{1},KX_{2})
-\zb(K[KX_{1},X_{2}]+K[X_{1},KX_{2}],X_{3})
+\zb([KX_{1},X_{3}],KX_{2})-\zb([KX_{2},X_{3}],KX_{1})=
\zb(N_{K}(X_{1},X_{2}),X_{3})+d\zb(KX_{1},KX_{2},X_{3})$.
$\square$
\bigskip

{\bf Corollary 1.3.1.} {\it Assume $\zb,\zb_{1}$ symplectic and
set $\zt=\zb_{1}((K+tI)^{-1},\quad)=\zb((I+tK^{-1})^{-1},\quad)$.
Then

\noindent $d\zt((K+tI)X_{1},(K+tI)X_{2},(K+tI)X_{3})=
t\zb(N_{K}(X_{1},X_{2}),X_{3})=-td\zb_{2}(X_{1},X_{2},X_{3})$.}
\bigskip

{\bf Remark.} At each point $(K+tI)^{-1}$ and $(I+tK^{-1})^{-1}$
are linear combination of powers of $K$, so $\zt$ is a $2$-form on
its domain of definition.

{\bf Proof.} From lemma 1.3 applied to $\zb,\zb_{1}$ and $K$
follows $d\zb_{2}=-\zb(N_{K}(\quad,\quad),\quad)$.

On the other hand, applying this lemma to $\zt$, $\zb_{1}$ and $K+tI$ and taking into account
that $N_{K}=N_{(K+tI)}$ and $\zt((K+tI)^{2},\quad)=\zb_{1}((K+tI),\quad)$ yields:

\noindent $\zt(N_{K}(X_{1},X_{2}),X_{3})+d\zt((K+tI)X_{1},(K+tI)X_{2},X_{3})
=-d\zb_{2}(X_{1},X_{2},X_{3})=\zb(N_{K}(X_{1},X_{2}),X_{3})$.

Hence by replacing $X_{3}$ by $(K+tI)X_{3}$ follows:

\noindent $d\zt((K+tI)X_{1},(K+tI)X_{2},(K+tI)X_{3})=
\zb(N_{K}(X_{1},X_{2}),(K+tI)X_{3})-\zt(N_{K}(X_{1},X_{2}),(K+tI)X_{3})=
\zb(N_{K}(X_{1},X_{2}),(K+tI-(I+tK^{-1})^{-1}(K+tI))X_{3})=
t\zb(N_{K}(X_{1},X_{2}),X_{3})=-td\zb_{2}(X_{1},X_{2},X_{3})$. $\square$

{\it Let us prove proposition 1.1.} Locally always there exists
$t\znoi 0$ such that $I+tH^{-1}$ is invertible. Since
$\zL_{\zW}+t\zL_{\zW_{1}}$ is the dual bivector of
$\zW((I+tH^{-1})^{-1},\quad)$, it projects in $\zL'+t\zL'_{1}$ and
$\zp':M\zfl \zgran {\frac {M\zpor{\mathbb K}^{r}} {{\mathcal
A}_{0}}}$ transforms $\zL+t\zL_{1}$ in $\zL'+t\zL'_{1}$, the
bivector $\zL+t\zL_{1}$ is given by the restriction to $M$ [always
identified to $M\zpor \{0\}$] of $\zW((I+tH^{-1})^{-1},\quad)$ and
$\za_{1}\zci (H+tI)^{-1},...,\za_{1}\zci (H+tI)^{-r}$. Indeed if
$\zW_{1}(Y_{j},\quad)=\za_{j}$, $j=1,...,r$, then
$Y_{1},...,Y_{r}$ span ${\mathcal A}_{0}$ and
$\zW((I+tH^{-1})^{-1}Y_{j},\quad)=\zW_{1}((H+tI)^{-1}Y_{j},\quad)=\za_{j}\zci
(H+tI)^{-1}$, $j=1,...,r$.

On the other hand each $\za_{j}\zci (H+tI)^{-1}$ is the pull-back
of $\za_{j}\zci (J+tI)^{-1}$ and $\za_{1}\zci
(J+tI)^{-1},...,\za_{r}\zci (J+tI)^{-1}$ define the foliation
${\bar w}(t)$. So $\za_{1}\zci (H+tI)^{-1},...,\za_{r}\zci
(H+tI)^{-1}$ define a foliation on $M\zpor{\mathbb K}^{r}$ and, by
restriction, on $M$. Thus $\zL+t\zL_{1}$ is  a Poisson structure,
that is $(\zL,\zL_{1})$  bihamiltonian, if and only if
$\zW((I+tH^{-1})^{-1},\quad)$ is closed modulo
$dy_{n-r+1},...,dy_{n},\za_{1}\zci (H+tI)^{-1},...,\za_{r}\zci
(H+tI)^{-1}$ when $y_{n-r+1}=...=y_{n}=0$.

But the coefficients of $H$ do not depend on $y$ and
$(I+tH^{-1})^{-1}\zpar /\zpar y_{j}$ equals
$(1+ta_{j}^{-1})^{-1}\zpar /\zpar y_{j}$ if $j\zmei n-r$ and
$(1+ta^{-1})^{-1}\zpar /\zpar y_{j}$ if $j\zmai n-r+1$, so
$\zW((I+tH^{-1})^{-1},\quad)=\sum_{j=1}^{n-r}(1+ta_{j}^{-1})^{-1}dx_{j}\zex
dy_{j}+\sum_{j=n-r+1}^{n}(1+ta^{-1})^{-1}dx_{j}\zex dy_{j}
+\zW''_{t}$ where $\zW''_{t}$ do not contain any term involving
$dy_{1},...,dy_{n}$ and its coefficients do not depend on $y$.

Therefore, since $d\zW((I+tH^{-1})^{-1},\quad)=d\zW''_{t}$, the
pair $(\zL,\zL_{1})$ is  bihamiltonian if and only if
$(\za_{1}\zci (H+tI)^{-1})\zex...\zex(\za_{r}\zci (H+tI)^{-1})\zex
d\zW((I+tH^{-1})^{-1},\quad)=0$.

From corollary 1.3.1 applied to $\zW$, $\zW_{1}$ and $H$ follows
$d\zW((I+tH^{-1})^{-1},\quad)((H+tI)\quad,(H+tI)\quad,(H+tI)\quad)
=-td\zW_{2}$, so the above condition holds if and only if
$\za_{1}\zex...\zex\za_{r}\zex d\zW_{2}=0$, {\it which finishes
the proof of proposition 1.1}.
\bigskip

{\bf Lemma 1.4.} {\it If $\zL,\zL_{1}$ are compatible then
$\za_{1}\zex...\zex\za_{r}\zex N_{G}=0$.}
\bigskip

{\bf Proof.} Since $N_{G}$ is the projection of $N_{H}$ it suffices to show that
$N_{H}(X_{1},X_{2})$ is a functional combination of $\zpar /\zpar y_{1},...,
\zpar /\zpar y_{n}$ if $X_{1},X_{2}\zpe Ker(\za_{1}\zex...\zex\za_{r})$.
By proposition 1.1, $d\zW_{2}=\zsu_{j=1}^{r}\zl_{j}\zex\za_{j}$ so
$d\zW_{2}(X_{1},X_{2},\quad)$ is a functional combination of $\za_{1},...,\za_{r}$.

From lemma 1.3 applied to $\zW$, $\zW_{1}$ and $H$ follows
$\zW(N_{H}(X_{1},X_{2}),\quad)=-d\zW_{2}(X_{1},X_{2},\quad)$,
which implies that $\zW(N_{H}(X_{1},X_{2}),\quad)$ is a functional
combination of $\za_{1},...\za_{r}$. Therefore
$N_{H}(X_{1},X_{2})$ has to be a functional combination of $\zpar
/\zpar y_{1},...,\zpar /\zpar y_{n}$. $\square$
\bigskip

{\bf Lemma 1.5.} {\it Assume $\zL,\zL_{1}$ compatible. Then $G$ is
a prolongation of $\zlma$.}
\bigskip

{\bf Proof.} As before $M$ is identified to $M\zpor\{0\}\zco M\zpor{\mathbb K}^{r}$ and
$\zp_{P}\zci\zp_{1}=\zp_{P}$ on $M$. First note that
$H(Im\zL_{1})\zco Im\zL$ modulo $\zpar /\zpar y_{n-r+1},...,\zpar /\zpar y_{n}$ since
on $M$ forms $\za_{1},...,\za_{r}$ define $Im\zL_{1}$ and forms
$\za_{1}\zci J^{-1},...,\za_{r}\zci J^{-1}$ define $Im\zL$. On the other hand if $Y$
belongs to $Im\zL$ modulo $\zpar /\zpar y_{1},...,\zpar /\zpar y_{n}$ then
$\zL((i_{Y}\zW)_{\zbv TM},\quad)$ equals $-Y$ modulo $\zpar /\zpar y_{1},...,
\zpar /\zpar y_{n}$.

Now consider $X\zpe Im\zL_{1}$; then $\zL_{1}((i_{X}\zW_{1})_{\zbv
TM},\quad)=-X$. But $i_{X}\zW_{1}=i_{HX}\zW$ and $HX$ belongs to
$Im\zL$ modulo $\zpar /\zpar y_{1},..., \zpar /\zpar y_{n}$, so
$\zL((i_{X}\zW_{1})_{\zbv TM},\quad)$ equals $-HX$ modulo $\zpar
/\zpar y_{1},..., \zpar /\zpar y_{n}$; that is to say
${\tilde\zlma}([X])=[(\zp_{1})_{*}HX]$. Finally projecting on $P$
via $\zp_{P}$ yields
$\zlma({(\zp_{P})}_{*}X)={(\zp_{P})}_{*}{((\zp_{1})}_{*}HX)
={(\zp_{P}\zci\zp_{1})}_{*}(HX)
=G({(\zp_{P}\zci\zp_{1})}_{*}X)=G({(\zp_{P})}_{*}X)$. $\square$

Let us see that $({\mathcal F},\zlma)$ is a weak Veronese flag.
Since $G$ projects in $J$, the morphism $\zlma:{\mathcal F}\zfl
TP$ projects in ${\bar\zlma}=J_{\zbv{\bar w}(\zinf)}$ and
$\zlma{\mathcal A}\zco{\mathcal A}$. As ${\bar w}$ is a Veronese
web there is no $\bar\zlma$-invariant vector subspace of positive
dimension and ${\bar\zlma}^{*}{\bar\za}$ is closed on ${\bar
w}(\zinf)$ for any closed $1$-form ${\bar\za}$ such that
$Ker{\bar\za}\zcco{\bar w}(\zinf)$. Pulling-back via $\zp:P\zfl N$
shows that conditions 1) and 3) hold. Finally as ${\mathcal
F}=Ker(\za_{1}\zex...\zex\za_{r})$ on $P$, lemmas 1.4 and 1.5
imply $N_{\zlma}=0$.

When we pointed out the existence, at each, point of an algebraic
projection of $\tilde\zlma$ it was showed that
$\zw_{1}=\zw(\zlma,\quad)$ [more exactly that
$\zw_{1}=\zw(\zf,\quad)$]. Therefore $({\mathcal
F},\zlma,\zw,\zw_{1})$ will be a Veronese flag if condition 3) of
this second definition holds.
\bigskip

{\bf Lemma 1.6.} {\it On an open set $P'$ of $P$ consider
functions $f,f_{1},...,f_{k}$, $k\zmai 0$, such that
$\za_{1}\zex...\zex\za_{r}\zex...\zex df_{1}\zex...\zex df_{k}$
has no zero. Assume closed $\zlma^{*}df$ along the foliation
$Ker(\za_{1}\zex...\zex\za_{r}\zex...\zex df_{1}\zex...\zex
df_{k})$. Then
$(L_{X_{f}}\zlma)(Ker(\za_{1}\zex...\zex\za_{r}\zex...\zex
df_{1}\zex...\zex df_{k}))(q)$ is contained in the vector subspace
of $T_{q}P$ spanned by $X_{f_{1}}(q),...,X_{f_{k}}(q)$ whenever
$q\zpe P'$.}
\bigskip

{\bf Proof.} Let ${\widetilde X}_{f},{\widetilde
X}_{f_{1}},...,{\widetilde X}_{f_{k}}$ be the $\zW$-hamiltonians
of $f,f_{1},...,f_{k}$ regarded as functions on an open set of
$M\zpor{\mathbb K}^{r}$. A straightforward calculation shows that
${\widetilde X}_{f},{\widetilde X}_{f_{1}},...,{\widetilde
X}_{f_{k}}$ project in $X_{f},X_{f_{1}},...,X_{f_{k}}$; in
particular $L_{{\tilde X}_{f}}H$ projects in $L_{X_{f}}G$.

On $M\zpor{\mathbb K}^{r}$ one has $\za_{1}\zex...\zex\za_{r}\zex
df_{1}\zex...\zex df_{k}\zex d(df\zci H)=0$ since
$\za_{1}\zex...\zex\za_{r}\zex df_{1}\zex...\zex df_{k}\zex
d(df\zci G)=0$ on $P$. Therefore $\za_{1}\zex...\zex\za_{r}\zex
df_{1}\zex...\zex df_{k}\zex L_{{\tilde X}_{f}}\zW_{1}=0$ whence
locally $L_{{\tilde
X}_{f}}\zW_{1}=\zsu_{j=1}^{r}\zl_{j}\zex\za_{j}
+\zsu_{i=1}^{k}\zm_{i}\zex df_{i}$. But $L_{{\tilde
X}_{f}}\zW_{1}=L_{{\tilde X}_{f}}(\zW(H,\quad)) =\zW(L_{{\tilde
X}_{f}}H,\quad)$ so $\zW(L_{{\tilde
X}_{f}}H,\quad)=\zsu_{j=1}^{r}\zl_{j}\zex\za_{j}
+\zsu_{i=1}^{k}\zm_{i}\zex df_{i}$;  this implies that $L_{{\tilde
X}_{f}}H=\zsu_{j=1}^{r}X'_{j}\zte\za_{j}+
\zsu_{j=1}^{r}X''_{j}\zte\zl_{j} +\zsu_{i=1}^{k}Y_{i}\zte df_{i}
-\zsu_{i=1}^{k}{\widetilde X}_{f_{i}}\zte \zm_{i}$ where
$X''_{1},...,X''_{r}$ are functional combination of $\zpar /\zpar
y_{1},..., \zpar /\zpar y_{n}$.  Thus the projection on $P$ of
$L_{{\tilde X}_{f}}H$ sends
$Ker(\za_{1}\zex...\zex\za_{r}\zex...\zex df_{1}\zex...\zex
df_{k})(q)$ into the vector subspace of $T_{q}P$ spanned by
$X_{f_{1}}(q),...,X_{f_{k}}(q)$. $\square$

When $k=0$ from lemma 1.6 follows the third condition of the
definition of Veronese flag.
{\it Thus $({\mathcal F},\zlma,\zw,\zw_{1})$is a Veronese
flag.}
\bigskip

{\bf 1.2. The bihamiltonian structure over a $(1,1)$-tensor field and a foliation.}

This second part contains a kind of inverse construction of that
of sub-section 1.1. Here, under some assumptions detailed later on, one will associate a
bihamiltonian structure defined on a quotient of the cotangent bundle to a
$(1,1)$-tensor field and a foliation.

Let $N$ be a $n$-manifold. Recall that on $\zL^{r}T^{*}N$ it is
defined a $r$-form $R$, called the Liouville $r$-form, as follows:
if $v_{1},...,v_{r}\zpe T_{\zm}(\zL^{r}T^{*}N)$ then
$R(v_{1},...,v_{r})=\zm(\zp_{*}v_{1},...,\zp_{*}v_{r})$ where
$\zp:\zL^{r}T^{*}N\zfl N$ is the canonical projection. In turn
$\zW=dR$ will be named the Liouville $(r+1)$-form of
$\zL^{r}T^{*}N$. When $r=1$, that is on the cotangent bundle, the
Liouville forms will be denoted $\zr$ and $\zw$ respectively.

Given a skew-symmetric $(1,r)$-tensor field $H$ on $N$, in other
words a section of $TN\otimes\zL^{r}T^{*}N$, let
$\zf_{H}:T^{*}N\zfl\zL^{r}T^{*}N$ be the morphism of vector
bundles defined by $\zf_{H}(\zt)=\zt\zci H$, that is
$\zf_{H}(\zt)(v_{1},...,v_{r})=\zt(H(v_{1},...,Hv_{r}))$. Set
$\zw_{1}=\zf_{H}^{*}\zW$.
\bigskip

{\bf Lemma 1.7.} {\it On a real or complex vector space $V$ of
dimension $2n$, consider a $2$-form $\za$ of rank $2n$ and a
$(r+1)$-form $\zb$. Then there exists $h\zpe V\zte \zL^{r}V^{*}$
connecting $\za$ and $\zb$, that is to say such that
$\zb(v_{1},...,v_{r+1})=\za(h(v_{1},...,v_{r}),v_{r+1})$,
$v_{1},...,v_{r+1}\zpe V$. Moreover $h$ is unique and
$\za(h(v_{1},...,v_{r}),v_{r+1})$

\noindent $=\za(v_{r},h(v_{1},...,v_{r-1},v_{r+1}))$,
$v_{1},...,v_{r+1}\zpe V$.

Conversely, given a $2$-form $\za$ and $h\zpe V\zte \zL^{r}V^{*}$
such that
$\za(h(v_{1},...,v_{r}),v_{r+1})$ $=\za(v_{r},h(v_{1},...,v_{r-1},v_{r+1}))$,
$v_{1},...,v_{r+1}\zpe V$, then setting
$\zb(v_{1},...,v_{r+1})=\za(h(v_{1},...,v_{r}),v_{r+1})$,
$v_{1},...,v_{r+1}\zpe V$, defines a $(r+1)$-form $\zb$.}
\bigskip

The foregoing lemma gives rise to a skew-symmetric $(1,r)$-tensor
field $H^{*}$ on $T^{*}N$ connecting $\zw$ and
$(-1)^{r+1}\zw_{1}$, which will be called {\it the prolongation of
$H$} (to the cotangent bundle).

Given coordinates $x=(x_{1},...x_{n})$ on $N$ let
$(x,y)=(x_{1},...x_{n},y_{1},...,y_{n})$ be the associated coordinates
on $T^{*}N$. Denote by $m(r)$ the set of all the $r$-multi-index
$K:k_{1}<...<k_{r}$ whereas $dx_{K}$ will mean
$dx_{k_{1}}\zex...\zex dx_{k_{r}}$ (as usual elements of $m(1)$
will be represented by small letters). On the other hand $K(j)$,
where $1\zmei j\zmei r$ and $r\zmai 2$, will be the element of
$m(r-1)$ obtained by deleting the term $k_{j}$ of $K$. Assume that
$H=\zsu_{j\zpe m(1), K\zpe m(r)}h_{jK}(\zpar /\zpar x _{j})\zte
dx_{K}$, then:
\bigskip

$\zgran H^{*}=\sum_{j\zpe m(1),K\zpe
m(r)}h_{jK}\zpizq{\frac{\zpar}{\zpar x_{j}}}\zte dx_{K}
-\zsu_{a=1}^{r}(-1)^{a}{\frac{\zpar}{\zpar y_{k_{a}} }}\zte
dy_{j}\zex dx_{K(a)}\zpder$
\smallskip

\hskip 2truecm $-\zgran\sum_{j\zpe m(1),K\zpe m(r+1)}y_{j} \zpizq
\zsu_{a,b=1}^{r+1}(-1)^{a+b}{\frac{\zpar h_{jK(a)}}{\zpar
x_{k_{a}}}}{\frac{\zpar}{\zpar y_{k_{b}}}}\zte dx_{K(b)}\zpder$
\bigskip

Therefore one has:

\noindent (a) $H^{*}$ projects in $H$.

\noindent (b) Let $\xi$ be the radial vector field on $TN^{*}$ [in
coordinates $\xi=\zsu_{j=1}^{n}y_{j}\zpar /\zpar y_{j}$] then
$L_{\zx}H^{*}=0$.

\noindent (c) If $v_{1},v_{2}$ are vertical vectors then
$H^{*}(v_{1},v_{2},...)=0$.

\noindent (d) Set
$\zl(X_{1},...,X_{r+1})=\zw(H^{*}(X_{1},...,X_{r}),X_{r+1})$,
which defines a $(0,r+1)$-tensor field. Then $\zl$ is a closed
$(r+1)$-form.

These four properties characterize the prolongation of $H$ to
$TN^{*}$. More exactly:
\bigskip

{\bf Proposition 1.2.} {\it If a $(1,r)$-tensor field $H'$ defined
on $TN^{*}$ satisfies (a), (b), (c) and (d), then $H'=H^{*}$.}
\bigskip

{\bf Proof.} The tensor field $H_{1}=H'-H^{*}$ satisfies (b), (c)
and (d), and its projection on $N$ vanishes. So in coordinates $(x,y)$:

$\zgran H_{1}=\zsu_{j,a\zpe m(1),K\zpe
m(r-1)}f_{jaK}{\frac{\zpar}{\zpar y_{j}}}\zte dy_{a}\zex dx_{K}
\zgran +\zsu_{j\zpe m(1),L\zpe m(r)}g_{jL}{\frac{\zpar}{\zpar
y_{j}}}\zte dx_{L}$.

But $L_{\zx}H_{1}=0$ therefore $\zx\zpu f_{jaK}=0$ and $\zx\zpu
g_{jL}=g_{jL}$. In other words, each function $f_{jaK}$ only
depend on $x$ and $g_{jL}(x,0)=0$ for every $j\zpe m(1),L\zpe
m(r)$.

Let $\zl_{1}$ be the closed $(r+1)$-form defined by
$\zl_{1}(X_{1},...,X_{r+1})$

\noindent $=\zw(H_{1}(X_{1},...,X_{r}),X_{r+1})$.
Then if $K$ is the multi-index $k_{1}<...<k_{r-1}$ one has [recall
that $\zw=\zsu_{j=1}^{n}dy_{j}\zex dx_{j}$]:

\noindent $\zl_{1}(\zpar/\zpar y_{a},\zpar/\zpar
x_{k_{1}},...,\zpar/\zpar x_{k_{r-1}},\zpar/\zpar x_{j})(x,0)$

\hskip 1truecm $=\zw(H_{1}(\zpar/\zpar y_{a},\zpar/\zpar
x_{k_{1}},...,\zpar/\zpar x_{k_{r-1}}),\zpar/\zpar
x_{j})(x,0)=f_{jaK}(x)$, whereas

\noindent $\zl_{1}(\zpar/\zpar x_{j},\zpar/\zpar
x_{k_{1}},...,\zpar/\zpar x_{k_{r-1}},\zpar/\zpar y_{a})(x,0)$

\hskip 2truecm $=\zw(H_{1}(\zpar/\zpar x_{j},\zpar/\zpar
x_{k_{1}},...,\zpar/\zpar x_{k_{r-1}}),\zpar/\zpar y_{a})(x,0)=0$.

Therefore $f_{jaK}=0$ and $\zl_{1}=\zsu_{S\zpe m(r+1)}h_{S}dx_{S}$
where each function $h_{S}$ only depend on $x$ since $d\zl_{1}=0$,
which implies $L_{\zx}\zl_{1}=0$. But $L_{\zx}\zw=\zw$ and
$L_{\zx}H_{1}=0$ so $L_{\zx}\zl_{1}=\zl_{1}$. Thus $\zl_{1}$ has
to vanish and $H_{1}=0$. $\square$
\bigskip

{\bf Proposition 1.3.} {\it Given a $(1,1)$-tensor field $H$ on
$N$ then the prolongation of $N_{H}$ equals $N_{H^{*}}$.}
\bigskip

{\bf Proof.} By construction $N_{H^{*}}$ satisfies (a), (b) and
(c) with respect to $N_{H}$. Therefore it suffices to show that
setting $\zl(X_{1},X_{2},X_{3})=\zw(N_{H^{*}}(X_{1},X_{2}),X_{3})$
defines a closed $3$-form, which immediately follows from lemma
1.3 applied to $\zw$, $\zw_{1}$ and $H^{*}$. $\square$

Now suppose that $H$ is an invertible $(1,1)$-tensor field and ${\mathcal G}$
a $r$-codimensional foliation both of them defined on $N$. Assume that:

\noindent 1) $\za\zci H$ is closed on $\mathcal G$ whichever $\za$ is a closed $1$-form
such that $Ker\za\zcco \mathcal G$,

\noindent 2) the restriction of $N_H$ to ${\mathcal G}$ vanishes.

Then $(H+tI){\mathcal G}$, $t\zpe {\mathbb K}$, is a $r$-codimensional foliation
on the open set $A_t$ of all points of $N$ where $H+tI$ is invertible. Indeed, reason
as in the first paragraph after lemma 1.1.

Let ${\mathcal G}_{0}$ be the $\zw$-orthogonal of the foliation
$\zp_{*}^{-1}(H{\mathcal G})=\{ v\zpe T(T^{*}N)\zbv \zp_{*}v\zpe H{\mathcal G}\}$,
which equals the $\zw_{1}$-orthogonal of the foliation
$\zp_{*}^{-1}({\mathcal G})=\{ v\zpe T(T^{*}N)\zbv \zp_{*}v\zpe {\mathcal G}\}$
because $\zw_{1}=\zw(H^{*},\quad)$ and $H^*$ projects in $H$. Note that
${\mathcal G}_{0}$ is a symplecticly complete foliation for $\zw$ and $\zw_1$.
On the other hand the quotient $M$ of $T^{*}N$ by ${\mathcal G}_{0}$ is globally
defined and there is a projection $\zp':M\zfl N$ such that
$\zp'\zci{\tilde\zp}=\zp$, where ${\tilde\zp}:T^{*}N\zfl M$ is the canonical
projection. In fact, $M$ can be regarded as the quotient of $T^{*}N$ by a vector
sub-bundle and $\zp':M\zfl N$ as its quotient vector bundle.

Since ${\mathcal G}_{0}$ is both $\zw$ and $\zw_1$ symplecticly complete, the Poisson
structures $\zL_{\zw}$ and $\zL_{\zw_{1}}$, respectively associated to $\zw$ and $\zw_1$,
project in two Poisson structures $\zL$ and $\zL_1$ on $M$.
\bigskip

{\bf Proposition 1.4.} {\it The pair $(\zL,\zL_{1})$ is a bihamiltonian structure.}
\bigskip

{\bf Proof.} The proof is very similar to that of proposition 1.1. As the question
is local one may suppose ${\mathcal G}$ defined by closed $1$-forms
$\za_{1},...,\za_{r}$; of course we will regard $\za_{1},...,\za_{r}$ as forms on $T^{*}N$
by identifying $\za_j$ and $\zp^{*}\za_{j}$, $j=1,...,r$.
Let $\{Y_{1},...,Y_{r}\}$ the basis of ${\mathcal G}_{0}$ defined by
$\zw_{1}(Y_{j},\quad)=\za_{j}$, $j=1,...,r$. Given a point $p\zpe T^{*}N$
consider a scalar $t\znoi 0$ and a small
transversal $P$ to ${\mathcal G}_{0}$, passing through this point, such that
$I+t{(H^{*})}^{-1}$ is invertible around $p$ [that is $(I+tH^{-1})(\zp(p))$ is
invertible], ${\tilde\zp}(P)$ is an open set of $M$ and ${\tilde\zp}:P\zfl{\tilde\zp}(P)$
a diffeomorphism. It suffices to prove that the bivector $\zL+t\zL_{1}$ is a Poisson
structure. Note that it is the projection of $\zL_{\zw}+t\zL_{\zw_{1}}$ which, in turns, is
the dual bivector of $\zw((I+t{(H^{*})}^{-1})^{-1},\quad)$.
But regarded on $P$ by means
of ${\tilde\zp}:P\zfl{\tilde\zp}(P)$, the bivector $\zL+t\zL_{1}$ is given by the restriction
to this transversal of $\zw((I+t{(H^{*})}^{-1})^{-1},\quad)$ and
$\zw((I+t{(H^{*})}^{-1})^{-1}Y_{j},\quad)=\za_{j}\zci (H^{*}+tI)^{-1}$,
$j=1,...,r$.

On the other hand each $\za_{j}\zci (H^{*}+tI)^{-1}$ is the pull-back of
$\za_{j}\zci (H+tI)^{-1}$, and
$\za_{1}\zci (H+tI)^{-1},...,\za_{r}\zci (H+tI)^{-1}$ define the foliation
$(H+tI){\mathcal G}$. So
$\za_{1}\zci (H^{*}+tI)^{-1},...,\za_{r}\zci (H^{*}+tI)^{-1}$ define a
foliation on $T^{*}N$ and, by restriction, on $P$ as well. Thus $\zL+t\zL_{1}$
is a Poisson structure if and only if $\zw((I+t{(H^{*})}^{-1})^{-1},\quad)$
restricted to $P$ is closed modulo
$\za_{1}\zci (H^{*}+tI)^{-1},...,\za_{r}\zci (H^{*}+tI)^{-1}$. Therefore for
finishing the proof it is enough to show that $\zw((I+t{(H^{*})}^{-1})^{-1},\quad)$
is closed on $T^{*}N$ modulo
$\za_{1}\zci (H^{*}+tI)^{-1},...,\za_{r}\zci (H^{*}+tI)^{-1}$.

From corollary 1.3.1, applied to $\zw$, $\zw_{1}$ and $H^*$, follows that

$d(\zw((I+t{(H^{*})}^{-1})^{-1},\quad))((H^{*}+tI)\quad,(H^{*}+tI)\quad,(H^{*}+tI)\quad)
=-td\zw_{2}$

\noindent where $\zw_{2}=\zw((H^{*})^{2},\quad)$. So
the above condition holds if
$\za_{1}\zex...\zex\za_{r}\zex d\zw_{2}=0$.

By lemma 1.3 applied to $\zw$, $\zw_{1}$ and $H^*$ one has
$d\zw_{2}=-\zw(N_{H^{*}}(\quad,\quad),\quad)$. Therefore $d\zw_{2}={\tilde\zw}_{1}$,
where ${\tilde\zw}_{1}=(\zf_{N_{H}})^{*}\zW$ and $\zW$ is the Liouville
$3$-form of $\zL^{2}T^{*}N$ since, by proposition 1.3, the prolongation of
$N_H$ is $N_{H^{*}}$.

On the other hand $\za_{1}\zex...\zex\za_{r}\zex(\zf_{N_{H}})^{*}R=0$, where $R$
is the Liouville $2$-form of $\zL^{2}T^{*}N$, because
$\za_{1}\zex...\zex\za_{r}\zex N_{H}=0$ [calculate $(\zf_{N_{H}})^{*}R$ on coordinates
$(x,y)$ such that $\za_{1}=dx_{1}$,..., $\za_{r}=dx_{r}$]. Hence
$\za_{1}\zex...\zex\za_{r}\zex{\tilde\zw}_{1}=0$, as $\zW=dR$ and
$\za_{1},...,\za_{r}$ are closed, and finally
$\za_{1}\zex...\zex\za_{r}\zex d\zw_{2}=0$. $\square$

{\bf Examples.} 1) On $N={\mathbb K}^{n}$, $n\zmai 1$, consider the foliation given by the closed
$1$-form $\za=\zsu_{j=1}^{n}dx_{j}$ and the $(1,1)$-tensor field
$H=\zsu_{j=1}^{n}h_{j}(x_{j})(\zpar /\zpar x_{j})\zte dx_{j}$ where the functions
$h_{1},...,h_{n}$ never vanish. Then the associated bihamiltonian structure
$(\zL,\zL_{1})$, defined on $M=T^{*}({\mathbb K}^{n})/{\mathcal G}_{0}$, has a symplectic factor
of positive dimension at a point $p\zpe M$ if and only if ${\tilde h}(\zp'(p))=0$ where
${\tilde h}=\zpr_{1\zmei j<k\zmei n}(h_{j}-h_{k})$.
In other words $(\zL,\zL_{1})$ is Kronecker just on the open set
$({\tilde h}\zci\zp')^{-1}({\mathbb K}-\{0\})$.

\noindent 2) Now on $N={\mathbb R}^{n}-\{0\}$, $n\zmai 1$, consider the foliation
${\mathcal G}$ defined by $\za=\zsu_{j=1}^{n}x_{j}^{a_{j}}dx_{j}$, where
$a_{1},...,a_{r}$ are positive natural numbers, and the $(1,1)$-tensor field
$H=\zsu_{j=1}^{n}j(\zpar /\zpar x_{j})\zte dx_{j}$. Then the associated bihamiltonian structure
$(\zL,\zL_{1})$, defined on $M=T^{*}({\mathbb R}^{n}-\{0\})/{\mathcal G}_{0}$, has non-trivial
symplectic factor on the closed set $(h\zci\zp')^{-1}(0)$, where $h=x_{1}\zpu\zpu\zpu x_{n}$, and
is Kronecker on the open set $(h\zci\zp')^{-1}({\mathbb R}-\{0\})$.

Let $\zfi_{t}$ be the flow of the vector field
$\zx=\zsu_{j=1}^{n}(a_{j}+1)^{-1}x_{j}\zpar/\zpar x_{j}$.
As $L_{\zx}\za=\za$ and $L_{\zx}H=0$, the foliation ${\mathcal G}$ and the $(1,1)$-tensor
field  $H$ project in a foliation ${\tilde\mathcal G}$ and a $(1,1)$-tensor
field  $\tilde H$ respectively, defined on the quotient manifold
${\tilde N}=({\mathbb R}^{n}-\{0\})/G$ where  $G=\{\zfi_{k}\zbv k\zpe{\mathbb Z}\}$.
Obviously ${\tilde\mathcal G}$ and $\tilde H$ satisfy 1) and 2), which gives rise
to a bihamiltonian structure on ${\tilde M}=(T^{*}{\tilde N})/{\tilde\mathcal G}_{0}$.
Moreover ${\tilde N}$ is diffeomorphic to $S^{1}\zpor S^{n-1}$.
\bigskip

{\bf 2. Some properties of Veronese flags}

The aim of this section is to establish two results on Veronese flags useful later on. Given a
vector bundle $E$ over a manifold $P$ and a morphism $H:E\zfl E$, we will say that
$H$ is {\it $0$-deformable} if for any points $p,q\zpe P$ there exists an isomorphism between their fibers
$\zf:E(p)\zfl E(q)$ such that $H(p)=\zf^{-1}\zci H(q)\zci\zf$.

By technical reasons parameters are needed. Therefore consider a
foliation ${\cal F}_{1}$ on a manifold $P$, a second foliation
${\cal F}\zco{\cal F}_{1}$ and a morphism $\zlma:{\cal F}\zfl{\cal
F}_{1}$, such that $({\cal F},\zlma)$ is a weak Veronese flag
along ${\cal F}_{1}$; set $r=dim{\cal F}_{1}-dim{\cal F}$.
Let ${\cal A}$ be the foliation of the
largest $\zlma$-invariant vector subspaces (as in section 1) and
$\zp:P\zfl N$ a local quotient of $P$ by ${\cal A}$. Then $N$ is
endowed with the quotient foliations ${\cal F'}_{1}={\cal
F}_{1}/{\cal A}$ and ${\cal F'}={\cal F}/{\cal A}$. Unless another
thing is stated, the Lie and the exterior derivatives of tensor
fields defined along a foliation, for example ${\cal F}_{1}$ on
$P$ or ${\cal F'}_{1}$ on $N$, will be considered along this
foliation. By definition (a system of) coordinates along a
$m$-dimensional foliation ${\cal G}$ will mean a family of
functions $y_{1},...,y_{m}$, on an open set of the support manifold, such that
$dy_{1}\zex,...,\zex dy_{m}$ is a volume form along ${\cal G}$; in this case
$\{\zpar/\zpar y_{1},...,\zpar/\zpar y_{m}\}$ will be the dual basis of
$\{dy_{1},...,dy_{m}\}$.

Consider functions $x_{1},...x_{n}$ on $N$, such that $dx_{1}\zex...\zex dx_{n}$ is a volume form on
${\cal F'}_{1}$, and functions $a_{1},...,a_{n}$ constant along ${\cal F'}_{1}$. Set
$J=\zsu_{j=1}^{n}a_{j}(\zpar/\zpar x_{j})\zte dx_{j}$ where $\{\zpar/\zpar x_{1},...,\zpar/\zpar x_{n}\}$ is
the dual basis of $\{dx_{1},...,dx_{n}\}$. One has:
\bigskip

{\bf Proposition 2.1.} {\it Let $G:{\cal F}_{1}\zfl{\cal F}_{1}$ be a morphism which extends $\zlma$ and
projects in $J$. Assume that:

\noindent  (a) $\zlma\zbv_{\mathcal A}$ is $0$-deformable, nilpotent and flat on each leaf of ${\mathcal A}$,

\noindent (b) $a_{1},...,a_{n}$ never vanish,

\noindent then around every point of $P$ there exists a morphism $G':{\cal F}_{1}\zfl{\cal F}_{1}$, which extends
$\zlma$ and projects in $J$, such that $N_{G'}=0$.}
\bigskip

From proposition 2.1 follows:
\bigskip

{\bf Lemma 2.1.} {\it Consider a morphism ${\widetilde H}:{\widetilde {\cal F}} \zfl {\widetilde {\cal F}}$ where ${\widetilde {\cal F}}$
is a $m$-dimensional foliation on a manifold ${\widetilde P}$. Suppose that ${\widetilde H}$ is $0$-deformable and only has
one eigenvalue. If ${\widetilde H}$ is flat on each leaf of ${\widetilde {\cal F}}$ then, around every point of $\widetilde P$,
there exists a system of coordinates $(z_{1},...,z_{m})$ along
${\widetilde {\cal F}}$such that ${\widetilde H}=\zsu_{j,k=1}^{m}a_{jk}(\zpar/\zpar z_{j})\zte dz_{k}$ where
$a_{jk}\zpe {\mathbb K}$.}
\bigskip

{\bf Proof.} Assume $m<dim\tilde P$ otherwise the result is obvious.
Consider coordinates $(x,y)=(x_{1},...,x_{n},y_{1},...,y_{m})$ defined on an open set $B$ around
a point of ${\widetilde P}$, such that $dx_{1}=...=dx_{n}=0$ gives ${\widetilde {\cal F}}$. Let $a$ be the eigenvalue of
${\widetilde H}$; by taking ${\widetilde H}-aI$ instead ${\widetilde H}$ we may suppose $a=0$. Set ${\widetilde \zlma}={\widetilde H}$
and $J=\zsu_{j=1}^{n}a_{j}(\zpar/\zpar x_{j})\zte dx_{j}$ where $a_{1},...,a_{n}\zpe {\mathbb K}-\{0\}$. By means of
coordinates $(x,y)$, $J$ and ${\widetilde H}$ can be regarded too as tensor fields on $B$ in an obvious way.
Set $G=J+{\widetilde H}$. It easily seem that $({\widetilde {\cal F}},{\widetilde \zlma})$ is a weak Veronese flag on $B$
for which ${\widetilde {\cal A}}={\widetilde {\cal F}}$ and the projected Veronese web is defined by $J$ and
$dx_{1},...,dx_{n}$.

Let $G'$ be the $(1,1)$-tensor field given by proposition 2.1. The characteristic polynomial of both $G$ and $G'$
equals $(\zpr_{j=1}^{n}(t-a_{j}))t^{m}$; even more
$Im(\zpr_{j=1}^{n}(G'-a_{j}I))=Im(\zpr_{j=1}^{n}(G-a_{j}I))={\widetilde {\cal F}}$ [here product means composition].
On the other hand, as $N_{G'}=0$ and $\zpr_{j=1}^{n}(t-a_{j})$ and $t^m$ are relatively prime, locally $B$
splits into a product following the foliations ${\widetilde {\cal F}}=Im(\zpr_{j=1}^{n}(G'-a_{j}I))=Ker((G')^{m})$
and $Im((G')^{m})=Ker(\zpr_{j=1}^{n}(G'-a_{j}I))$. Thus one may consider coordinates
$(x,u)=(x_{1},...,x_{n},u_{1},...,u_{m})$ such that ${\widetilde {\cal F}}$ is given by $dx_{1}=...=dx_{n}=0$ and
$Im((G')^{m})$ by $du_{1}=...=du_{m}=0$ respectively. Moreover
$G'=J+\zsu_{j,k=1}^{m}f_{jk}(u)(\zpar/\zpar u_{j})\zte du_{k}$ since $N_{G'}=0$. But ${\widetilde H}$ is flat on the
leaves of ${\widetilde {\cal F}}$ and $G'_{\zbv{\widetilde {\cal F}}}={\widetilde H}$, so
$\zsu_{j,k=1}^{m}f_{jk}(u)(\zpar/\zpar u_{j})\zte du_{k}$ is flat and one can choose functions $z_{1},...,z_{m}$
of $u$ such that

\centerline{$\zsu_{j,k=1}^{m}f_{jk}(u)(\zpar/\zpar u_{j})\zte du_{k}=\zsu_{j,k=1}^{m}a_{jk}(\zpar/\zpar z_{j})\zte dz_{k}$,
$a_{jk}\zpe{\mathbb K}$. $\square$}
\bigskip

{\bf Lemma 2.2.} {\it Consider a $m$-dimensional foliation
${\widetilde {\cal F}}$ on a manifold ${\widetilde P}$ and a
morphism ${\widetilde H}:{\widetilde {\cal F}} \zfl {\widetilde
{\cal F}}$. Suppose that ${\widetilde H}$ is $0$-deformable and
$N_{\widetilde H}=0$. Then along ${\widetilde {\cal F}}$, given a
function $f$ such that $Kerdf\zcco Ker{\widetilde H}$ and
$d(df\zci {\widetilde H})=0$, locally there exists a function $g$
such that $dg\zci {\widetilde H}=df$.}
\bigskip

{\bf Proof.} As $N_{\widetilde H}=0$ and ${\widetilde H}$ is $0$-deformable, $Im{\widetilde H}$ is a foliation contained in
${\widetilde {\cal F}}$; moreover there exists a vector sub-bundle $E$ of ${\widetilde {\cal F}}$ and a morphism
$\zr:{\widetilde {\cal F}}\zfl{\widetilde {\cal F}}$ such that ${\widetilde {\cal F}}=E\zdi Ker {\widetilde H}$ and
$(\zr\zci {\widetilde H})_{\zbv E}=I$. Set $\za=df\zci\zr$; then $\za\zci {\widetilde H}=df$. From lemma 1.1 applied along
${\widetilde {\cal F}}$ follows that  $d\za(Im{\widetilde H},Im{\widetilde H})=0$, that is $\za_{\zbv Im{\widetilde H}}$ is closed.
Therefore locally there is a function $g$ such that $(dg-\za)_{\zbv Im{\widetilde H}}=0$ so
$dg\zci{\widetilde H}=\za\zci{\widetilde H}=df$. $\square$

{\it One will prove proposition 2.1 by induction on $m=dim {\cal A}$.}
If $m=0$ the result is obvious; now suppose the proposition true up to
dimension $m-1$. Note that in this case lemma 2.1 is also true if
$dim{\widetilde {\cal F}}\zmei m-1$. As the problem is local we
may assume that ${\cal F}'$ is defined by $r$ closed $1$-forms
$\za_{1},...,\za_{r}$ along ${\cal F}'_{1}$, that is
$J,\za_{1},...,\za_{r}$ describe the associated Veronese web.
Functions $x_{1},...,x_{n}$ and forms $\za_{1},...,\za_{r}$ can be
regarded as defined on $P$ in the obvious way (via $\zp$). This
allows us to consider coordinates
$(x,z)=(x_{1},...,x_{n},z_{1},...,z_{m})$ along ${\cal F}_{1}$
such that $dx_{1}=...=dx_{n}=0$ defines ${\cal A}$ and, by means
of $(x,z)$, regard $J$ and $H=\zlma_{\zbv\cal A}$ as
$(1,1)$-tensor fields along ${\cal F}_{1}$. Moreover as
$KerG=Ker(H_{\zbv{\cal A}})\zco{\cal A}$ is a foliation since
$H_{\zbv{\cal A}}$ is flat, coordinates $(x,z)$ can be chosen in
such a way that $KerG$ is defined by
$dx_{1}=...=dx_{n}=dz_{1}=...=dz_{m-s}=0$ where $s=dimKerG$. Then
$G=J+H+\zsu_{j=1}^{m}(\zpar/\zpar z_{j})\zte\zb_{j}$ where every
$\zb_{j}$ is a functional combination of $dx_{1},...,dx_{n}$ and
$H=\zsu_{j=1}^{m}\zsu_{k=1}^{m-s}f_{jk}(\zpar/\zpar z_{j})\zte
dz_{k}$.

But when $i=m-s+1,...,m$ one has:

$\,$

\noindent $\zgran -N_{G}\zpizq{\frac{\zpar}{\zpar
z_{i}}},\quad\zpder= G\zci L_{\frac{\zpar}{\zpar z_{i}}}G-
L_{G\zpizq{\frac{\zpar}{\zpar z_{i}}}\zpder}G=G\zci
L_{\frac{\zpar}{\zpar z_{i}}}G$

$\,$

\hskip 2truecm $\zgran =\zsu_{j=1}^{m}\zsu_{k=1}^{m-s}{\frac
{\zpar f_{jk}} {\zpar z_{i} }} H\zpizq{\frac{\zpar}{\zpar
z_{j}}}\zpder \zte dz_{k}
+\zsu_{j=1}^{m}H\zpizq{\frac{\zpar}{\zpar z_{j}}}\zpder \zte
{\frac {\zpar\zb_{j}} {\zpar z_{i} }}$

$\,$

\noindent therefore $\zpar f_{jk}/\zpar z_{i} =0$,
$j,k=1,...,m-s$, and
$\za_{1}\zex...\zex\za_{r}\zex(\zpar\zb_{j}/\zpar z_{i})=0$,
$j=1,...,m-s$, since $\za_{1}\zex...\zex\za_{r}\zex N_{G}=0$.
Observe that
it is the same proving proposition 2.1 for $G$ or for
$G+\zsu_{j=1}^{r}X_{j}\zte\za_{j}$ where $X_{1},...,X_{r}$ are
vector fields tangent to ${\cal A}$. So, by choosing suitable
vector fields $X_{1},...,X_{r}$, one may suppose
$\zpar\zb_{j}/\zpar z_{i}=0$, $j=1,...,m-s$, without loss of
generality.

In this case $Im(L_{(\zpar /\zpar z_{i})}G)\zco KerG$, $i=m-s+1$,
which allows us to project $G$ in a $(1,1)$-tensor field ${\bar
G}$ defined on the local quotient ${\bar P}$ of $P$ by $KerG$.
Besides ${\cal F}_{1}$, ${\cal F}$, ${\cal A}$ and $\zlma$ project
in similar objects ${\bar {\cal F}}_{1}$, ${\bar {\cal F}}$,
${\bar {\cal A}}$ and $\bar\zlma$ on ${\bar P}$,
$(x,z_{1},...,z_{m-s})$ can be regarded as coordinates along
${\bar {\cal F}}_{1}$, and $N$ is still the quotient of ${\bar P}$
by ${\bar {\cal F}}$; in particular $\za_{1},...,\za_{r}$ may be
seem as forms on ${\bar P}$. Obviously all these objects satisfy
the hypothesis of proposition 2.1 and, by the induction hypothesis,
there exists ${\bar G}':{\bar {\cal F}}_{1}\zfl{\bar {\cal
F}}_{1}$, which extends $\bar\zlma$ and projects in $J$, such that
$N_{{\bar G}'}=0$. Since ${\bar G}'-{\bar G}=\zsu_{j=1}^{r}{\bar
X}_{j}\zte\za_{j}$ where ${\bar X}_{1},...,{\bar X}_{r}$ are
tangent to ${\bar {\cal A}}$ by considering
$G+\zsu_{j=1}^{r}X_{j}\zte\za_{j}$ instead of $G$, where
$X_{1},...,X_{r}$ are tangent to $\cal A$ and project in ${\bar
X}_{1},...,{\bar X}_{r}$, one may suppose $N_{\bar G}=0$ without
loss of generality.

The characteristic polynomial of ${\bar G}$ equals
$(\zpr_{j=1}^{n}(t-a_{j}))t^{m-s}$ since $=dim{\bar {\cal
A}}={m-s}$. As $\zpr_{j=1}^{n}(t-a_{j})$ and $t^{m-s}$ are
relatively prime and $N_{\bar G}=0$, locally ${\bar {\cal F}}_{1}$
splits into a product of two foliations ${\bar{\cal
A}}=Im(\zpr_{j=1}^{n}({\bar G}-a_{j}I))=Ker({\bar G}^{m-s})$ and
${\cal G}=Im({\bar G}^{m-s})=Ker(\zpr_{j=1}^{n}({\bar
G}-a_{j}I))$. Thus we may consider coordinates
$(v,x,u)=(v_{1},...,v_{b},x_{1},...,x_{n},u_{1},...,u_{m-s})$ on
$\bar P$ such that ${\bar {\cal F}}_{1}$ is defined by
$dv_{1}=...=dv_{b}=0$, ${\bar{\cal A}}$ by
$dv_{1}=...=dv_{b}=dx_{1}=...=dx_{n}=0$, and ${\cal G}$ by
$dv_{1}=...=dv_{b}=du_{1}=...=du_{m-s}=0$; moreover

\centerline{${\bar G}=\zsu_{j=1}^{n}a_{j}(\zpar /\zpar x_{j})\zte
dx_{j} +\zsu_{j,k=1}^{m-s}f_{jk}(v,u)(\zpar /\zpar u_{j})\zte
du_{k}$.}

Now from lemma 2.1, applied to coordinates $(v,u)$ and the
$(1,1)$-tensor field $\zsu_{j,k=1}^{m-s}f_{jk}(v,u)(\zpar /\zpar
u_{j})\zte du_{k}$ on ${\bar{\cal A}}$, follows the existence of
coordinates $(v,{\bar z}_{1},...,{\bar z}_{m-s})$ such that

\centerline{$\zsu_{j,k=1}^{m-s}f_{jk}(v,u)(\zpar /\zpar u_{j})\zte
du_{k} =\zsu_{j,k=1}^{m-s}a_{jk}(\zpar /\zpar {\bar z}_{j})\zte
d{\bar z}_{k}$, $a_{jk}\zpe \mathbb {K}$.}

Thus $dx_{j}\zci{\bar G}=a_{j}dx_{j}$, $j=1,...,n$, and every
$d{\bar z}_{k}\zci{\bar G}$, $k=1,...,m-s$, is a linear
combination with constant coefficients of $d{\bar
z}_{1},...,d{\bar z}_{m-s}$. Consequently if
$x_{1},...,x_{n},{\bar z}_{1},...,{\bar z}_{m-s}$ are regarded as
functions on $P$, since $G$ projects in ${\bar G}$, then
$dx_{j}\zci G=a_{j}dx_{j}$, $j=1,...,n$, and each $dz_{k}\zci{\bar
G}$, $k=1,...,m-s$, is a linear combination with constant
coefficients of $d{\bar z}_{1},...,d{\bar z}_{m-s}$. On the other
hand, as $N_{(H_{\zbv \mathcal{A}})}=0$, $Kerd{\bar z}_{k}\zcco
Ker(H_{\zbv \mathcal{A}})$ and $d(d{\bar z}_{k}\zci H)_{\zbv\mathcal
A}=0$, by lemma 2.2 there exists a function $g_{k}$ such that
$(dg_{k}\zci H)_{\zbv\mathcal A}={(d{\bar z}_{k})}_{\zbv
\mathcal{A}}$.

As $Im((H_{\zbv \mathcal{A}})^{*})$ is the annihilator in
$\mathcal{A}$ of $Ker(H_{\zbv \mathcal{A}})$ and $H_{\zbv
\mathcal{A}}$ is nilpotent and
$0$-deformable, around any point and among
$g_{1},...,g_{m-s}$, we may choose functions ${\bar
z}_{m-s+1},...,{\bar z}_{m-{\bar s}}$, where $s-{\bar s}=
dim(Im(H_{\zbv {\mathcal A}})\zin Ker(H_{\zbv {\mathcal A}}) )$,
such that $Ker(H_{\zbv \mathcal{A}})=(Im(H_{\zbv {\mathcal
A}})\zin Ker(H_{\zbv {\mathcal A}}))\zdi Ker({(d{\bar
z}_{1}\zex...\zex d{\bar z}_{m-{\bar s}})}_{\zbv{\mathcal A}})$.
Now if ${\bar z}_{m-{\bar s}+1},...,{\bar z}_{m}$ are functions
such that $Kerd{\bar z}_{j}\zcco Im(H_{\zbv {\mathcal A}})$,
$j=m-{\bar s}+1,...,m$, and $d{\bar z}_{m-{\bar
s}+1}\zex...\zex{d\bar z}_{m}$ restricted to $Ker({(d{\bar
z}_{1}\zex...\zex d{\bar z}_{m-{\bar s}})}_{\zbv{\mathcal A}})$
does not vanish anywhere, then ${\bar z}=({\bar z}_{1},...,{\bar
z}_{m})$ is a system of coordinates on ${\mathcal A}$ and
$(x,{\bar z})$ a system of coordinates on ${\mathcal F}_{1}$. By
construction $d{\bar z}_{j}\zci G$, $j=m-s+1,...,m$, equals a
linear combination with constant coefficients of $d{\bar
z}_{1},...,d{\bar z}_{m-s}$ plus a functional combination of
$dx_{1},...,dx_{n}$.

In short, naming $z_{k}$ every function ${\bar z}_{k}$ allows us
to suppose that in coordinates $(x,z)$

$G=\zsu_{j=1}^{n}a_{j}(\zpar /\zpar x_{j})\zte
dx_{j}+\zsu_{j=1}^{m}\zsu_{k=1}^{m-s}a_{jk}(\zpar /\zpar
z_{j})\zte dz_{k}$

\hskip 5truecm $+\zsu_{j=m-s+1}^{m}(\zpar /\zpar z_{j})\zte
\zb_{j}$

\noindent where every $a_{jk}\zpe\mathbb{K}$, each $\zb_{j}$ is a
functional combination of $dx_{1},...,dx_{n}$ and $\{\zpar/\zpar
z_{m-s+1},...,\zpar/\zpar z_{m}\}$ a basis of $KerG$.

Besides, by linearly rearranging $z_{1},...,z_{m}$ if necessary,
one may suppose that ${\{\zpar/\zpar z_{\zl}\}}_{\zl\zpe L}$, for
some subset $L$ of $\{1,...,m\}$, is a basis of $G({\mathcal A})$.

But now $N_{G}(\zpar/\zpar z_{k},\quad)=L_{G(\zpar/\zpar z_{k})}G
-G\zci L_{\zpar/\zpar z_{k}}G=L_{G(\zpar/\zpar z_{k})}G$ and
$\za_{1}\zex...\zex\za_{r}\zex N_{G}=0$, therefore
$\za_{1}\zex...\zex\za_{r}\zex(\zpar\zb_{j}/\zpar z_{\zl})=0$,
$\zl\zpe L$, $j=m-s+1,...,m$. Thus considering
$G+\zsu_{j=1}^{r}X_{j}\zte\za_{j}$ instead of $G$, where
$X_{1},...,X_{r}$ are suitable functional combinations of
$\zpar/\zpar z_{m-s+1},...,\zpar/\zpar z_{m}$, and calling it $G$
again allows us to suppose $\zpar\zb_{j}/\zpar z_{\zl}=0$,
$\zl\zpe L$, $j=m-s+1,...,m$ without loss of generality.

By lemma 1.1, $dx_{j}\zci N_{G}=dz_{k}\zci N_{G}=0$, $j=1,...n$,
$k=1,...,m-s$. Therefore one has to study $dz_{j}\zci N_{G}$ when
$j=m-s+1,...m$. Note that each $(\zb_{j}\zci J^{-1})\zci N_{G}=0$
[here $J^{-1}=\zsu_{i=1}^{n}a_{i}^{-1}(\zpar/\zpar x_{i})\zte
dx_{i}$ ], so $dz_{j}\zci N_{G}=(dz_{j}-\zb_{j}\zci J^{-1})\zci
N_{G}$, $j=m-s+1,...m$, and from lemma 1.1 applied to
$dz_{j}-\zb_{j}\zci J^{-1}$ and $G$ follows $(d(dz_{j}-\zb_{j}\zci
J^{-1}))(G,G)+(dz_{j}-\zb_{j}\zci J^{-1})\zci N_{G}=0$ since
$(dz_{j}-\zb_{j}\zci J^{-1})\zci G$ and $(dz_{j}-\zb_{j}\zci
J^{-1})\zci G^{2}$ equal zero or a linear combination with
constant coefficients of $dz_{1},...,dz_{m-s}$.

Hence

$\,$

\noindent $\zgran dz_{j}\zci N_{G}
=(d(\zb_{j}\zci J^{-1}))(G,G)
=(d_{x}(\zb_{j}\zci J^{-1}))(J,J)+\zsu_{k=m-s+1}^{m}\zb_{k}\zex
\zpizq{\frac{\zpar(\zb_{j}\zci J^{-1})}{\zpar z_{k}}}\zci J\zpder
=\zpizq d_{x}(\zb_{j}\zci J^{-1})+\zsu_{k=m-s+1}^{m}(\zb_{k}\zci
J^{-1})\zex {\frac{\zpar(\zb_{j}\zci J^{-1})}{\zpar z_{k}}}\zpder
(J,J)$

$\,$

\noindent where $d_{x}$ is the exterior derivative with respect to
$x=(x_{1},...,x_{n})$ only [recall that  ${\{\zpar/\zpar
z_{\zl}\}}_{\zl\zpe L}$ is a basis of $G({\mathcal A})$ and
$\zpar\zb_{j}/\zpar z_{\zl}=0$].

Therefore the equation $\za_{1}\zex...\zex\za_{r}\zex N_{G}=0$ is
equivalent to the system

$\,$

\noindent (1) \hskip .1truecm $\cases{\zgran\zpizq
d_{x}(\zb_{j}\zci J^{-1})+\zsu_{k=m-s+1}^{m}(\zb_{k}\zci
J^{-1})\zex {\frac{\zpar(\zb_{j}\zci J^{-1})}{\zpar
z_{k}}}\zpder\cr \,\cr \hskip 4.5truecm \zex(\za_{1}\zci
J^{-1})\zex...\zex(\za_{r}\zci J^{-1})=0 \cr\zgran j=m-s+1,...,m
\cr}$

$\,$

By the same reason if $G'=G+\zsu_{j=m-s+1}^{m}(\zpar/\zpar
z_{j})\zte\zb'_{j}$, where $\zb'_{m-s+1},...,\zb'_{m}$ are
functional combinations of $\za_{1},...,\za_{r}$ whose coefficient
functions do not depend on $z_{\zl}$, $\zl\zpe L$, the equation
$N_{G'}=0$ is equivalent to the system

$\,$

\noindent (2) \hskip .1truecm $\cases{\zgran
d_{x}((\zb_{j}+\zb'_{j})\zci
J^{-1})+\zsu_{k=m-s+1}^{m}((\zb_{k}+\zb'_{j})\zci J^{-1})\zex
{\frac{\zpar((\zb_{j}+\zb'_{j})\zci J^{-1})}{\zpar
z_{k}}}=0\cr \,\cr j=m-s+1,...,m \cr}$

$\,$

In other words we need to show that given forms
$\zb_{m-s+1}...,\zb_{m}$ satisfying system (1), there exist forms
$\zb'_{m-s+1}...,\zb'_{m}$ such that system (2) is satisfied too.

On $N$ forms $\za_{1}\zci J^{-1},...,\za_{r}\zci J^{-1}$ define a
foliation contained in ${\mathcal F}'_{1}$ since
$\za_{1},...,\za_{r},J$ give rise to a Veronese web along ${\mathcal
F}'_{1}$; moreover around every point of $N$ there exist indices
$1\zmei k_{1}<...<k_{n-r}\zmei n$ such that $dx_{k_{1}}\zex...\zex
dx_{k_{n-r}}\zex(\za_{1}\zci J^{-1})\zex...\zex (\za_{r}\zci
J^{-1})$ does not vanish anywhere. As the order of functions
$x_{1},...,x_{n}$ is arbitrary, we may assume $dx_{1}\zex...\zex
dx_{n-r}\zex(\za_{1}\zci J^{-1})\zex...\zex (\za_{r}\zci J^{-1})$
non-singular and consider coordinates $y=(y_{1},...,y_{n})$ along
${\mathcal F}'_{1}$ such that $y_{1}=x_{1}$,..., $y_{n-r}=x_{n-r}$
and $Ker(dy_{n-r+1}\zex...\zex dy_{n})=Ker((\za_{1}\zci
J^{-1})\zex...\zex (\za_{r}\zci J^{-1}))$; thus $\za_{1}\zci
J^{-1},...,\za_{r}\zci J^{-1}$ and each $\zb'_{j}\zci J^{-1}$ are
functional combination of $dy_{n-r+1},...,dy_{n}$; in the first
case the coefficients only depend on $y$ and in the second one
they do not depend on $z_{\zl}$, $\zl\zpe L$. Moreover one can
assume that every $\zb_{j}\zci J^{-1}$ is only combination of
$dy_{1},...,dy_{n-r}$; indeed if $\zb_{j}\zci
J^{-1}=\zg_{j}+\zr_{j}$ with $\zg_{j}\zex dy_{1}\zex...\zex
dy_{n-r}=0$ and $\zr_{j}\zex dy_{n-r+1}\zex...\zex dy_{n}=0$, it
suffices replacing $G$ by $G-\zsu_{j=m-s+1}^{m}(\zpar/\zpar
z_{j})\zte(\zr_{j}\zci J)$.

On the other hand, linearly rearranging coordinates $z$ allows us
to suppose that $\{1,...,m-s\}-L=\{1,...,m'\}$ and
$\{m-s+1,...,m\}-L=\{m-s+1,...,m-s+s'\}$ where $m'\zmei m-s$ and
$s'\zmei s$
(here $m'=0$ means $\{1,...,m-s\}\zco L$ and
$s'=0$ that $\{m-s+1,...,m\}\zco L$).
Now on $P$ take a system of coordinates $(v,y,u,w)=
(v_{1},...,v_{a+m'},y_{1},...,y_{n},u_{1},...,u_{s},w_{1},...,w_{m-m'-s})$
such that $dv_{1}=...=dv_{a}=0$ defines ${\mathcal F}_{1}$,
$v_{a+k}=z_{k}$, $k=1,...,m'$, $u_{j}=z_{m-s+j}$, $j=1,...,s$, and
$w_{k}=z_{m'+k}$, $k=1,...,m-m'-s$. Set $\zt_{j}=\zb_{m-s+j}\zci
J^{-1}$, $j=1,...,s$. Since $\zb_{m-s+j}$ and $\zb'_{m-s+j}$ do
not depend on $z_{\zl}$, $\zl\zpe L$, our problem may be stated in
coordinates $(v,y,u)$, that is on a manifold $P'$ of dimension
$a+m'+s$ and along the foliation ${\mathcal G}'$ defined by
$dv_{1}=...=dv_{a+m'}=0$, as follows:

Given $1$-forms $\zt_{j}=\zsu_{k=1}^{n-r}f_{jk}dy_{k}$,
$j=1,...,s$, where functions $f_{jk}$ do not depend on
$u_{s'+1},...,u_{s}$ such that

$\,$

\noindent $(3)\hskip .5truecm \zgran
d_{(y_{1},...,y_{n-r})}\zt_{j}+\zsu_{i=1}^{s}\zt_{i}\zex{\frac{\zpar\zt_{j}}{\zpar
u_{i} } }=0$, $\,j=1,...,s$,

$\,$

\noindent find forms
${\widetilde\zt}_{j}=\zt_{j}+\zsu_{k=n-r+1}^{n}f_{jk}dy_{k}$,
$\,j=1,...,s$, where each $f_{jk}$ does not depend on
$u_{s'+1},...,u_{s}$ such that

$\,$

\noindent $(4)\hskip .5truecm \zgran
d_{y}{\widetilde\zt}_{j}+\zsu_{i=1}^{s}{\widetilde\zt}_{i}\zex
{\frac{\zpar{\widetilde\zt}_{j}}{\zpar u_{i} } }=0$,
$\,j=1,...,s$,

$\,$

\noindent (here $d_{(y_{1},...,y_{n-r})}$ and $d_{y}$ are the
exterior derivative in $(y_{1},...,y_{n-r})$ or
$y=(y_{1},...,y_{n})$ respectively).
\bigskip

{\bf Lemma 2.3.} {\it Forms
${\widetilde\zt}_{1},...,{\widetilde\zt}_{s}$ always exist
locally.}
\bigskip

 {\bf Proof.} As a straightforward calculation
shows, system (3) is equivalent to say that vector fields
$X_{k}=\zpar/\zpar y_{k}+\zsu_{j=1}^{s}f_{jk}\zpar/\zpar u_{j}$,
$k=1,...,n-r$, commute among them.

An analogous statement holds
for system (4).

In turn, functions $f_{jk}$ do not depend on $u_{s'+1},...,u_{s}$
if and only if $X_{1},...,X_{n-r}$ commute with vector fields
$Y_{1},...,Y_{s-s'}$, where each $Y_{i}=\zpar/\zpar u_{s'+i}$.

Since by hypothesis $X_{1},...,X_{n-r},Y_{1},...,Y_{s-s'}$ commute
and are linearly independent everywhere, along ${\mathcal G}'$ and
around every point, there exist coordinates
$\zn_{1},...,\zn_{n+s}$ such that
$\zn_{1}=y_{1},...,\zn_{n}=y_{n}$, $X_{k}=\zpar/\zpar\zn_{k}$,
$k=1,...,n-r$, and  $Y_{i}=\zpar/\zpar\zn_{n+s'+i}$,
$i=1,...,s-s'$.

Set  $X_{k}=\zpar/\zpar\zn_{k}$ when $k=n-r+1,...,n$. Then in
coordinates $(y,u)$ one has $X_{k}=\zpar/\zpar
y_{k}+\zsu_{j=1}^{s}f_{jk}\zpar/\zpar u_{j}$, $k=n-r+1,...,n$.
Moreover by construction $X_{1},...,X_{n},Y_{1},...,Y_{s-s'}$
commute among them, so forms
${\widetilde\zt}_{j}=\zsu_{k=1}^{n}f_{jk}dy_{k}
=\zt_{j}+\zsu_{k=n-r+1}^{n}f_{jk}dy_{k}$, $\,j=1,...,s$, satisfy
system (4) and functions $f_{jk}$ do not depend on
$u_{s'+1},...,u_{s}$. $\square$

{\it Now proposition 2.1 is proved.}

The next step will be extending this result to Veronese flags.
Therefore let $\zw,\zw_{1}$ be a symplectic form and a closed
$2$-form, respectively, defined on ${\mathcal A}$. Suppose that
$({\cal F},\zlma,\zw,\zw_{1})$ is a Veronese flag on $P$
or at some point of $P$ both along ${\cal
F}_{1}$ [in the second case by definition
condition 3') has to hold on neighbourhoods on $P$
of this point].
Set $dim{\mathcal A}=2m$ (now the dimension of ${\mathcal A}$ has
to be even since $\zw$ is symplectic).
\bigskip

{\bf Theorem 2.1.} {\it Let $G:{\cal F}_{1}\zfl{\cal F}_{1}$ be a
morphism which extends $\zlma$ and projects in $J$. Assume that:

\noindent  (a) $\zlma\zbv_{\cal A}$ is $0$-deformable and its
characteristic polynomial equals $(t-a)^{2m}$ where
$a\zpe\mathbb{K}$,

\noindent (b) functions $a_{1},...,a_{n}$ never take the value $a$,

\noindent then around every point of $P$
such that $({\cal F},\zlma,\zw,\zw_{1})$ is a Veronese flag at it
there exist a morphism
$G':{\cal F}_{1}\zfl{\cal F}_{1}$, which extends $\zlma$ and
projects in $J$, and functions $z_{1},...,z_{2m}$ such that
$(x_{1},...,x_{n},z_{1},...,z_{2m})$ is a system of coordinates
along ${\cal F}_{1}$,

\centerline{ $G'=\zsu_{j=1}^{n}a_{j}(\zpar/\zpar x_{j})\zte dx_{j}
+\zsu_{j,k=1}^{2m}a_{jk}(\zpar/\zpar z_{j})\zte dz_{k}$,}

\noindent where every $a_{jk}\zpe\mathbb{K}$, and $\zw,\zw_{1}$
are expressed with constant coefficients relative to
${(dz_{1})}_{\zbv{\cal A}},...,{(dz_{2m})}_{\zbv{\cal A}}$.}
\bigskip

Again, one will prove theorem 2.1 by induction on $m$. If $m=0$
the result is obvious; now suppose the theorem true up to $m-1$.
Note that we may assume $\zlma_{\zbv \mathcal{A}}$ nilpotent by
considering $G-aI$, $\zlma-aI$ and $\zw_{1}-a\zw$ instead of $G$,
$\zlma$ and $\zw_{1}$. Then $KerG=Ker\zw_{1}\zco \mathcal{A}$; so
$KerG$ is a foliation since $\zw_{1}$ is closed. Consider
coordinates $(x,z)=(x_{1},...,x_{n},z_{1},...,z_{2m})$ along
${\mathcal F}_{1}$ such that
$dx_{1}=...=dx_{n}=dz_{1}=...=dz_{2(m-s)}=0$ defines $KerG$, where
$dimKerG=2s$. Reasoning as in the proof of proposition 2.1 allows
us to assume $G$ projectable in a tensor field ${\bar G}$, defined
on the local quotient ${\bar P}$ of $P$ by $KerG$, and consider
the objects ${\bar{\mathcal F}}_{1}$, ${\bar{\mathcal F}}$,
${\bar{\mathcal A}}$ and $\bar\zlma$ with the obvious difference
that now $dim{\bar{\mathcal A}}$ is even. Thus $({\bar{\mathcal
F}},{\bar{\mathcal A}})$ is a weak Veronese flag along
${\bar{\mathcal F}}_{1}$.

On the other hand $\zw_{1}$ projects in a symplectic form
${\bar\zw}$ on ${\bar{\mathcal A}}$. By lemma 1.3 the $2$-form
$\zw_{2}(X,Y)=\zw(\zlma X,Y)$ is closed and $Ker\zw_{2}\zcco
Ker\zw_{1}=KerG$; therefore it projects in a closed $2$-form
${\bar\zw}_{1}$ on ${\bar{\mathcal A}}$ such that
${\bar\zw}_{1}={\bar\zw}({\bar\zlma},\quad)$, and $({\bar{\mathcal
F}},{\bar\zlma},{\bar\zw},{\bar\zw}_{1})$ will be a Veronese flag
if we are able to check the third condition of the definition of
this notion.

Let $h$ be a function on an open set of ${\bar P}$ such that
${\bar\zlma}^{*}dh$ is a closed form on ${\bar{\mathcal F}}$, that
is such that $\za_{1}\zex...\zex\za_{r}\zex (d(dh\zci {\bar
G}))=0$. Regarded on $P$ one has $dh(KerG)=0$ and
$\za_{1}\zex...\zex\za_{r}\zex (d(dh\zci G))=0$. In particular,
locally and along ${\mathcal F}_{1}$, $dh=\zb\zci G$ for some
$1$-form $\zb$ and, by lemma 1.1, one has $d\zb(G,G)+d(dh\zci
G)+\zb\zci N_{G}=0$. Hence, as $\za_{1}\zex...\zex\za_{r}\zex
N_{G}=\za_{1}\zex...\zex\za_{r}\zex (d(dh\zci G))=0$, results
$\za_{1}\zex...\zex\za_{r}\zex d\zb(G,G)=0$, that is $(\za_{1}\zci
J^{-1})\zex...\zex(\za_{r}\zci J^{-1})\zex d\zb=0$.

But $\za_{1}\zci J^{-1},...,\za_{r}\zci J^{-1}$ define a
foliation, therefore $\zb=dg$ modulo $\za_{1}\zci
J^{-1},...,\za_{r}\zci J^{-1}$ for some function $g$. Thus
$\za_{1}\zex...\zex\za_{r}\zex (dg\zci G-dh)=0$, whence
$\za_{1}\zex...\zex\za_{r}\zex (d(dg\zci G))=0$; in other words
$\zlma^{*}dg$ is closed on ${\mathcal F}$.

Let $X$ be the $\zw$-hamiltonian of $g$. From
$\zw_{1}(X,\quad)=\zw(GX,\quad)=-(dg\zci G)_{\zbv{\mathcal
A}}=-dh_{\zbv\mathcal A}$ follows that the projection $\bar X$ of
$X$ on$\bar P$ is the $\bar\zw$-hamiltonian of $h$. But
$L_{X}\zlma=0$ since $({\mathcal F},\zlma,\zw,\zw_{1})$ is a
Veronese flag, so $L_{\bar X}{\bar\zlma}=0$; that is to say
$({\bar{\mathcal F}},{\bar\zlma},{\bar\zw},{\bar\zw}_{1})$ is a
Veronese flag too (everywhere or at some point).

By the induction hypothesis, there exist a morphism ${\bar
G}':{\bar{\mathcal F}}_{1}\zfl{\bar{\mathcal F}}_{1}$ extending
${\bar\zlma}$ and projecting in $J$ and functions
$z_{1},...,z_{2(m-s)}$, such that
$(x_{1},...,x_{n},z_{1},...,z_{2(m-s)})$ is a system of
coordinates along ${\bar{\mathcal F}}_{1}$ in which ${\bar G}'$,
${\bar\zw}$ and ${\bar\zw}_{1}$ are written with constant
coefficients. But ${\bar G}'-{\bar G}=\zsu_{j=1}^{r}{\bar
X}_{j}\zte\za_{j}$ where ${\bar X}_{1},...,{\bar X}_{r}$ are
tangent to ${\bar{\mathcal A}}$. Therefore considering
$G+\zsu_{j=1}^{r}{\bar X}_{j}\zte\za_{j}$ instead of $G$, where
$X_{1},...,X_{r}$  are tangent to ${\mathcal A}$  and project in
${\bar X}_{1},...,{\bar X}_{r}$, allows us to suppose that $G$
projects in ${\bar G}'$; that is to say ${\bar G}'={\bar G}$.

On the other hand, proceeding as in the proof of proposition 2.1
shows the existence of vector fields ${\widetilde
X}_{1},...,{\widetilde X}_{r}$, tangent to $KerG$, such that the
Nijenhuis torsion of $G+\zsu_{j=1}^{r}{\widetilde
X}_{j}\zte\za_{j}$ vanishes; in other words one may assume
$N_{G}=0$. Indeed, see $({\bar{\mathcal F}},{\bar\zlma})$,
$({\mathcal F},{\zlma})$ as weak Veronese flags and ${\bar G}$, $G$
like suitable prolongations of ${\bar\zlma}$, $\zlma$
respectively.

In short, only case to consider: in coordinates
$(x_{1},...,x_{n},z_{1},...,z_{2(m-s)})$ ${\bar G}$, $\bar\zw$,
${\bar\zw_{1}}$ are written with constant coefficients
following theorem 2.1 and $N_{G}=0$.

Regarded like function on $P$ every $dz_{k}\zci G$,
$k=1,...,2(m-s)$, is a linear combination with constant
coefficients of $dz_{1},....dz_{2(m-s)}$; moreover $dx_{j}\zci G
=a_{j}dx_{j}$, $j=1,...,n$. By lemma 2.2 there exist functions
$g_{1},...,g_{2(m-s)}$ such that $dg_{k}\zci G=dz_{k}$,
$k=1,...,2(m-s)$.

Since $G_{\zbv\mathcal A}={\zlma}_{\zbv\mathcal A}$ is
$0$-deformable and nilpotent, around any point and among $g_{1},...,g_{2(m-s)}$,
one can choose functions ${\bar z}_{2(m-s)+1},...,{\bar z}_{2(m-{\bar
s})}$, where $2(s-{\bar s})=dim(G({\mathcal A})\zin KerG)$, such that

\centerline{$KerG=(G({\mathcal A})\zin KerG)\zdi
Ker({(dz_{1}\zex...\zex dz_{2(m-s)}\zex d{\bar
z}_{2(m-s)+1}\zex...\zex d{\bar z}_{2(m-{\bar
s})})}_{\zbv{\mathcal A}})$.}

\noindent Now if ${\bar z}_{2(m-{\bar s})+1},...,{\bar z}_{2m}$
are functions such that $Kerd{\bar z}_{j}\zcco ImG$, $j=2(m-{\bar
s})+1,...,2m$, and $d{\bar z}_{2(m-{\bar s})+1}\zex...\zex d{\bar
z}_{2m}$ restricted to

\centerline{$ Ker({(dz_{1}\zex...\zex dz_{2(m-s)}\zex d{\bar
z}_{2(m-s)+1}\zex...\zex d{\bar z}_{2(m-{\bar
s})})}_{\zbv{\mathcal A}})$}

\noindent does not vanish anywhere, then
$(x_{1},...x_{n},z_{1},...,z_{2(m-s)},{\bar
z}_{2(m-s)+1},...,{\bar z}_{2m})$ is a system of coordinates along
${\mathcal F}_{1}$. By construction $d{\bar z}_{j}\zci G$,
$j=2(m-s)+1,...,2m$, equals a linear combination with constant
coefficients of $dz_{1},...,dz_{2(m-s)}$. Thus

\centerline{$G=\zsu_{j=1}^{n}a_{j}(\zpar/\zpar x_{j})\zte dx_{j}+
\zsu_{j=1}^{2m}\zsu_{k=1}^{2(m-s)}Z_{jk}\zte dz_{k}$}

\noindent where each $Z_{jk}$ is a linear combination with
constant coefficients of

\centerline{$\zpar/\zpar z_{1},...,\zpar/\zpar
z_{2(m-s)},\zpar/\zpar {\bar z}_{2(m-s)+1},...,\zpar/\zpar {\bar
z}_{2m}$.}

Moreover, in these coordinates, $\zw_{1}$ and $\zw_{2}$ are
written with constant coefficients since $\bar\zw$ and
$\bar\zw_{1}$ are in coordinates
$(x_{1},...,x_{n},z_{1},...,z_{2(m-s)})$.

Let $X_{k}$ be the $\zw$-hamiltonian of $z_{k}$, $k=1,...2(m-s)$,
or ${\bar z}_{k}$, $k=2(m-s)+1,...2m$. Then
$\zw_{1}(X_{k},\quad)=\zw(GX_{k},\quad)$ equals $-dz_{k}\zci G$ or
$-d{\bar z}_{k}\zci G$; in both cases a linear combination with
constant coefficients of $dz_{1},...,dz_{2(m-s)}$ because
$\zw_{1}$ projects in $\bar\zw$. In other words

\centerline{$X_{k}=\zsu_{i=2(m-s)+1}^{2m}f_{ki}\zpar/\zpar {\bar
z}_{i} +\zsu_{j=1}^{2(m-s)}b_{kj}\zpar/\zpar z_{j}$}

\noindent where each $b_{kj}\zpe\mathbb{K}$.

In particular $\{z_{j},{ z}_{k}\}_{\zw}$, $j,k=1,...,2(m-s)$,
and $\{z_{j},{\bar z}_{k}\}_{\zw}$, $j=1,...,2(m-s)$,
$k=2(m-s)+1,...,2m$, are constant [here $\{\quad,\quad \}_{\zw}$
and $\zL_{\zw}$ are respectively the Poisson structure and the
dual bivector on $\mathcal A$ associated to $\zw$].

On the other hand, everywhere or close to some point,
$L_{X_{k}}\zlma=0$ since $dz_{k}\zci G$, or
$d{\bar z}_{k}\zci G$, is closed. Hence
$\za_{1}\zex...\zex\za_{r}\zex L_{X_{k}}G=0$. A straightforward
calculation shows that
$L_{X_{k}}G=-\zsu_{i=2(m-s)+1}^{2m}(\zpar/\zpar {\bar z}_{i})\zte
(df_{ki}\zci G)$, so $\za_{1}\zex...\zex\za_{r}\zex
(d_{x}f_{ki}\zci J)=0$ where $d_{x}$ is the exterior derivative
with respect to $x=(x_{1},...,x_{n})$; that is $(\za_{1}\zci
J^{-1})\zex...\zex(\za_{r}\zci J^{-1})\zex d_{x}f_{ki}=0$,
$i=2(m-s)+1,...,2m$. In other words functions $f_{ki}$ are basic
for the foliation ${\mathcal G}''\zco{\mathcal F}_{1}$ defined by
$\za_{1}\zci J^{-1},...,\za_{r}\zci
J^{-1},dz_{1},...,dz_{2(m-s)},d{\bar z}_{2(m-s)+1},...,d{\bar
z}_{2m}$.

But $\{{\bar z}_{k},{\bar z}_{i}\}_{\zw}=f_{ki}$, therefore
$\zL_{\zw}$ and by consequence $\zw$ are written with coefficients
which are ${\mathcal G}''$-basic functions.
\bigskip

{\bf Lemma 2.3.} {\it Along a foliation ${\widetilde{\mathcal F}}$
 of dimension $2{\widetilde m}$
defined on a manifold ${\widetilde P}$, consider
a symplectic form $\zl$
and functions $f_{1},...,f_{k}$ such that $df_{1}\zex...\zex
df_{k}$ has no zeros. Assume constant every function
$\{f_{i},f_{j}\}$, $i,j=1,...,k$. Then locally there are functions
$g_{1},...,g_{2{\widetilde m}-k}$ such that
$(f_{1},...,f_{k},g_{1},...,g_{2{\widetilde m}-k})$ is a system of
coordinates along ${\widetilde{\mathcal F}}$ and $\zl$ is written
with constant coefficients relative to it.}
\bigskip

{\bf Proof.} It is just  one of the version of Darboux theorem.
$\square$

Pulling-back the functions given by lemma 2.3, applied to the
projections on the local quotient $P/{\mathcal G}''$ of ${\mathcal
A}$, $z_{1},...,z_{2(m-s)}$ and $\zw$, yields ${\mathcal
G}''$-basic functions $g_{1},...,g_{2s}$ such that
$(x_{1},...,x_{n},z_{1},...,z_{2(m-s)},g_{1},...,g_{2s})$ is a
system of coordinates along ${\mathcal F}_{1}$. In this system
$\zw$ and $\zw_{1}$ are written with constant coefficients [recall
that $\zw_{1}$ is a  constant linear combination of $dz_{j}\zex
dz_{k}$, $1\zmei j<k\zmei 2(m-s)$]; by consequence the restriction
of $G$ to ${\mathcal A}$ is written with constant coefficients
too and every $(dg_{i}\zci G)_{\zbv\mathcal A}$ is a constant linear
combination of ${dz_{1}}_{\zbv\mathcal A},...,{dz_{2(m-s)}}_{\zbv\mathcal A}$.
Therefore in coordinates  $(x_{1},...,x_{n},z_{1},...,z_{2(m-s)},
{\bar z}_{2(m-s)+1},...,{\bar z}_{2m})$
each $dg_{i}\zci G$ equals $d_{x}g_{i}\zci J$ plus
a constant linear combination of $dz_{1},...,dz_{2(m-s)}$.

But $g_{i}$ is ${\mathcal G}''$-basic, so $d_{x}g_{i}\zci J$ is a
functional combination of $\za_{1},...,\za_{r}$. Thus in
coordinates $(x,z)=(x_{1},...,x_{n},z_{1},...,z_{2m})$ where
$z_{2(m-s)+i}=g_{i}$, $i=1,...,2s$, one has:

$G=\zsu_{j=1}^{n}a_{j}(\zpar/\zpar x_{j})\zte dx_{j}
+\zsu_{j=1}^{2m}\zsu_{k=1}^{2(m-s)}c_{jk}(\zpar/\zpar z_{j})\zte
dz_{k}$
\medskip

\hskip 6truecm $+\zsu_{i=1}^{2s}(\zpar/\zpar
z_{2(m-s)+i})\zte\zb_{i}$
\medskip

\noindent where every $c_{jk}\zpe\mathbb{K}$ and
$\za_{1}\zex...\zex\za_{r}\zex\zb_{i}=0$, $i=1,...,2s$.

Now it suffices to set $G'=G-\zsu_{i=1}^{2s}(\zpar/\zpar
z_{2(m-s)+i})\zte\zb_{i}$  {\it for finishing the proof of theorem
2.1.}
\bigskip

{\bf 3. The case of an eigenvalue function}

In the foregoing section one has studied Veronese flags with
parameters when $\zlma_{\zbv\mathcal A}$ is $0$-deformable and
nilpotent (theorem 2.1). Here we will consider Veronese flags for
more general tensor field $\zlma_{\zbv\mathcal A}$, which will be
the main tool for establishing the splitting theorem of
bihamiltonian structures.

One starts introducing the notion of regular open set. Let
${\mathbb K}_{P}[t]$ be the polynomial algebra in one variable
over the ring of differentiable functions on a manifold $P$. A
polynomial $\zf\zpe{\mathbb K}_{P}[t]$ is said {\it irreducible}
if it is irreducible at every point of $P$. Two polynomials
$\zf,\zq\zpe{\mathbb K}_{P}[t]$ are called {\it relatively prime}
if they are at each point. Given a vector bundle $E$ over $P$, of
dimension ${\widetilde m}$, and a morphism $H:E\zfl E$ its
characteristic polynomial $\zf=\zsu_{j=0}^{\widetilde
m}h_{j}t^{j}$ belongs to ${\mathbb K}_{P}[t]$. Set
$g_{j}=trace(H^{j})$. Since $h_{0},...,h_{{\widetilde m}-1}$ are,
up to sign, the elementary symmetric polynomials of the roots and
each $g_{j}$ the sum of their $j$-th powers, every function
$g_{j}$ may be expressed as a rational polynomial of
$h_{0},...,h_{{\widetilde m}-1}$, and each function $h_{j}$ like a
rational polynomial of $g_{1},...,g_{\widetilde m}$. In particular
$g_{j}$ when $j\zmai{\widetilde m}+1$ equals a rational polynomial
of $g_{1},...,g_{\widetilde m}$.

One will say that $H:E\zfl E$ {\it has constant algebraic type} if
there exist relatively prime irreducible polynomials
$\zf_{1},...,\zf_{s}\zpe{\mathbb K}_{P}[t]$ and positive integers
$a_{jk}$, $j=1,...,r_{k}$, $k=1,...,s$, such that at each point
$p\zpe P$ the family $\{\zf^{a_{jk}}_{k}(p)\}$, $j=1,...,r_{k}$,
$k=1,...,s$, is that of elementary divisors of $H(p)$. Let
$f_{1},...,f_{\widetilde n}$ be the family of all significant
coefficient functions of $\zf_{1},...,\zf_{s}\zpe{\mathbb
K}_{P}[t]$; that is $f$ when $\zf_{k}=t+f$ and $f,g$ if
$\zf_{k}=t^{2}+ft+g$. Obviously $h_{0},...,h_{{\widetilde m}-1}$,
and by consequence each $g_{j}$, are rational polynomials of
$f_{1},...,f_{\widetilde n}$. Conversely, for every
point of $P$ there exist analytic functions
$\zl_{k}(u_{1},...,u_{\widetilde n})$ such that close to this
point $f_{k}=\zl_{k}(g_{1},...,g_{\widetilde n})$,
$k=1,...,{\widetilde n}$ (note that ${\widetilde n}\zmei
{\widetilde m}$). Indeed, assume the degree of every $\zf_{k}$
equals one (otherwise complexify $E$ and $H$); then ${\widetilde
n}=s$ and it suffices to remark that the polynomial map
$F:{\mathbb K}^{s}\zfl{\mathbb K}^{s}$, defined by
$F(z)=(\zsu_{k=1}^{s}b_{k}z_{k},\zsu_{k=1}^{s}b_{k}z_{k}^{2}
,....,\zsu_{k=1}^{s}b_{k}z_{k}^{s})$
where each $b_{k}=\zsu_{j=1}^{r_{k}}a_{jk}$, is a local
diffeomorphism on the open set $\{z\zpe{\mathbb K}^{s}\zbv
z_{j}\znoi z_{k} {\rm if} j\znoi k\}$ since the determinant of its
Jacobian matrix equals $c\zpr_{1\zmei j<k\zmei j}(z_{j}-z_{k})$
with $c\zpe {\mathbb K}-\{0\}$.

Let $B_{H}$ be the set of all points such that around them $H$ has
constant algebraic type.
\bigskip

{\bf Lemma 3.1.} {\it The set $B_{H}$ is open and dense.}
\bigskip

{\bf Proof.} One may suppose ${\mathbb K}={\mathbb C}$ by
complexifying $E$ and $H$ if necessary. Given $p\zpe P$ let $a$ be
a root of $\zf(p)$ of multiplicity $b$. Then if $\ze>0$ is small
enough and $q$ close to $p$, the sum of multiplicities of the
roots of $\zf(q)$ belonging to the disk $D_{\ze}(a)$ equals $b$;
indeed, this sum is the degree of the map $e^{i\zh}\zpe
S^{1}\zfl\zf(q)(a+\ze e^{i\zt})\zdbv \zf(q)(a+\ze
e^{i\zh})\zdbv^{-1}\zpe S^{1}$. Therefore the number of different
roots of $\zf$ is locally constant on a dense open set $P'$ of
$P$.

Now assume $p\zpe P'$. Then there exist $\ze>0$ and an open set
$p\zpe B\zco P'$ such that $\zf(q)$, $q\zpe B$, has just one root
on $D_{\ze}(a)$ and its multiplicity is $b$. Let $\zl$ be the
$(b-1)$-derivative of $\zf$ with respect to $t$. Then
$(\zpar\zl/\zpar t)(p,a)\znoi 0$ and by the implicit function
theorem applied to $\zl$ and $0\zpe{\mathbb C}$, shrinking $B$ if
necessary, there is a differentiable (holomorphic or $C^{\zinf}$)
function $f:B\zfl{\mathbb C}$ such that $-f(q)$, $q\zpe B$, is the
root of $\zf(q)$ on $D_{\ze}(a)$. Thus
$\zf=\zpr_{k=1}^{s}(t+f_{k})^{b_{k}}$ around $p$ where
$f_{1},...,f_{s}$ are differentiable functions,
$b_{1},...,b_{s}$ integers $\zmai 1$ and $\zpr_{1\zmei
j<k\zmei s }(f_{j}-f_{k})$ never vanishes.

Finally, remark that the functions $dimKer((H+f_{k}I)^{j})$ are
locally decreasing, so locally constant on a dense open set
$B'\zco B$. $\square$

Suppose that $E$ is a foliation and $N_{H}=0$; then
$jdg_{j+1}=(j+1)dg_{j}\zci H$.
Indeed, consider $(E,H)$ as a weak Veronese flag (if $codimE=0$
regard the problem on ${\mathbb K}\zpor P$ in the obvious way)
and apply lemma 1.2.
Therefore $\zin_{j=1}^{\widetilde
m}Kerdg_{j}(p)=\zin_{j=0}^{{\widetilde m}-1}Kerdh_{j}(p)$ is a
$H$-invariant vector subspace of $T_{p}P$ because each $g_{j}$,
$j\zmai{\widetilde m}+1$, is a function of
$g_{1},...,g_{{\widetilde  m}}$.

One will be say that a point $p\zpe P$ is {\it regular} if there
exists an open neighbourhood $B$ of $p$ such that:

\noindent (1) $H$ has constant algebraic type on $B$,

\noindent (2) $\zin_{j=1}^{\widetilde m}Kerdg_{j}$, restricted to
$B$, is a vector sub-bundle of $E$ and therefore a foliation,

\noindent (3) $H$ restricted to $\zin_{j=1}^{\widetilde
m}Kerdg_{j}$ has constant algebraic type on $B$.

By lemma 3.1, applied to $H$ and its restriction to
$\zin_{j=1}^{\widetilde m}Kerdg_{j}$, the set of all regular
points is a dense open set $P$, called the {\it regular open set}.

To remark that if $H$ has constant algebraic type on an open set
$D$, then $\zin_{j=1}^{\widetilde
m}Kerdg_{j}=\zin_{j=1}^{\widetilde n}Kerdf_{j}$ on it where
$f_{1},...,f_{\widetilde n}$ are the significant coefficient
functions of $\zf_{1},...,\zf_{s}$.

Consider a Veronese flag $({\mathcal F},\zlma,\zw,\zw_{1})$ on a manifold $P$
or at some point of $P$.
Let ${\mathcal A}$ be the foliation of the largest $\zlma$-invariant vector
subspace (as in section 1) and $\zp:P\zfl N$ a local quotient of $P$ by ${\mathcal A}$.
Set $codim{\mathcal F}=r$, $dim{\mathcal A}=2m$ and $dimN=n$. Then $N$ is endowed
with a $r$-codimensional Veronese web whose limit when $t\zfl\zinf$ equals the quotient
foliation ${\mathcal F}'={\mathcal F}/{\mathcal A}$; moreover $\zlma$ projects in the morphism
$\zlma'$ associated to this Veronese web.

Suppose that $\zf=(t-f)^{2m}$ is the characteristic polynomial of $\zlma_{\zbv{\mathcal A}}$;
then from lemma 1.2 follows $df\zci\zlma=fdf$ on ${\mathcal F}$.
Now assume that $df_{\zbv {\mathcal A}}$ never vanishes. Let $X_{f}$ be the
$\zw$-hamiltonian of $f$; then $L_{X_{f}}\zlma=0$ and
$\zlma(X_{f})=fX_{f}$ since $df\zci\zlma=fdf$ on ${\mathcal F}$. Denoted by ${\bar P}$ and
${\bar\zp}:P\zfl {\bar P}$, respectively, the local quotient of $P$ by $X_{f}$ and its canonical
projection. Consider coordinates $(y,z)=(y_{1},...,y_{n},z_{1},...,z_{2m})$ on $P$ such that
$dy_{1}=...=dy_{r}=0$ defines ${\mathcal F}$, $dy_{1}=...=dy_{n}=0$ the foliation ${\mathcal A}$, $f=z_{2m}$
and $X_{f}=-\zpar/\zpar z_{2m-1}$. Thus $(y_{1},...,y_{n})$ can be regarded as coordinates on $N$,
$(y_{1},...,y_{n},z_{1},...,z_{2m-2},z_{2m})$ as coordinates on $\bar P$
and $f$ as a function on this last manifold. Now it is obvious that
$Kerdf$ and ${\mathcal F}\zin Kerdf$ project in two foliations ${\bar {\mathcal F}}_{1}$ and  ${\bar {\mathcal F}}$
on $\bar P$, respectively, and $\zlma_{\zbv{\mathcal F}\zin Kerdf}$ projects in a morphism
${\bar\zlma}:{\bar {\mathcal F}}\zfl {\bar {\mathcal F}}_{1}$; moreover $({\bar {\mathcal F}},{\bar\zlma})$
is a weak Veronese flag along ${\bar {\mathcal F}}_{1}$ (locally any extension of $\bar\zlma$ can be
lifted to an extension of $\zlma$), whose foliation ${\bar {\mathcal A}}$ of the largest
$\bar\zlma$-invariant vector subspaces equals the projection of ${\mathcal A}\zin Kerdf$,
${\bar P}/{\bar{\mathcal A}}$ is identified
to $N\zpor B$, where $B$ is an open neighbourhood of $f(p)$ on
${\mathbb K}$, and ${\bar {\mathcal F}}_{1}$ projects in the
foliation of $N\zpor B$ by the first factor. Moreover, the
Veronese web induced by $({\bar{\mathcal F}},{\bar\zlma})$ on each
leaf $N\zpor \{b\}$ of this last foliation equals the pull-back,
by the first projection $\zp_{1}:N\zpor B\zfl N$, of that induced
by $({\mathcal F},\zlma)$.

On the other hand, since $i_{X_{f}}\zw=-df_{\zbv{\mathcal A}}$ and
$i_{X_{f}}\zw_{1}=-(df\zci\zlma)_{\zbv{\mathcal A}}=-fdf_{\zbv{\mathcal A}}$,
the vector field $X_{f}$ belongs to
$Ker(\zw_{\zbv {\mathcal A}\zin Kerdf})$ and $Ker({\zw_{1}}_{\zbv {\mathcal A}\zin Kerdf})$,
so $\zw_{\zbv {\mathcal A}\zin Kerdf}$  projects in a symplectic form $\bar\zw$ and
${\zw_{1}}_{\zbv {\mathcal A}\zin Kerdf}$ in a closed $2$-form ${\bar\zw}_{1}$, both on ${\bar{\mathcal A}}$;
besides ${\bar\zw}_{1}={\bar\zw}({\bar\zlma},\quad)$.
The family $({\bar{\mathcal F}},{\bar\zlma},\bar\zw,{\bar\zw}_{1})$
will be called the {\it symplectic reduction of} $({\mathcal F},\zlma,\zw,\zw_{1})$. For proving that
this family is a Veronese flag it suffices to check the third condition of the definition.

On $N$ consider coordinates $(x_{1},...,x_{n})$ and a $(1,1)$-tensor field

\noindent $J=\zsu_{j=1}^{n}a_{j}(\zpar/\zpar x_{j})\zte dx_{j}$ where $a_{1},...,a_{n}$ are scalars.
\bigskip

{\bf Theorem 3.1.} {\it Let $G$ be a $(1,1)$-tensor field, which
extends $\zlma$ and projects in $J$, defined around a point $p$ of
$P$ such that $({\mathcal F},\zlma,\zw,\zw_{1})$ is a Veronese flag
at this point. Assume that:

\noindent (a) the characteristic polynomial of $\zlma_{\zbv
\mathcal A}$ equals $(t-f)^{2m}$ where $df_{\zbv \mathcal A}$
never vanishes,

\noindent (b) the function $f$ does not take the values
$a_{1},...,a_{n}$,

\noindent (c) $p$ is a regular point of $\zlma_{\zbv \mathcal A}$,

\noindent (d) the symplectic reduction of $({\mathcal
F},\zlma,\zw,\zw_{1})$ is a Veronese flag at  ${\bar\zp}(p)$,

\noindent then around $p$ there exist a $(1,1)$-tensor field $G'$
extending $\zlma$ and projecting in $J$ and functions
$z_{1},...,z_{2m}$ such that
$(x,z)=(x_{1},...,x_{n},z_{1},...,z_{2m})$ is a system of
coordinates,

\centerline{$G'=\zsu_{j=1}^{n}a_{j}(\zpar/\zpar x_{j})\zte dx_{j}+
\zsu_{j,k=1}^{2m}h_{jk}(z)(\zpar/\zpar z_{j})\zte dz_{k}$}

\noindent and $\zw,\zw_{1}$ are expressed relative to ${dz_{1}}
_{\zbv \mathcal A},...,{dz_{2m}} _{\zbv \mathcal A}$ with
coefficient functions only depending on $z$.}
\bigskip

Lets us prove theorem 3.1. Consider closed $1$-forms
$\za_{1},...,\za_{r}$ defining ${\mathcal F}$; by modifying the order of
variables $x_{1},...,x_{n}$ if necessary one may suppose that
$dx_{1}\zex...\zex dx_{n-r}\zex\za_{1}\zex...\zex\za_{n-r}$ has
no zeros. Since
$df\zci\zlma=fdf$ on ${\mathcal F}$, one has $df\zci G=fdf
+\zsu_{j=1}^{r}h_{j}\za_{j}$. Now consider a vector field $Y\zpe
{\mathcal A}$ such that $Yf=1$ and
set $G_{1}=G-Y\zte (\zsu_{j=1}^{r}h_{j}\za_{j})$;
then $df\zci G_{1}=fdf$, which
allows us to assume $df\zci G=fdf$ by considering $G_{1}$
instead of $G$ and calling it $G$. On the other hand from
$d(df\zci G)=0$ follows $L_{X_{f}}\zlma=0$, that is
$\za_{1}\zex...\zex\za_{r}\zex L_{X_{f}}G=0$, whence
$L_{X_{f}}G=\zsu_{j=1}^{r}X_{j}\zte\za_{j}$; moreover
$X_{1},...X_{r}\zpe{\mathcal A}\zin Kerdf$. Indeed,
$0=L_{X_{f}}(fdf)
=L_{X_{f}}(df\zci G)=df\zci L_{X_{f}}G$ and the projection on
$N$ of $L_{X_{f}}G$ vanishes since $G$ projects in $J$ and $X_{f}$
in zero.

Around $p$ there exist vector fields $Y_{1},...Y_{r}\zpe{\mathcal
A}\zin Kerdf$ such that $[X_{f},Y_{j}]=-X_{j}$, $j=1,...,r$. Then
$L_{X_{f}}(G+\zsu_{j=1}^{r}Y_{j}\zte\za_{j})=0$ and, by the
same reason as before, we can suppose $df\zci G=fdf$ and
$L_{X_{f}}G=0$. Thus $G_{\zbv Kerdf}$ projects in a $(1,1)$-tensor
field ${\bar G}$ defined along ${\bar {\mathcal F}}_{1}$, which
extend ${\bar \zlma}$ and projects in $J$ (regarded along the
foliation of $N\zpor B$ by the first factor in the obvious way).

On the other hand $({\bar{\mathcal
F}},{\bar\zlma}-fI,{\bar\zw},{\bar\zw}_{1}-f{\bar\zw})$ is a
Veronese flag along ${\bar{\mathcal F}}_{1}$ because $f$ is basic
for this last foliation. Moreover near ${\bar \zp}(p)$ the tensor
field $({\bar\zlma}-fI)_{\zbv{\bar\mathcal A}}$ is nilpotent
(obvious!) and $0$-deformable. Indeed, $p$ regular implies the
existence of positive integers $k_{1},...,k_{s}$ such that, around
this point, $\{t^{k_{j}},t^{k_{j}}\}$, $j=1,...,s$, is the family
of elementary divisors of $(\zlma-fI)_{\zbv\mathcal A}$ [recall
that every elementary divisor of $\zlma_{\zbv {\mathcal A}}$
occurs an even number of times because $\zw_{1}=\zw(\zlma,\quad)$]
whereas
$\{\{t^{k_{j}},t^{k_{j}}\}_{j=1,...,s-1},t^{k_{s}},t^{k_{s}-1}\}$ is
that of $(\zlma-fI)_{\zbv{\mathcal A}\zin Kerdf}$; now a
straightforward calculation shows that, close to ${\bar \zp}(p)$,
the family of elementary divisors of
$({\bar\zlma}-fI)_{\zbv{\bar\mathcal A}}$ has to be
$\{\{t^{k_{j}},t^{k_{j}}\}_{j=1,...,s-1},t^{k_{s}-1},t^{k_{s}-1}\}$.

Thus theorem 2.1 applied to ${\bar G}-fI$, $({\bar{\mathcal
F}},{\bar\zlma}-fI,{\bar\zw},{\bar\zw}_{1}-f{\bar\zw})$ and
$\zsu_{j=1}^{n}(a_{j}-f)(\zpar/\zpar x_{j})\zte dx_{j}$ yields
coordinates $(x_{1},...,x_{n},z_{1},...,z_{2m-2})$ along
${\bar{\mathcal F}}_{1}$, which become coordinates coordinates
$(x_{1},...,x_{n},z_{1},...,z_{2m-2},z_{2m})$ on $\bar P$ by
setting $z_{2m}=f$, such that $dz_{2m}=0$ defines ${\bar{\mathcal
F}}_{1}$,

${\bar G}-z_{2m}I=\zsu_{j=1}^{n}(a_{j}-z_{2m})(\zpar/\zpar
x_{j})\zte dx_{j}$

\hskip 4truecm $+\zsu_{j,k=1}^{2m-2}a_{jk}(\zpar/\zpar z_{j})\zte
dz_{k}+\zsu_{j=1}^{r}{\bar X}_{j}\zte\za_{j}$

\noindent where ${\bar X}_{1},...{\bar X}_{r}\zpe{\bar{\mathcal
A}}$ and each $a_{jk}\zpe{\mathbb K}$, and
${\bar\zw},{\bar\zw}_{1}-z_{2m}{\bar\zw}$ are written with
constant coefficients. Moreover, considering
$G-\zsu_{j=1}^{r}X_{j}\zte\za_{j}$ instead of $G$ where each
$X_{j}\zpe {\mathcal A}\zin Kerdf$, commutes with $X _{f}$ and
projects in ${\bar X}_{j}$, and linearly rearranging coordinates
$(z_{1},...,z_{2m-2})$ allows us to suppose
${\bar\zw}=(\zsu_{k=1}^{m-1}dz_{2k-1}\zex dz_{2k})_{\zbv
{\bar{\mathcal A}}}$ and

\centerline{${\bar
G}-z_{2m}I=\zsu_{j=1}^{n}(a_{j}-z_{2m})(\zpar/\zpar x_{j})\zte
dx_{j}+\zsu_{j,k=1}^{2m-2}a_{jk}(\zpar/\zpar z_{j})\zte dz_{k}$.}

If we regard $z_{1},...,z_{2m-2},z_{2m}$ as functions on $P$ one
has $\zw=(\zsu_{k=1}^{m-1}dz_{2k-1}\zex dz_{2k}+\zm\zex
dz_{2m})_{\zbv {\mathcal A}}$. But $\zw$ is closed so $(d\zm\zex
dz_{2m})_{\zbv {\mathcal A}}=0$ and one may choose a function
$z_{2m-1}$ such that $dz_{2m-1}\zex dz_{2m}$ equals $\zm\zex
dz_{2m}$ on ${\mathcal A}$; that is to say
$\zw=(\zsu_{k=1}^{m}dz_{2k-1}\zex dz_{2k})_{\zbv {\mathcal A}}$
and $(x_{1},...,x_{n},z_{1},...,z_{2m})$ is a system of
coordinates around $p$. Now $f=z_{2m}$, $X_{f}=-\zpar/\zpar
z_{2m-1}$ and $G=\zsu_{j=1}^{n}a_{j}(\zpar/\zpar x_{j})\zte
dx_{j}+(\zpar/\zpar z_{2m-1})\zte\zt+T$ where $\zt$ is a
functional combination of $dx_{1},...,dx_{n}$ and $T$ of
$(\zpar/\zpar z_{j})\zte dz_{k}$, $j,k=1,...,2m$; moreover the
coefficient function of these combinations do not depend on
$z_{2m-1}$ since $L_{X_{f}}G=0$.

On the other hand
$\zw_{1}=z_{2m}\zw+\zW_{1}+\zb\zex({dz_{2m}}_{\zbv\mathcal A})$
where $\zW_{1}$ is a constant linear combination of $(dz_{j}\zex
dz_{k})_{\zbv\mathcal A}$, $1\zmei j<k\zmei 2m-2$, whose
restriction to $Kerdf$ projects in ${\bar\zw}_{1}-f{\bar\zw}$, and
$\zb$ a functional combination of ${dz_{1}}_{\zbv\mathcal
A},...,{dz_{2m-2}}_{\zbv\mathcal A}$ whose coefficient functions
do not depend on $z_{2m-1}$ (recall that
$i_{X_{f}}\zw_{1}=-(df\zci\zlma)_{\zbv\mathcal
A}=-fdf_{\zbv\mathcal A}$; in particular $L_{X_{f}}\zw_{1}=0$).
Therefore as $\zw_{1}=\zw(\zlma,\quad)=\zw(G,\quad)=\zw(T,\quad)$
one has:

\centerline{$T=z_{2m}I_{z}+H+(\zpar/\zpar
z_{2m-1})\zte\zb^{*}-Z\zte dz_{2m}$}

\noindent where $I_{z}=\zsu_{k=1}^{2m}(\zpar/\zpar z_{k})\zte
dz_{k}$, $Z$ is is the vector field functional combination of
$\zpar/\zpar z_{1},...,\zpar/\zpar z_{2m-2}$ defined by the
equation $\zw(Z,\quad)=\zb$, $\zb^{*}$ the extension of $\zb$ to
$TP$ such that $\zb^{*}(\zpar/\zpar x_{j})=0$, $j=1,...,n$, and
$H$ the constant linear combination of $(\zpar/\zpar z_{j})\zte
dz_{k}$, $j,k=1,...,2m-2$, satisfying the relation
$\zw(H,\quad)=\zW_{1}$.

Set $\zW=\zw-(dz_{2m-1}\zex dz_{2m})_{\zbv\mathcal A}$. Then
$\zW_{1}=\zW(H,\quad)$ and $i_{Z}\zW=\zb$.

Note that any tensor field on $P$, partial or not, without terms
including $\zpar/\zpar z_{2m-1}$, $\zpar/\zpar z_{2m}$,
$dz_{2m-1}$ or $dz_{2m}$, whose coefficients functions do not
depend on $z_{2m-1}$, can be regarded as a tensor field along the
foliation $Kerdz_{2m}$ on $\bar P$ by projecting, via $\bar\zp$,
its restriction to $Kerdz_{2m}$. In coordinates $(x,z)$ this is
equivalent to delete coordinate $z_{2m-1}$ and consider the tensor
field along the foliation $dz_{2m}=0$. For the sake of simplicity both
tensor fields will be denoted by the same letter. Among others
that will be the case of $Z,H,\zW,\zW_{1},\zb$ already defined and
${\widetilde Z}, {\widetilde{\zb}}$ defined later on.

Until the end of the proof of theorem 3.1 and for making
calculations easier, $D$ will denote the exterior derivative with
respect of variables $(z_{1},...,z_{2m-2})$, along ${\mathcal A}$ on $P$ or
$\bar{\mathcal A}$ on $\bar P$,
and ${\mathcal L}$ the
Lie derivative on $\bar P$. As $\zw_{1}=z_{2m}(\zW+(dz_{2m-1}\zex
dz_{2m})_{\zbv\mathcal A})
+\zW_{1}+\zb\zex({dz_{2m}}_{\zbv\mathcal A})$ is closed one has
$D\zb={\mathcal L}_{Z}\zW=-\zW$.

On the other hand $N_{G}(\zpar/\zpar
z_{2m},\quad)=L_{G(\zpar/\zpar z_{2m})}G-G\zci L_{(\zpar/\zpar
z_{2m})}G=-L_{Z}G-H+(\zpar/\zpar z_{2m-1})\zte\zl+Z'\zte dz_{2m}$
with $\zl\zex dx_{1}\zex...\zex dx_{n}\zex dz_{1}\zex...\zex
dz_{2m-2}=0$ and $Z'\zex (\zpar/\zpar
z_{1})\zex...\zex(\zpar/\zpar z_{2m-2})=0$. Now from
$\za_{1}\zex...\zex\za_{r}\zex N_{G}=0$ follows
$L_{Z}G=-H+\zsu_{j=1}^{r}Y_{j}\zte\za_{j} +(\zpar/\zpar
z_{2m-1})\zte\zl+Z'\zte dz_{2m}$, where
$Y_{1},...,Y_{r}\zpe{\mathcal A}$ because $G$ projects in $J$ on
$N$ and $Z$ in zero.

 Hence projecting on $\bar P$, that is to say
considering variables $(x_{1},...,x_{n},z_{1},...,z_{2m-2})$ and
parameter $z_{2m}$, yields ${\mathcal L}_{Z}{\bar G}={\mathcal
L}_{Z}H=-H+\zsu_{j=1}^{r}{\bar Y}_{j}\zte\za_{j}$ where ${\bar
Y}_{1},...,{\bar Y}_{r}\zpe{\bar{\mathcal A}}$, which is the
foliation defined by $dx_{1}=...dx_{n}=dz_{2m}=0$.

Set ${\widetilde G}={\bar
G}-z_{2m}I=\zsu_{j=1}^{n}(a_{j}-z_{2m})(\zpar/ \zpar x_{j}\zte
dx_{j})+H$. Then on $\bar P$ one has $\zW(Z,\quad)=\zb$,
$\zW_{1}=\zW({\widetilde G},\quad)=\zW(H,\quad)$ and ${\mathcal
L}_{Z}{\widetilde G}={\mathcal L}_{Z}H=-H+\zsu_{j=1}^{r}{\bar
Y}_{j}\zte\za_{j}$. Moreover $\za_{1}\zex...\zex\za_{r}\zex
N_{\widetilde G}= \za_{1}\zex...\zex\za_{r}\zex N_{\bar G}=0$.

Given an endomorphism $S$ of a vector space $V$ and a vector
$v\zpe V$, one will say that {\it $v$ is $S$-generic} if $v$ and
$S$ have the same minimal polynomial; in particular $v\znoi 0$ if $V\znoi \{0\}$.
\bigskip

{\bf Lemma 3.2.} {\it Close to $p$ on $P$ and to ${\bar\zp}(p)$ on
$\bar P$ the vector field $Z$ is $H$-generic.}
\bigskip

{\bf Proof.} First remark that, for any $q\zpe P$, $Z(q)$ is
$H$-generic if and only if $Z({\bar\zp}(q))$ is $H$-generic on
$\bar P$.

Assume $dim{\bar{\mathcal A}}\zmai 2$, otherwise the result is obvious.
Along ${\bar{\mathcal A}}$ each $\zW(H^{k},\quad)$, $k\zmai 0$, is
a constant $2$-form and ${\mathcal
L}_{Z}(\zW(H^{k},\quad))=-(k+1)\zW(H^{k},\quad)$; on the other
hand ${\mathcal L}_{Z}(\zW(H^{k},\quad))=D(\zW(H^{k}Z,\quad))$.
Now suppose the minimal polynomial of $H$ equals $t^{k+1}$. Then
$\zW(H^{k},\quad)$ never vanishes and by consequence $H^{k}Z$ only does
on a closed set of empty interior; otherwise ${\mathcal
L}_{Z}(\zW(H^{k},\quad))=0$ on some non-empty open set. Therefore $Z$
is $H$-generic almost everywhere around ${\bar\zp}(p)$ on $\bar P$
and $p$ on $P$.

A straightforward calculation shows that the minimal polynomial of
$(\zlma_{\zbv {\mathcal A}})(q)$ is $(t-f(q))^{k+2}$ if and only
if $Z(q)$ is $H$-generic. So $Z$ has to be $H$-generic close to
$p$ since this point is regular. $\square$.

Let ${\widetilde Z}$ be a second vector field defined around $p$,
functional combination of $\zpar/\zpar z_{1},...,\zpar/\zpar
z_{2m-2}$ with coefficient only depending on
$(z_{1},...,z_{2m-2})$, such that ${\widetilde Z}(p)=Z(p)$,
$L_{\widetilde Z}\zW=-\zW$ and $ L_{\widetilde
Z}\zW_{1}=-2\zW_{1}$. Then $L_{\widetilde Z}H=-H$ on $P$ and
${\cal L}_{\widetilde Z}{\widetilde G}={\cal L}_{\widetilde
Z}H=-H$ on $\bar P$. The existence of a such vector field is clear
since $\zW$ and $\zW_{1}$ are written with constant coefficients;
for example take as $\widetilde Z$ a suitable linear vector field
(just a linear algebra exercise) plus a constant one.

Set ${\widetilde \zb}=\zW({\widetilde Z},\quad)$. Then on $\bar P$
one has $D{\widetilde \zb}=-\zW$ and there is a function
$g=g(x,z_{1},...,z_{2m-2},z_{2m})$ such that $\zb={\widetilde
\zb}+Dg$ and $Dg({\bar\zp}(p))=0$; moreover $Z={\widetilde
Z}-X_{g}$ where $X_{g}$ is the $\zW$-hamiltonian of $g$ (recall
that $\zW$ is symplectic on ${\bar{\mathcal A}}$).

Given a $1$-form $\zm$ defined on a vector sub-bundle $E$ by
$\zm(KerH^{k})=0$ one means $\zm(v)=0$ for every $v\zpe E\zin
KerH^{k}$. It is clear that $Dg(KerH^{k})=0$ for some integer
$k\zmai 0$. The next step will be to show that our problem
reduces to the case $Dg(KerH^{k+1})=0$.

Set $\zb_{t}={\tilde{\zb}}+tDg$ and $Z_{t}=(1-t){\tilde Z}+tZ$,
$t\zpe{\mathbb K}$. Then $\zW(Z_{t},\quad)=\zb_{t}$, ${\mathcal
L}_{Z_{t}}\zW=-\zW$ and $Z_{t}({\bar\zp}(p))=Z({\bar\zp}(p))$, so
$Z_{t}({\bar\zp}(p))$ is $H$-generic. Moreover $({\mathcal
L}_{Z_{t}}{\widetilde G}+H)\zex\za_{1}\zex...\zex\za_{r}=
((1-t){\mathcal L}_{Z_{t}}{\widetilde G}+t{\mathcal
L}_{Z}{\widetilde G}+H)\zex\za_{1}\zex...\zex\za_{r}=0$ and
$({\mathcal L}_{X_{g}}{\widetilde
G})\zex\za_{1}\zex...\zex\za_{r}=0$ since $X_{g}={\widetilde
Z}-Z$. By technical reasons let us decompose the manifold $\bar P$
into a product of a $n$-manifold [variables $(x_{1},...,x_{n})$],
a $(2m-2)$-manifold [variables $(z_{1},...,z_{2m-2})$] a
$1$-manifold [variable $(z_{2m})$], and set
${\bar\zp}=(\zp_{1},\zp_{2},\zp_{3})$ following this
decomposition.

On ${\bar P}\zpor{\mathbb K}$ consider coordinates
$(x_{1},...,x_{n},z_{1},...,z_{2m-2},z_{2m},t)$ and the foliation
$dz_{2m}=dt=0$, that is variables
$(x_{1},...,x_{n},z_{1},...,z_{2m-2})$ and parameters
$(z_{2m},t)$. From proposition 4.1, applied to
$Z_{t}$, $\widetilde G$,
$\zW$, $\zW_{1}$, ${\bar{\mathcal A}}$ all of them regarded along
$dz_{2m}=dt=0$, $\zp_{1}(p)$, $\zp_{2}(p)$, the compact set
$K=\{\zp_{3}(p)\}\zpor [0,1]$ and $g$ when $a=c=c'=-1$, follows
the existence of a function $f_{t}$ defined around
$K'=\{{\bar\zp}(p)\}\zpor [0,1]$ such that:

\noindent (I) ${Df_{t}}_{\zbv K'}=0$,

\noindent (II) $({\mathcal L}_{X_{t}}{\widetilde
G})\zex\za_{1}\zex...\zex\za_{r}=0$ where $X_{t}$ is the
$\zW$-hamiltonian of $f_{t}$,

\noindent (III) $Df_{t}(KerH^{k})=0$ and $D(Z_{t}f_{t}+f_{t}-g)(KerH^{k+1})=0$.

Therefore:

\noindent (1) ${X_{t}}_{\zbv K'}=0$,

\noindent (2) $X_{t}$ is tangent to $ImH^{k}$,

 \noindent (3) ${\mathcal L}_{X_{t}}\zW={\mathcal
L}_{X_{t}}\zW_{1}=0$,

\noindent (4) $({\mathcal L}_{X_{t}}\zb_{t}-Dg)(KerH^{k+1})=0$.

Indeed, assertions (1) is clear and (3) follows from he fact that
$\zW_{1}=\zW({\widetilde G},\quad)$.
For checking (2)
remark that $ImH^{k}$ is the $\zW$-orthogonal of ${\bar{\mathcal
A}}\zin KerH^{k}$ and $Ker\zW(X_{t},\quad)=KerDf_{t}\zcco
{\bar{\mathcal A}}\zin KerH^{k}$. Finally, for assertion (4) take
into account that ${\mathcal L}_{X_{t}}\zb_{t}={\mathcal
L}_{X_{t}}(\zW(Z_{t},\quad))=\zW([X_{t},Z_{t}],\quad)=-{\mathcal
L}_{Z_{t}}(\zW(X_{t},\quad))+({\mathcal
L}_{Z_{t}}\zW)(X_{t},\quad)=D(Z_{t}f_{t}+f_{t})$ and apply (III).

Let $\zQ_{s}$ be the flow of the time depending vector field $-X_{t}$ (on
${\bar P}\zpor{\mathbb K}$ one considers the vector field $\zpar/\zpar t-X_{t}$).
As ${X_{t}}_{\zbv K'}=0$, $\zQ_{1}$ is defined around ${\bar\zp}(p)$ and it can
be regarded as a germ  of diffeomorphism at this point. By construction
$\zQ_{1}$ preserves ${\bar\zp}(p)$, $\zW$, $\zW_{1}$, $\za_{1},...,\za_{r}$,
${\widetilde G}_{\zbv{\bar{\mathcal A}}}=H_{\zbv{\bar{\mathcal A}}}$ and
${\widetilde G}\zex\za_{1}\zex...\zex\za_{r}$.

Since $X_{t}$ is tangent to ${\bar{\mathcal A}}$ one has

$\zQ_{1}(x,z_{1},...,z_{2m-2},z_{2m})=(x,\zF(x,z_{1},...,z_{2m-2},z_{2m}),z_{2m})$.

\noindent Thus the pull-back by $\zQ_{1}$ of $\widetilde G$ equals
${\widetilde G}+\zsu_{j=1}^{r}{\bar Y}_{j}\zte\za_{j}$ with
${\bar Y}_{1},...,{\bar Y}_{r}\zpe{\bar{\mathcal A}}$,
and that of $\bar G$ equals ${\bar G}+\zsu_{j=1}^{r}{\bar Y}_{j}\zte\za_{j}$

Now we construct a germ of diffeomorphism $F$, at point $p$, by setting

\centerline{$F(x,z)=(x,\zF(x,z_{1},...,z_{2m-2},z_{2m}),z_{2m-1}
+\zf(x,z_{1},...,z_{2m-2},z_{2m}),z_{2m})$}

\noindent such that $F(p)=p$ and $F^{*}\zw=\zw$. Indeed,
$F^{*}\zw=\zW+{(dz_{2m-1}\zex dz_{2m})}_{\zbv{\mathcal A}}
+\zr\zex({dz_{2m}}_{\zbv{\mathcal A}})
+{(d\zf\zex dz_{2m})}_{\zbv{\mathcal A}}$ where $\zr$ is a
$1$-form along ${\mathcal A}$; as $0=d(F^{*}\zw)
=d\zr\zex({dz_{2m}}_{\zbv{\mathcal A}})$ one may choose $\zf$ in such a way that
${(d\zf\zex dz_{2m})}_{\zbv{\mathcal A}}=-\zr\zex({dz_{2m}}_{\zbv{\mathcal A}})$.

On the other hand:

\centerline {$F^{*}(z_{2m}\zw+\zb\zex({dz_{2m}}_{\zbv{\mathcal
A}}))=z_{2m}\zw+({\widetilde\zb}+Dh_{1})\zex({dz_{2m}}_{\zbv{\mathcal
A}})$}

\noindent where $h_{1}=h_{1}(x,z_{1},...,z_{2m-2},z_{2m})$ and
$Dh_{1}(KerH^{k+1})=0$ since $D\zb=D{\widetilde\zb}$ and $\zF$
transforms ${\widetilde\zb}_{\zbv KerH^{k+1}}$ in ${\zb}_{\zbv
KerH^{k+1}}$, while

\centerline{$F^{*}\zW_{1}=\zW_{1}+Dh_{2}\zex({dz_{2m}}_{\zbv{\mathcal
A}})$}

\noindent where $h_{2}=h_{2}(x,z_{1},...,z_{2m-2},z_{2m})$ because
$\zW_{1}$ is closed.

Let us see that $Dh_{2}(KerH^{k+1})=0$. Set

\centerline{ $\zQ_{s}(x,z_{1},...,z_{2m-2},z_{2m})
=(x,\zq_{s}(x,z_{1},...,z_{2m-2},z_{2m}),z_{2m})$.}

Since $X_{t}$
is tangent to ${\bar{\mathcal A}}\zin ImH^{k}=ImH^{k}$, the flow
$\zQ_{s}$ respects each leaf of the foliation  ${\bar{\mathcal
A}}\zin ImH^{k}$. Therefore $(z_{1},...,z_{2m-2})$ and
$\zq_{s}(x,z_{1},...,z_{2m-2},z_{2m})$ belong to the same leaf of
the foliation $ImH^{k}$ regarded on the variables
$(z_{1},...,z_{2m-2})$ only; by consequence
${(\zq_{s})}_{*}(\zpar/\zpar z_{2m})\zpe ImH^{k}$ and, in
particular, $\zF_{*}(\zpar/\zpar z_{2m})\zpe ImH^{k}$, whence
$(F_{*}(\zpar/\zpar z_{2m})-\zpar/\zpar z_{2m})\zpe ImH^{k}$; so
set $F_{*}(\zpar/\zpar z_{2m})=\zpar/\zpar z_{2m}+H^{k}V$. As $F$
respects ${\mathcal A}$, $Kerdz_{2m}$ and $H_{\zbv {\mathcal
A}\zin Kerdz_{2m}}$ one has $F_{*}({\mathcal A}\zin KerH^{k+1}\zin
Kerdz_{2m}) ={\mathcal A}\zin KerH^{k+1}\zin Kerdz_{2m}$. But on
$P$, $\zW_{1}(\zpar/\zpar z_{2m},\quad)=0$, $\zW_{1}=\zW(H,\quad)$
and $\zW(H,\quad)=\zW(\quad,H)$, therefore
\medskip

 $(F^{*}\zW_{1})(\zpar/\zpar z_{2m},{\mathcal A}\zin
KerH^{k+1}\zin Kerdz_{2m})$

\hskip 1.2truecm $=\zW_{1}(F_{*}(\zpar/\zpar z_{2m}),{\mathcal
A}\zin KerH^{k+1}\zin Kerdz_{2m})$

\hskip 2.4truecm $=\zW(H^{k+1}V,{\mathcal A}\zin KerH^{k+1}\zin
Kerdz_{2m})$

\hskip 3.6truecm $=\zW(V,H^{k+1}({\mathcal A}\zin KerH^{k+1}\zin
Kerdz_{2m}))=0$
\medskip

\noindent which implies $Dh_{2}(KerH^{k+1})=0$.

In short
$F^{*}\zw_{1}=z_{2m}\zw+\zW_{1}+({\widetilde\zb}+Dh)\zex({dz_{2m}}_{\zbv{\mathcal
A}})$ where $h=h(x,z_{1},...,z_{2m-2},z_{2m})$ and
$Dh(KerH^{k+1})=0$.

Set $\zg={\widetilde\zb}+Dh$. Let $\zg^{*}$ be the extension of
$\zg$ to $TP$ such that $\zg(\zpar/\zpar x_{j})=0$, $j=1,...,n$,
and $U$ the vector field functional combination of $\zpar/\zpar
z_{1},...,\zpar/\zpar z_{2m-2}$ defined by $\zw(U,\quad)=\zg$.
Since $G$, up to the term $(\zpar/\zpar z_{2m-1})\zte\zt$,
is determined by $\zw$, $\zw_{1}$ and $\bar G$, its
pull-back $G_{F}$ by $F$ is determined by $F^{*}\zw=\zw$,
$F^{*}\zw_{1}=z_{2m}\zw+\zW_{1}+\zg\zex ({dz_{2m}}_{\zbv\mathcal
A})$ and $\bar G+\zsu_{j=1}^{r}{\bar Y}_{j}\zte\za_{j}$; therefore
reasoning as before yields
\medskip

$G_{F}=\zsu_{j=1}^{n}a_{j}(\zpar/\zpar x_{j})\zte
dx_{j}+(\zpar/\zpar z_{2m-1})\zte\zs+z_{2m}I_{z}+H$

\hskip 3truecm $+(\zpar/\zpar z_{2m-1})\zte\zg^{*}-U\zte
dz_{2m}+\zsu_{j=1}^{r}Y_{j}\zte\za_{j}$
\medskip

\noindent where $Y_{1},...,Y_{r}\zpe\mathcal A$ and $\zg$ is a
functional combination of $dx_{1},...,dx_{n}$ whose coefficients
do not depend on $z_{2m-1}$.

Clearly it suffices proving theorem 3.1 for $G_F$, $F^{*}\zw$ and
$F^{*}\zw_{1}$; even more the term
$\zsu_{j=1}^{r}Y_{j}\zte\za_{j}$ is irrelevant and may be deleted.
Thus a change of notation (denote $G_F$ by $G$, $\zg$ by $\zb$
etc...) allows us to assume $\zb={\widetilde\zb}+Dh$ where
$h=h(x,z_{1},...,z_{2m-2},z_{2m})$ and $Dh(KerH^{k+1})=0$.

Now we start the process again with ${\widetilde Z}+(Z-{\widetilde
Z})(p)$. Finally, after a finite number of steps, we may suppose
$Z={\widetilde Z}+W$ where $W$ is a constant vector field linear
combination of $\zpar/\zpar z_{1},...,\zpar/\zpar z_{2m-2}$. Thus
$Z$ only depend on $(z_{1},...,z_{2m-2})$.
\bigskip

{\bf Lemma 3.3.} {\it On an open set of ${\mathbb K}^{k+1}$
endowed with coordinates $(v,u)=(v_{1},...,v_{k},u)$ consider a
point $q=(q_{1},...,q_{k},{\bar q})$ and a tensor field
${\widetilde T}=uI+ {\widetilde H}-U\zte du$ where ${\widetilde
H}=\zsu_{i,j=1}^{k}a_{ij}(\zpar/\zpar v_{i})\zte dv_{j}$, with
each $a_{ij}$ constant, and
$U=\zsu_{j=1}^{k}\zf_{j}(v)(\zpar/\zpar v_{j})$. Assume that
$U(q)$ is ${\widetilde H}$-generic, ${\widetilde H}$ nilpotent and
$L_{U}{\widetilde H}=c{\widetilde H}$, $c\zpe{\mathbb K}$. Then
around $q$ there exist functions $h_{1},...,h_{k}$ of $v$ such
that $d(dh_{j}\zci{\widetilde T})=0$, $j=1,...,k$, and
$(dh_{1}\zex...\zex dh_{k}\zex du)(q)\znoi 0$.}
\bigskip

{\bf Proof.} Given $h=h(v)$ one has $d(dh\zci{\widetilde
T})=d(d\zci {\widetilde H})-d(Uh+h)\zex du$. Close to $q$ and
for every $j=1,...,k$, consider a function $g_{j}$ of $v$ such
that $g_{j}(q_{1},...,q_{k})=0$ and $dg_{j}=d(Uv_{j}+v_{j})$.
Then $d(dg_{j}\zci{\widetilde H})=0$; indeed,
$d(dv_{j}\zci{\widetilde H})=0$ and $d(Uv_{j})\zci{\widetilde
H}=(L_{U}dv_{j})\zci{\widetilde H}=L_{U}(dv_{j}\zci{\widetilde
H})-cdv_{j}\zci{\widetilde H}$. By proposition 4.2, applied
in coordinates $v=(v_{1},...,v_{k})$ with a zero dimensional space
of parameters, close to $(q_{1},...,q_{k})$ there exists
$f_{j}=f_{j}(v)$ such that $df_{j}(q_{1},...,q_{k})=0$,
$d(df_{j}\zci{\widetilde H})=0$ and $Uf_{j}=-f_{j}+g_{j}$.

Set $h_{j}=v_{j}+f_{j}$; then $d(dh_{j}\zci{\widetilde T})=0$,
$j=1,...,k$, and $(dh_{1}\zex...\zex dh_{k}\zex du)(q)
=(dv_{1}\zex...\zex dv_{k}\zex du)(q)$. $\square$

Let us come back to the proof of theorem 3.1. If
$h=h(z_{1},...,z_{2m-2},z_{2m})$ one has $dh\zci G=z_{2m}dh+dh\zci
H-(Zh)dz_{2m}$. Thus we can apply lemma 3.3, in variables
$(z_{1},...,z_{2m-2},z_{2m})$ when ${\widetilde T}=z_{2m}I+H-Z\zte
dz_{2m}$, for concluding the existence close to $p$ of functions
$h_{j}(z_{1},...,z_{2m-2})$, $j=1,...,2m-2$, such that
$d(dh_{j}\zci G)=0$ and $(dh_{1}\zex...\zex dh_{2m-2}\zex
dz_{2m})(p)\znoi 0$.

Denote by $X_{j}$ the $\zw$-hamiltonian of $h_{j}$, $j=1,...,2m-2$. From the
third condition of Veronese flag follows
$(L_{X_{j}}G)\zex\za_{1}\zex...\zex\za_{r}=0$
(everywhere or close to point $p$),
which in turn implies
$(L_{X_{j}}\zt)\zex\za_{1}\zex...\zex\za_{r}=0$, $j=1,...,2m-2$. Now set
$\zt=\zsu_{k=1}^{n-r}\zf_{k}dx_{j}+\zsu_{k=n-r+1}^{n}\zf_{k}\za_{k+r-n}$;
then $X_{j}\zf_{k}=0$, $j=1,...,2m-2$, $k=1,...,n-r$. But
$(X_{1}\zex...\zex X_{2m-2}\zex(\zpar/\zpar z_{2m-1}))(p)\znoi 0$ because
$(dh_{1}\zex...\zex dh_{2m-2}\zex dz_{2m})(p)\znoi 0$, so each
$\zf_{k}$, $k=1,,,.n-r$, only depends on $(x,z_{2m})$. Besides the term
$(\zpar/\zpar z_{2m-1})\zte(\zsu_{k=n-r+1}^{n}\zf_{k}\za_{k+r-n})$ is
irrelevant for our purpose and it may be deleted. In short, one
can suppose that $\zt$ only depends on $(x,z_{2m})$.

On the other hand
\vskip .2truecm

\noindent $N_{G}(\zpar/\zpar x_{i},\zpar/\zpar x_{j})=
[a_{i}\zpar/\zpar x_{i}+\zt(\zpar/\zpar x_{i})\zpar/\zpar z_{2m-1},
a_{j}\zpar/\zpar x_{j}+\zt(\zpar/\zpar x_{j})\zpar/\zpar z_{2m-1}]$
\vskip .2truecm

\noindent $-G[\zpar/\zpar x_{i},a_{j}\zpar/\zpar x_{j}+\zt(\zpar/\zpar x_{j})\zpar/\zpar z_{2m-1}]
-G[a_{i}\zpar/\zpar x_{i}+\zt(\zpar/\zpar x_{i})\zpar/\zpar z_{2m-1},\zpar/\zpar x_{j}]$
\vskip .4truecm

$= \zgran\zpizq (a_{i}-z_{2m}) {\frac {\zpar(\zt(\zpar/\zpar x_{j})} {\zpar x_{i}}} -
(a_{j}-z_{2m}) {\frac {\zpar(\zt(\zpar/\zpar x_{i})} {\zpar x_{j}}}\zpder
{\frac {\zpar} {\zpar z_{2m-1}}}$
\vskip .2truecm

$=\zgran d_{x}(\zt\zci(J-z_{2m}I_{x})^{-1})
((J-z_{2m}I_{x})(\zpar/\zpar x_{i}),(J-z_{2m}I_{x})(\zpar/\zpar x_{j}))
{\frac {\zpar} {\zpar z_{2m-1}}}$
\vskip .4truecm

\noindent where $I_{x}=\zsu_{j=1}^{n}(\zpar/\zpar x_{j})\zte dx_{j}$ and $d_{x}$
is the exterior derivative with respect to $x$.

Therefore $N_{G}\zex\za_{1}\zex...\zex\za_{r}=0$ implies
$d_{x}(\zt\zci(J-z_{2m}I_{x})^{-1})\zex(\za_{1}\zci(J-z_{2m}I_{x})^{-1})
\zex...\zex(\za_{r}\zci(J-z_{2m}I_{x})^{-1})=0$. In other words,
$\zt\zci(J-z_{2m}I_{x})^{-1}$ is closed along the foliation in
variables $(x,z_{2m})$ defined by
$\za_{1}\zci(J-z_{2m}I_{x})^{-1},...,\za_{r}\zci(J-z_{2m}I_{x})^{-1},dz_{2m}$,
and near $p$ there exists a function $h=h(x,z_{2m})$ such that
$d_{x}h$ equals $\zt\zci(J-z_{2m}I_{x})^{-1}$ modulo
$\za_{1}\zci(J-z_{2m}I_{x})^{-1},...,\za_{r}\zci(J-z_{2m}I_{x})^{-1}$.
Thus adding a suitable functional combination of
$\za_{1},...,\za_{r}$ to $\zt$ allows us to suppose
$\zt\zci(J-z_{2m}I_{x})^{-1}=d_{x}h$. Then

\noindent $d(z_{2m-1}-h)\zci G=z_{2m}dz_{2m-1}+\zt-z_{2m}(\zpar
h/\zpar z_{2m})dz_{2m}-(d_{x}h)\zci J+\zb^{*}=z_{2m}dz_{2m-1}+
d_{x}h\zci(J-z_{2m}I_{x})-z_{2m}(\zpar h/\zpar
z_{2m})dz_{2m}-(d_{x}h)\zci
J+\zb^{*}=z_{2m}d(z_{2m-1}-h)+\zb^{*}$.

Finally, if the coordinate $z_{2m-1}$ is replaced by ${\widetilde
z}_{2m-1}=z_{2m-1}-h$ and next ${\widetilde z}_{2m-1}$ is called
$z_{2m-1}$, as $h$ only depends on $(x,z_{2m})$, then the
expression of $\zw$ and that of $\zw_{1}$ are not modified whereas

\centerline{$G=\zsu_{j=1}^{n}a_{j}(\zpar/\zpar
x_{j})+z_{2m}I_{z}+H+(\zpar/\zpar z_{2m-1})\zte\zb^{*}-Z\zte
dz_{2m}$}

\noindent {\it which proves theorem 3.1.}
\bigskip

{\bf 4. The equation $Zf=af+g$}

In this section, rather technical, one will establish the results
on the foregoing equation needed in the proof of theorem 3.1.
Consider three open sets $A\zco {\mathbb K}^{n}$, $A'\zco {\mathbb
K}^{2m}$ and $B\zco {\mathbb  K}^{\bar s}$, their product $A\zpor
A'\zpor B\zco{\mathbb  K}^{n+2m+\bar s}$ endowed with product
coordinates
$(x,z,u)=(x_{1},...,x_{n},z_{1},...,z_{2m},u_{1},...,u_{\bar s})$ and
the following objects on it:

\noindent ${\mathcal F}_{1}$: foliation defined by setting
$u_{1},...u_{\bar s}$ constant,

\noindent $d$: exterior derivative along ${\mathcal F}_{1}$,

\noindent ${\mathcal A}$: foliation contained in ${\mathcal
F}_{1}$ defined by $dx_{1}=...=dx_{n}=0$,

\noindent $D$: exterior derivative along ${\mathcal A}$,

\noindent $\zw,\zw_{1}$: couple of $2$-forms on ${\mathcal A}$,

\noindent $Z=\zsu_{j=1}^{2m}\zf_{j}\zpar/\zpar z_{j}$: vector
field tangent to ${\mathcal A}$.

Along ${\mathcal F}_{1}$, that is to say regarding $B$ as the
space of parameters, set $J=\zsu_{j=1}^{n}a_{j}(u)(\zpar/\zpar
x_{j})\zte dx_{j}$, $H=\zsu_{j,k=1}^{2m}a_{jk}(\zpar/\zpar
z_{j})\zte dz_{k}$ where each $a_{jk}\zpe{\mathbb K}$,
$\zx=\zsu_{j=1}^{r}X_{j}\zte\za_{j}$ where $X_{1},...,X_{r}\zpe
{\mathcal A}$ and $\za_{1},...,\za_{r}$ are closed $1$-forms
functional combination of $dx_{1},...,dx_{n}$ [so their
coefficients only depend on $(x,u)$] such that
$\za_{1}\zex...\zex\za_{r}$ never vanishes, and $G=J+H+\zx$.

Let ${\mathcal F}$ be the foliation contained in ${\mathcal
F}_{1}$ defined by $\za_{1},...,\za_{r}$. Suppose that
$\zw,\zw_{1}$ are written with constant coefficients which respect
to ${dz_{1}}_{\zbv {\mathcal   A}}, ...,{dz_{2m}}_{\zbv {\mathcal
A}}$, functions $a_{1}(u),...,a_{n}(u)$ never vanish on $B$, $H$
is nilpotent, $\zw_{1}=\zw(G,\quad)$
and $({\mathcal F},G_{\zbv{\mathcal F}})$
is a weak Veronese flag along ${\mathcal F}_{1}$ whose associated
$G_{\zbv{\mathcal   F}}$-invariant foliation equals ${\mathcal
A}$; therefore $\za_{1},...,\za_{r},J$ defines a Veronese web
along ${\mathcal F}_{1}/{\mathcal A}$ on $A\zpor B$

For any function $\zf$ one will denote $X_{\zf}$ its $\zw$-hamiltonian.
\bigskip

{\bf Lemma 4.1.} {\it Let $X_{f}$ be the $\zw$-hamiltonian of a
function $f$. Then $(L_{X_{f}}G)\zex\za_{1}\zex...\zex\za_{r}=0$
if and only if
$\za_{1}\zex...\zex\za_{r}\zex d(df\zci G)$ is
semi-basic for ${\mathcal A}$ (that is
$i_{U}(\za_{1}\zex...\zex\za_{r}\zex d(df\zci G))=0$
for any $U\zpe{\mathcal A}$).}

{\bf Proof.} As $(L_{X_{f}}\zx)\zex\za_{1}\zex...\zex\za_{r}=0$
and $d(df\zci\zx)\zex\za_{1}\zex...\zex\za_{r}=0$ we may suppose
$\zx=0$ without loss of generality. Now consider new coordinates
$(z_{j}^{i})$, $j=1,...,2m_{i}$, $i=1,...,s$, on $A'$ constant
linear combination of $(z_{1},...,z_{2m})$ such that
$\zw=(\zsu_{i=1}^{s}\zsu_{k=1}^{m_{i}}dz_{2k-1}^{i}\zex
dz_{2k}^{i})_{\zbv\mathcal A}$ and
$\zw_{1}=(\zsu_{i=1}^{s}\zsu_{k=1}^{m_{i}-1}dz_{2k-1}^{i}\zex
dz_{2k+2}^{i})_{\zbv\mathcal A}$. Then
$H=\zsu_{i=1}^{s}\zsu_{k=1}^{m_{i}-1}[(\zpar/\zpar z_{2k+1}^{i})\zte
dz_{2k-1}^{i}+(\zpar/\zpar z_{2k}^{i})\zte dz_{2k+2}^{i}]$
 since $\zw_{1}=\zw(G,\quad)
=\zw(H,\quad)$, and $L_{X_{f}}G=S+T$ where \vskip .3truecm

$\zgran S=\zsu_{i=1}^{s} {\frac{\zpar}{\zpar z_{1}^{i}}}\zte\zpizq
\zsu_{j=1}^{n}a_{j}{\frac{\zpar^{2}f}{\zpar z_{2}^{i}\zpar
x_{j}}}dx_{j}\zpder$\vskip .3truecm

\hskip 1.1truecm
$+\zgran\zsu_{i=1}^{s}\zsu_{k=1}^{m_{i}-1}{\frac{\zpar}{\zpar
z_{2k+1}^{i}}}\zte\zpizq \zsu_{j=1}^{n}\zcizq
a_{j}{\frac{\zpar^{2}f}{\zpar z_{2k+2}^{i}\zpar
x_{j}}}-{\frac{\zpar^{2}f}{\zpar z_{2k}^{i}\zpar x_{j}}}\zcder
dx_{j}\zpder$\vskip .3truecm

\hskip 2.1truecm $+\zgran\zsu_{i=1}^{s}
\zsu_{k=1}^{m_{i}-1}{\frac{\zpar}{\zpar z_{2k}^{i}}}\zte\zpizq
\zsu_{j=1}^{n}\zcizq {\frac{\zpar^{2}f}{\zpar z_{2k+1}^{i} \zpar
x_{j}}}-a_{j}{\frac{\zpar^{2}f}{\zpar z_{2k-1}^{i}\zpar
x_{j}}}\zcder dx_{j}\zpder$\vskip .3truecm

\hskip 3.1truecm $-\zgran\zsu_{i=1}^{s} {\frac{\zpar}{\zpar
z_{2m_{i}}^{i}}}\zte\zpizq
\zsu_{j=1}^{n}a_{j}{\frac{\zpar^{2}f}{\zpar z_{2m_{i}-1}^{i}\zpar
x_{j}}}dx_{j}\zpder$\vskip .3truecm

\noindent and $T$ does not involve any $\zpar/\zpar x_{j}$  nor
$dx_{j}$, $j=1,...,n$.

Note that $T=0$ if and only if  $(L_{X_{f}}G)_{\zbv\mathcal A}=0$.

On the other hand $df\zci G=\zsu_{j=1}^{n}a_{j}(\zpar f/\zpar
x_{j})dx_{j}$

\hskip 4truecm $+\zsu_{i=1}^{s}\zsu_{k=1}^{m_{i}-1}[ (\zpar
f/\zpar z_{2k+1}^{i})dz_{2k-1}^{i} +(\zpar f/\zpar
z_{2k}^{i})dz_{2k+2}^{i}]$.

Therefore $d(df\zci G)=\zr+\zl+\zm$ where \vskip .3truecm

\centerline{$\zr=\zgran\zsu_{i=1}^{s} \zsu_{k=1}^{m_{i}-1}\zpizq
\zsu_{j=1}^{n}\zcizq {\frac{\zpar^{2}f}{\zpar z_{2k+1}^{i} \zpar
x_{j}}}-a_{j}{\frac{\zpar^{2}f}{\zpar z_{2k-1}^{i}\zpar
x_{j}}}\zcder dx_{j}\zpder \zex dz_{2k-1}^{i}$}\vskip .3truecm

\centerline{$-\zgran\zsu_{i=1}^{s} \zpizq
\zsu_{j=1}^{n}a_{j}{\frac{\zpar^{2}f}{\zpar z_{2m_{i}-1}^{i}\zpar
x_{j}}}dx_{j}\zpder \zex dz_{2m_{i}-1}^{i}
-\zgran\zsu_{i=1}^{s}\zpizq
\zsu_{j=1}^{n}a_{j}{\frac{\zpar^{2}f}{\zpar z_{2}^{i}\zpar
x_{j}}}dx_{j}\zpder\zex dz_{2}^{i}$}\vskip .3truecm

\centerline{$+\zgran\zsu_{i=1}^{s}\zsu_{k=1}^{m_{i}-1} \zpizq
\zsu_{j=1}^{n} \zcizq {\frac{\zpar^{2}f}{\zpar z_{2k}^{i}\zpar
x_{j}}} -a_{j}{\frac{\zpar^{2}f}{\zpar z_{2k+2}^{i}\zpar x_{j}}}
\zcder dx_{j}\zpder\zex dz_{2k+2}^{i}\,$,}\vskip .3truecm

$\zl=\zsu_{jb}h_{jb}dx_{j}\zex dx_{b}$ and
$\zm=\zsu_{iajb}{\widetilde h}_{iajb}dz_{j}^{i}\zex dz_{b}^{a}$;
thus $(d(df\zci G))_{\zbv\mathcal A}=\zm_{\zbv\mathcal A}$.

But $L_{X_{f}}\zw_{1}=\zw(L_{X_{f}}G,\quad)$ and at the same time
$L_{X_{f}}\zw_{1}=D(\zw_{1}(X_{f},\quad))= -d(df\zci
G))_{\zbv\mathcal A}$; therefore $T=0$ if and only if $\zm=0$
since $\zw$ is symplectic.

Finally, remark that the $1$-forms functional combination of
$dx_{1},...,dx_{n}$ which are the coefficients of $\zpar/\zpar
z_{b}^{i}$ or $dz_{\bar b}^{i}$ in the expressions of $S$ and $\zr$,
respectively, are the same up to sign and change of order, so
$S\zex\za_{1}\zex...\zex\za_{r}=0$ if and only if
$\zr\zex\za_{1}\zex...\zex\za_{r}=0$. As $\zl$ is semi-basic for
${\mathcal A}$, the lemma is proved. $\square$

{\bf Remark.} From lemma 4.1 immediately follows that
$({\mathcal F},G_{\zbv{\mathcal F}},\zw,\zw_{1})$
is a Veronese flag along ${\mathcal F}_{1}$.
\bigskip

{\bf Proposition 4.1.} {\it Given an integer $k\zmai 0$, $p\zpe
A$, $q\zpe A'$, a compact set $K\zco B$, three scalars $a,c,c'$
and a function $g:A\zpor A'\zpor B\zfl {\mathbb K}$ such that:

\noindent (1) $(L_{Z}G-cH)\zex\za_{1}\zex...\zex\za_{r}=0$ and
$L_{Z}\zw=c'\zw$,

\noindent (2) $Z$ is $H$-generic on $\{(p,q)\}\zpor K$,

\noindent (3) $(L_{X_{g}}G)\zex\za_{1}\zex...\zex\za_{r}=0$ and
$Dg(KerH^{k})=0$,

\noindent then there exist a product open set $U\zpor U'\zpor
V\zco A\zpor A'\zpor B$,
which contains $\{(p,q)\}\zpor K$, and a function $f:U\zpor
U'\zpor V\zfl{\mathbb K}$ such that:

\noindent (I) $Df(KerH^{k})=0$ and $D(Zf-af-g)(KerH^{k+1})=0$,

\noindent (II) $(L_{X_{f}}G)\zex\za_{1}\zex...\zex\za_{r}=0$,

\noindent (III) $Df=0$ on $\{(p,q)\}\zpor V$.}
\bigskip

The next goal will be to prove proposition 4.1. Note that we may
assume $\zx=0$ since $(L_{X_{g}}\zx)\zex\za_{1}\zex...\zex\za_{r}
=(L_{X_{f}}\zx)\zex\za_{1}\zex...\zex\za_{r}=0$. Set
$\zf_{g}(x,z,u)=g(x,q,u)$, then $D\zf_{g}=0$ and $d(d\zf_{g}\zci
G)$ is ${\mathcal A}$-basic; thus any solution of our problem for
$g-\zf_{g}$  is a solution for $g$ too, which allows to suppose
$g(A\zpor \{q\}\zpor B)=0$ by considering $g-\zf_{g}$ instead of
$g$ if necessary. On the other hand by shrinking $A$ and modifying
the order of variables $(x_{1},...,x_{n})$ one may assume that
$dx_{1}\zex...\zex dx_{n-r}\zex\za_{1}\zex...\zex\za_{r}$ never
vanishes. By lemma 4.1, the first statement
(3) and part (II) of proposition 4.1 are
respectively equivalent to suppose $\za_{1}\zex...\zex\za_{r}\zex
d(dg\zci G)$ and $\za_{1}\zex...\zex\za_{r}\zex d(df\zci G)$
semi-basic for ${\mathcal A}$; {\it throughout the proof one will
use these statements instead of original ones}.

We start reducing the problem to the case $k=0$. Consider $H$ as a
tensor field on $A'$ and linearly rearrange coordinates $z$ in
such a way that $dz_{1}=...=dz_{2{\widetilde m}}=0$ defines
$KerH^{k}$. Let $A''$ be the quotient (close to $q$) of $A'$ by
$KerH^{k}$ endowed with coordinates $(z_{1},...,z_{2{\widetilde
m}})$, and $\zp:A\zpor A'\zpor B\zfl A\zpor A''\zpor B$ the
canonical projection. As $Dg(KerH^{k})=0$ there is a function
$\bar g$ on $A\zpor A''\zpor B$ such that $g={\bar g}\zci\zp$.
Obviously $Z,H,G,{\mathcal F},{\mathcal A}$ project in similar
object ${\bar Z},{\bar H},{\bar G},{\bar{\mathcal
F}},{\bar{\mathcal A}}$ defined on $A\zpor A''\zpor B$.

On the other hand $\zw(H^{k},\quad)$ and $\zw_{1}(H^{k},\quad)$
project in a couple of $2$-forms ${\bar\zw}$, ${\bar\zw}_{1}$
defined along ${\bar{\mathcal A}}$. It is easily checked that the
hypothesis of proposition 4.1 still hold on $A\zpor A''\zpor B$
for the scalars $a$, $c$ and $c'+kc$. Therefore, if the result is
proved for $k=0$, there exists a solution $\bar f$ and it suffices
to set $f={\bar f}\zci\zp$.

{\it In short $k=0$ is the only case to deal with.} We do that by
induction on the order ${\widetilde k}$ of nilpotency of $H$.
First consider the case $H=0$, that is ${\widetilde k}=1$ and
$G=J$. Assume $m\zmai 1$, otherwise it suffices setting $f=0$. As
$Z$ has no zeros on the compact set $\{(p,q)\}\zpor K$, we may
suppose that $\zf_{1}$ does not vanish on $A\zpor A'\zpor B$ by
shrinking these three factor and changing the order of variables
$z=(z_{1},...,z_{2m})$ if necessary. From (1) of proposition 4.1
follows $\za_{1}\zex...\zex\za_{r}\zex(d_{x}\zf_{j}\zci J)=0$,
$j=1,...,2m$, that is to say $(\za_{1}\zci
J^{-1})\zex...\zex(\za_{r}\zci J^ {-1})\zex d_{x}\zf_{j}=0$,
$j=1,...,2m$. Consider new coordinates $y=(y_{1},...,y_{n})$
around $p$ on $A$ such that $dy_{1}=...=dy_{r}=0$ defines the same
foliation as $\za_{1}\zci J^{-1},...,\za_{r}\zci J^{-1}$ (recall
that $\za_{1},...,\za_{r},J$ gives rise to a Veronese web).
Since every coordinate $y_{i}$ only depends on $x$ one has that
$d_{x}=d_{y}$ and each vector field $\zpar/\zpar z_{j}$,
$j=1,...,2m$, belongs to the dual basis of
$\{dy_{1},...,dy_{n},dz_{1},...,dz_{2m}\}$ as well. Thus
$\zf_{j}=\zf_{j}(y_{1},...,y_{r},z,u)$, $j=1,...,2m$.

On the other hand given a function $h$ then
$\za_{1}\zex...\zex\za_{r}\zex d(dh\zci G)$ is semi-basic for
${\mathcal A}$ if and only if $\za_{1}\zex...\zex\za_{r}\zex
(d_{x}(\zpar h/\zpar z_{j})\zci J)=0$, $j=1,...,2m$, that is
$(\za_{1}\zci J^{-1})\zex...\zex(\za_{r}\zci J^{-1})\zex
d_{x}(\zpar h/\zpar z_{j})=0$, $j=1,...,2m$, or in coordinates
$(y,z,u)$ if and only if $\zpar^{2} h/\zpar z_{j}\zpar y_{i}=0$,
$i=r+1,...,n$, $j=1,...,2m$. In these last coordinates consider the
open neighbourhoods of $p$ and $q$, respectively, $U=U_{1}\zpor
U_{2}$ and $U'$, where $U_{1}\zco{\mathbb K}^{r}$,
$U_{2}\zco{\mathbb K}^{n-r}$ and $U_{1},U_{2},U'$ are
polycylinders (that is product of open intervals if ${\mathbb K}={\mathbb R}$
or open disks if ${\mathbb K}={\mathbb C}$). Then
$\za_{1}\zex...\zex\za_{r}\zex d(dh\zci G)=0$ is semi-basic for
${\mathcal A}$ on $U\zpor U'\zpor B$ if and only if
$h=h_{1}(y_{1},...,y_{r},z,u)+h_{2}(y,u)$; moreover we may suppose
$h_{1}(p_{1},...,p_{r},q,u)=0$, $u\zpe B$, by taking
$h_{1}-\widetilde h$ and $h_{2}+\widetilde h$ where ${\widetilde
h}(y,z,u)=h_{1}(y_{1},...,y_{r},q,u)$ if necessary.

In particular as $\za_{1}\zex...\zex\za_{r}\zex d(dg\zci G)$ is
semi-basic for ${\mathcal A}$, on $U\zpor U'\zpor B$ one has
$g=g_{1}(y_{1},...,y_{r},z,u)+g_{2}(y,u)$ where
$g_{1}(p_{1},...,p_{r},q,u)=0$, $u\zpe B$. But $g(A\zpor
\{q\}\zpor B)=0$ so $g_{2}=0$, that is $g=g(y_{1},...,y_{r},z,u)$.

Let $f:U_{1}\zpor U'\zpor B\zfl{\mathbb K}$ be the function
defined by the ordinary differential equation $Zf=af+g$ and the
initial condition $f(U_{1}\zpor T\zpor B)=0$ where $T=\{z\zpe
U'\zbv z_{1}=q_{1}\}$ [again shrink $U_{1},U',B$ if necessary].
Then regarded on $U\zpor U'\zpor B$ in the obvious way
$\za_{1}\zex...\zex\za_{r}\zex d(df\zci G)$ is semi-basic for
${\mathcal A}$ and $D(Zf-af-g)=0$. By construction $(\zpar f/\zpar
z_{j})(\{(p,q)\}\zpor B)=0$, $j=2,...,2m$. But
$(Zf)(\{(p,q)\}\zpor B)=af(\{(p,q)\}\zpor B)+g(\{(p,q)\}\zpor
B)=0$ so  $(\zpar f/\zpar z_{1})(\{(p,q)\}\zpor B)=0$; in short
$Df=0$ on $\{(p,q)\}\zpor B$, {\it which proves proposition 4.1
when ${\widetilde k}=0$}.

{\bf Remark.} Observe that in this step proposition 4.1 was established without
making use of $\zw$ or $\zw_{1}$; therefore the result
stated in terms of being semi-basic for ${\mathcal A}$ is true
regardless the existence or not of $\zw,\zw_{1}$. This fact
implies that proposition 4.1 also holds if $Dg\zci H=0$;
even more in this case there exists $f$ satisfying (I), (II) and (III) such that
$Df\zci H=0$ and $D(Zf-af-g)=0$. Indeed,
regard $H$ as a tensor field on $A'$ and consider the quotient
$A''=A'/ImH$. Let $\zp:A\zpor A'\zpor B\zfl A\zpor A''\zpor B$ be
the canonical projection. Then $g={\bar g}\zci\zp$ for some ${\bar
g}: A\zpor A''\zpor B\zfl {\mathbb K}$. Since all the relevant
objects project on $A\zpor A''\zpor B$ and $H$ does in the zero
tensor field, from the case ${\widetilde k}=0$ follows the
existence of a suitable function ${\bar f}$ for ${\bar g}$ and it
suffices setting $f={\bar f}\zci\zp$.

Now suppose true proposition 4.1 up to the order of nilpotency
${\widetilde k}-1\zmai 1$ and for any scalars $a,c,c'$. One will
need the following result.
\bigskip

{\bf Lemma 4.2.} {\it Given a function $h:A\zpor A'\zpor B\zfl
{\mathbb K}$ such that $Dh(KerH)=0$ and
$\za_{1}\zex...\zex\za_{r}\zex d(dh\zci G)$ is semi-basic for
${\mathcal A}$, then there exist a product open set $U\zpor
U'\zpor V\zco A\zpor A'\zpor B$,
which contains $\{(p,q)\}\zpor K$, and a function
$\zf:U\zpor U'\zpor V\zfl{\mathbb K}$ such that:

\noindent (I) $D\zf\zci H=Dh$ and $\za_{1}\zex...\zex\za_{r}\zex
d(d\zf\zci G)$ is semi-basic for ${\mathcal A}$,

\noindent (II) $\zf(\{(p,q)\}\zpor V)=0$ and $D\zf(p,q,u)=0$ for
every $u\zpe V$ such that $Dh(p,q,u)=0$.}
\bigskip

{\bf Proof.} Consider coordinates $(z_{j}^{i})$,
$j=1,...,m_{i}$, $i=1,...,s$, on $A'$ as in the proof of lemma 4.1
and shrink this open set in such a way that, in these coordinates,
$A'$ is a polycylinder. Then
$dz_{2k+1}^{i}\zci H=dz_{2k-1}^{i}$, $dz_{2k}^{i}\zci
H=dz_{2k+2}^{i}$, $k=1,...,m_{i}-1$, $dz_{1}^{i}\zci
H=dz_{2m_{i}}^{i}\zci H=0$, $i=1,...,s$.

Since $Dh(KerH)=0$ one has $\zpar h/\zpar z_{2m_{i}-1}^{i}=\zpar
h/\zpar z_{2}^{i}=0$, $i=1,...,s$ and $\zb\zci H=Dh$ where\vskip
.3truecm

\centerline{$\zb=\zgran\zsu_{i=1}^{s} \zsu_{k=1}^{m_{i}-1}\zpizq
{\frac{\zpar h }{\zpar z_{2k-1}^{i}}} dz_{2k+1}^{i}+ {\frac{\zpar
h }{\zpar z_{2k+2}^{i}}} dz_{2k}^{i} \zpder_{\zbv\mathcal A}$.}
\vskip .3truecm

Note that $\za_{1}\zex...\zex\za_{r}\zex d(dh\zci G)$ semi-basic
for ${\mathcal A}$ implies $D(Dh\zci H)=0$. Now from lemma 1.1
follows that $D\zb(ImH,ImH)=0$, so ${\zb}_{\zbv ImH}$ is closed
and there exists a function $\zq:A\zpor A'\zpor B\zfl{\mathbb K}$
such that ${(D\zq-\zb)}_{\zbv ImH}=0$; therefore $D\zq\zci H=Dh$.
Hence $\zpar\zq/\zpar z_{2k+1}^{i}=\zpar h/\zpar z_{2k-1}^{i}$,
$\zpar\zq/\zpar z_{2k}^{i}=\zpar h/\zpar z_{2k+2}^{i}$,
$k=1,...,m_{i}-1$.

In general $\za_{1}\zex...\zex\za_{r}\zex d(d\zq\zci G)$ is not
semi-basic for ${\mathcal A}$ and we need to modify $\zq$.

If following the terminology of the proof of lemma 4.1 we set
$d(dh\zci G)={\zr}_{h}+{\zl}_{h}+{\zm}_{h}$  and $d(\zq\zci
G)={\zr}_{\zq}+{\zl}_{\zq}+{\zm}_{\zq}$ then
$\za_{1}\zex...\zex\za_{r}\zex{\zr}_{h}$ and
${\zm}_{h}={\zm}_{\zq}=0$ because $\za_{1}\zex...\zex\za_{r}\zex
d(dh\zci G)$ is semi-basic for ${\mathcal A}$ and $D(D\zq\zci
H)=D(Dh)=0$. Therefore $\za_{1}\zex...\zex\za_{r}\zex d(d\zq\zci
G)$ is semi-basic for ${\mathcal A}$ if and only if
$\za_{1}\zex...\zex\za_{r}\zex{\zr}_{\zq}=0$.

When $k=1,...,m_{i}-1$, $i=1,...,s$, the $1$-form coefficient of
$dz_{2k+1}^{i}$ in the expression of ${\zr}_{\zq}$ equals that of
$dz_{2k-1}^{i}$ in the expression of ${\zr}_{h}$, and the
coefficient of $dz_{2k}^{i}$ that of $dz_{2k+2}^{i}$ (recall that
$\zpar h/\zpar z_{2m_{i}-1}^{i}=\zpar h/\zpar z_{2}^{i}=0$), so
they vanish modulo $\za_{1},...,\za_{r}$. Thus we have only to
study the coefficients of $dz_{1}^{i}$ and $dz_{2m_{i}}^{i}$,
$i=1,...,s$, denoted by $\zb_{2i-1}$, $\zb_{2i}$ hereafter, which
are \vskip .3truecm

$\zb_{2i-1}=\zgran\zsu_{j=1}^{n}\zpizq {\frac{\zpar^{2}\zq}{\zpar
z_{3}^{i} \zpar x_{j}}}-a_{j}{\frac{\zpar^{2}\zq}{\zpar
z_{1}^{i}\zpar x_{j}}}\zpder dx_{j} =\zsu_{j=1}^{n}\zpizq
{\frac{\zpar^{2}h}{\zpar z_{1}^{i} \zpar
x_{j}}}-a_{j}{\frac{\zpar^{2}\zq}{\zpar z_{1}^{i}\zpar
x_{j}}}\zpder dx_{j}$\vskip .3truecm

\noindent and \vskip .3truecm

 $\zb_{2i}=\zgran\zsu_{j=1}^{n}\zpizq {\frac{\zpar^{2}\zq}{\zpar
z_{2m_{i}-2}^{i} \zpar x_{j}}} -a_{j}{\frac{\zpar^{2}\zq}{\zpar
z_{2m_{i}}^{i}\zpar x_{j}}}\zpder dx_{j}$ \vskip .3truecm

\hskip 3truecm $\zgran=\zsu_{j=1}^{n}\zpizq
{\frac{\zpar^{2}h}{\zpar z_{2m_{i}}^{i} \zpar
x_{j}}}-a_{j}{\frac{\zpar^{2}\zq}{\zpar z_{2m_{i}}^{i}\zpar
x_{j}}}\zpder dx_{j}\,\,$ respectively. \vskip .3truecm

For the sake of simplicity, set $z_{1}^{i}={\bar z}_{2i-1}$ and
$z_{2m_{i}}^{i}={\bar z}_{2i}$, $i=1,...,s$. From the expression
of $\zb_{j}$ and $\zb_{k}$ immediately follows
$\zpar\zb_{j}/\zpar{\bar z}_{k}=\zpar\zb_{k}/\zpar{\bar z}_{j}$,
$j,k=1,...,2s$. Moreover
$\za_{1}\zex...\zex\za_{r}\zex(\zpar\zb_{\bar a}/\zpar z_{\bar
b}^{i})=0 $, ${\bar a}=1,...,2s$, ${\bar b}=2,...,2m_{i}-1$,
$i=1,...,s$. Indeed, \vskip .3truecm

\centerline{$\zgran{\frac{\zpar\zb_{\bar a}}{\zpar z_{\bar
b}^{i}}} =\zsu_{j=1}^{n}{\frac{\zpar}{\zpar {\bar z}_{\bar a}}}
\zpizq {\frac{\zpar^{2}h}{\zpar z_{\bar b}^{i} \zpar
x_{j}}}-a_{j}{\frac{\zpar^{2}\zq}{\zpar z_{\bar b}^{i}\zpar
x_{j}}}\zpder dx_{j}$}\vskip .3truecm

\noindent which is the derivative with respect to ${\bar z}_{\bar
a}$ of the coefficient of $dz_{{\bar b}-2}^{i}$, if $\bar b$ is
odd, or $dz_{{\bar b}+2}^{i}$, if $\bar b$ is even, in the
expression of ${\zr}_{h}$. Thus, if we set $\zb_{\bar
a}={\bar\zb}_{\bar a}+\zb'_{\bar a}$ where ${\bar\zb}_{\bar a}$ is
a functional combination of $dx_{1},...,dx_{n-r}$ and $\zb'_{\bar
a}$ a functional combination of $\za_{1},...,\za_{r}$ (recall that
$dx_{1},...,dx_{n-r},\za_{1},...,\za_{r}$ are linearly independent
everywhere) one has $\zpar{\bar\zb}_{\bar a}/\zpar z_{\bar
b}^{i}=0$; that is every ${\bar\zb}_{\bar a}$, ${\bar
a}=1,...,2s$, only depends on $x$, ${\bar z}=({\bar
z}_{1},...,{\bar z}_{2s})$ and $u$.

Of course $\zpar{\bar\zb}_{j}/\zpar{\bar z}_{k}=\zpar{\bar\zb}_{k}/
\zpar{\bar z}_{j}$, $j,k=1,...,2s$.

The coefficient of $d{\bar z}_{k}$, $k=1,...,2s$, in the expression of $\zr_{h}$
equals
\vskip .3truecm

\centerline{$\zgran\zsu_{j=1}^{n}
\zpizq {\frac{\zpar^{2}h}{\zpar z_{\bar b}^{i} \zpar
x_{j}}}
-a_{j}{\frac{\zpar^{2}h}{\zpar{\bar  z}_{k}\zpar
x_{j}}}\zpder dx_{j}
=d_{x}\zpizq {\frac{\zpar h}{\zpar z_{\bar b}^{i}}}\zpder
-d_{x}\zpizq {\frac{\zpar h}{\zpar {\bar z}_{k}}}\zpder
\zci J$}\vskip .3truecm

\noindent where $i$ and $\bar b$ depend on $k$. This coefficient is a functional
combination of $\za_{1},...,\za_{r}$ therefore, as $\za_{1},...,\za_{r}$ define a
foliation, one has
$\za_{1}\zex...\zex\za_{r}\zex d_{x}(d_{x}(\zpar h/
\zpar{\bar z}_{k})\zci J)=0$. In turn from lemma 1.1 applied
to $d_{x}(\zpar h/\zpar{\bar z}_{k})\zci J^{-1}$ and $J$ follows
$(\za_{1}\zci J^{-1})\zex...\zex(\za_{r}\zci J^{-1})\zex d_{x}(d_{x}(\zpar h/
\zpar{\bar z}_{k})\zci J^{-1})=0$, whence taking into account the expression of
$\zb_{k}$ given before results
$(\za_{1}\zci J^{-1})\zex...\zex(\za_{r}\zci J^{-1})\zex
d_{x}(\zb_{k}\zci J^{-1})=0$, $k=1,...,2s$. Finally, since
$\za_{1}\zex...\zex\za_{r}\zex(\zb_{k}-{\bar\zb}_{k})=0$ and
$\za_{1}\zci J^{-1},...,\za_{r}\zci J^{-1}$ define a foliation,
one has $(\za_{1}\zci J^{-1})\zex...\zex(\za_{r}\zci J^{-1})\zex
d_{x}({\bar\zb}_{k}\zci J^{-1})=0$, $k=1,...,2s$.

After shrinking $A$ and $A'$ if necessary, we may suppose that $A$
in coordinates $y=(y_{1},...,y_{n})$ and $A'$ in coordinates
$(z_{j}^{i})$ are polycylinders. Set
$A=A_{1}\zpor A_{2}\zco{\mathbb K}^{r}\zpor{\mathbb K}^{n-r}$ and
$p=(p_{1},...,p_{n})$ following coordinates $(y_{1},...,y_{n})$.
Then there exist functions $f_{k}:A\zpor A'\zpor B\zfl{\mathbb
K}$, $k=1,...,2s$, only depending on $x$, ${\bar z}$ and $u$ such
that $f_{k}(A_{1}\zpor \{(p_{r+1},...,p_{n})\}\zpor A'\zpor B)=0$
and $(\za_{1}\zci J^{-1})\zex...\zex(\za_{r}\zci J^{-1})\zex
(d_{x}f_{k}- {\bar\zb}_{k}\zci J^{-1})=0$. Moreover $d_{\bar
z}(\zsu_{k=1}^{2s}f_{k}d{\bar z}_{k})=0$ where $d_{\bar z}$ is the
exterior derivative with respect to ${\bar z}=({\bar
z}_{1},...,{\bar z}_{2s})$. Indeed, $d_{x}(\zpar f_{k}/\zpar {\bar
z}_{j}-\zpar f_{j}/\zpar {\bar z}_{k})$, $j,k=1,...,2s$, equals
$(\zpar {\bar \zb}_{k}/\zpar {\bar z}_{j}-\zpar {\bar
\zb}_{j}/\zpar {\bar z}_{k})\zci J^{-1}=0$ modulo $\za_{1}\zci
J^{-1},...,\za_{r}\zci J^{-1}$, that is modulo
$dy_{1},...,dy_{r}$; in other words $\zpar f_{k}/\zpar {\bar
z}_{j}-\zpar f_{j}/\zpar {\bar z}_{k}$ does not depend on
$(y_{r+1},...,y_{n})$. But clearly $\zpar f_{k}/\zpar {\bar
z}_{j}-\zpar f_{j}/\zpar {\bar z}_{k}$ vanishes on $A_{1}\zpor
\{(p_{r+1},...,p_{n})\}\zpor A'\zpor B$ so $\zpar f_{k}/\zpar
{\bar z}_{j}-\zpar f_{j}/\zpar {\bar z}_{k}=0$, $j,k=1,...,2s$.

Thus there is a function $\zq_{1}:A\zpor A'\zpor B\zfl{\mathbb K}$
only depending on $x$, ${\bar z}$ and $u$ such that $\zpar
\zq_{1}/\zpar {\bar z}_{k}=f_{k}$, $k=1,...,2s$. Therefore
$\za_{1}\zex...\zex\za_{r}\zex (d_{x}(\zpar \zq_{1}/\zpar {\bar
z}_{k})\zci J-{\bar \zb}_{k})=0$.

Now set ${\widetilde \zf}=\zq+\zq_{1}$. Then $D{\widetilde
\zf}\zci H=D\zq\zci H=Dg$ and $\za_{1}\zex...\zex\za_{r}\zex
d(d{\widetilde \zf}\zci G)$ is semi-basic for ${\mathcal A}$, that
is ${\widetilde \zf}$ satisfies (I).

Finally, let ${\widetilde \zf}_{1}$ be the function given by
\vskip .3truecm

\centerline{$\zgran{\widetilde
\zf}_{1}(x,z,u)=\zsu_{k=1}^{2s}({\bar z}_{k}-{\bar
z}_{k}(q)){\frac {\zpar{\widetilde \zf}} {\zpar {\bar
z}_{k}}}(p,q,u)+{\widetilde \zf}(p,q,u)$;} \vskip .3truecm

\noindent then $\zf={\widetilde \zf}-{\widetilde \zf}_{1}$
satisfies (I) and (II). $\square$

A {\it box} [of coordinates $(z_{j}^{i})$] will mean a block of
coordinates $(z_{1}^{i},...,z_{2m_{i}}^{i})$ for any fixed $i$; so
one has $s$ boxes. A box will be named {\it short} if $m_{i}=1$
and {\it long} otherwise.

Consider a function $h:A\zpor A'\zpor B\zfl{\mathbb K}$ such that
$\za_{1}\zex...\zex\za_{r}\zex d(dh\zci G)$ is semi-basic for
${\mathcal A}$. Then there exists a function ${\widetilde h}$,
perhaps after shrinking $A'$, such that $Dh\zci H=D{\widetilde
h}$. From lemma 1.1, applied to $Dh$ and $H$ along ${\mathcal A}$,
follows $D(D{\widetilde h}\zci H)=D(Dh\zci H^{2})=0$, that is
$\zm_{\widetilde h}=0$ in the terminology of the proof of lemma
4.1. Moreover $\za_{1}\zex...\zex\za_{r}\zex d({\widetilde dh}\zci
G)$ is semi-basic for ${\mathcal A}$ since the coefficient of each
$dz_{j}^{i}$ in the expression of $\zr_{\widetilde h}$ equals that
of some $dz_{\bar b}^{\bar a}$ in the expression of $\zr_{h}$.

Observe that ${\widetilde h}$ does not depend on the short boxes.
By lemma 4.2, applied to ${\widetilde h}$ but considering long
boxes only, there exists a function $h_{1}:A\zpor A'\zpor
B\zfl{\mathbb K}$ [after shrinking we identify $A\zpor A'\zpor B$
and $U\zpor U'\zpor V$ for the sake of simplicity of the
notation], which does not depend on the short boxes, such that
$\za_{1}\zex...\zex\za_{r}\zex d(dh_{1}\zci G)$ is semi-basic for
${\mathcal A}$ and $Dh_{1}\zci H=Dh\zci H$. Now set
$h_{2}=h-h_{1}$; then $\za_{1}\zex...\zex\za_{r}\zex d(dh_{2}\zci
G)$ is semi-basic for ${\mathcal A}$ and $Dh_{2}\zci H=0$.

In other words, after shrinking $A$, $A'$ and $B$ if necessary,
the function $h$ decompose into a sum $h=h_{1}+h_{2}$ in such a
way that $\za_{1}\zex...\zex\za_{r}\zex d(dh_{1}\zci G)$ and
$\za_{1}\zex...\zex\za_{r}\zex d(dh_{2}\zci G)$ are semi-basic for
${\mathcal A}$, $h_{1}$ only depend on $x$, $u$ and the long
boxes, and $h_{2}$ does on $(x,{\bar z}_{1},...,{\bar z}_{2s},u)$
that is $Dh_{2}\zci H=0$.

Moreover $h_{1}(\{(p,q)\}\zpor B)=0$ and $Dh_{1}(p,q,u)=0$
whenever $(Dh\zci H)(p,q,u)=0$.

Consider a function $\zf:A\zpor A'\zpor B\zfl{\mathbb K}$ such
that $D\zf\zci H=0$ and $\za_{1}\zex...\zex\za_{r}\zex d(d\zf\zci
G)$ is semi-basic for ${\mathcal A}$. If there is one long box at
least, then {\it $Z$ and $Z+X_{\zf}$ are equivalent for the
purpose of proposition 4.1}. Let us see that. On $A\zpor A'\zpor
(B\zpor{\mathbb K})$ consider the vector field $Z_{t}=Z+tX_{\zf}$
and the obvious extensions of $G$, $\zw$, $\zw_{1}$ and ${\mathcal
A}$, where now the space of parameters is $B\zpor{\mathbb K}$
endowed with coordinates $(u_{1},...,u_{\bar s},t)$. Note that the
vector field $Z_{t}$ is $H$-generic on
$\{(p,q)\}\zpor(K\zpor\mathbb{K})$ since $X_{\zf}\zpe KerH$ and
$H\znoi \{0\}$.

By the remark preceding lemma 4.2 applied to $a=c'$, $c$, $c'$,
$Z_{t}$, $\zf$ and the compact set $K'=K\zpor [0,1]$, as $D\zf\zci
H=0$ there exists a function $f$ such that
$\za_{1}\zex...\zex\za_{r}\zex d(df\zci G)$ is semi-basic for
${\mathcal A}$, $Df\zci H=0$, $D(Z_{t}f-c'f-\zf)=0$ and $Df=0$ on
$\{(p,q)\}\zpor K'$. Then
$(L_{X_{f}}G)\zex\za_{1}\zex...\zex\za_{r}=0$,
$L_{X_{f}}\zw=L_{X_{f}}\zw_{1}=0$ since $\zw_{1}=\zw(G,\quad)$,
and $[X_{f},Z_{t}]=-X_{\zf}$ because $L_{Z_{t}}\zw=L_{Z}\zw=c'\zw$
and
$i_{[X_{f},Z_{t}]}\zw=i_{X_{f}}L_{Z_{t}}\zw-L_{Z_{t}}(i_{X_{f}}\zw)
=D(Z_{t}f-c'f)=D\zf$.

Let $\zQ_{\bar t}$ be the flow of the time depending vector field
$X_{f}$. As ${X_{f}}_{\zbv \{(p,q)\}\zpor K'}=0$, $\zQ_{1}$ is
defined around $\{(p,q)\}\zpor K$, preserves ${\mathcal
A},\zw,\zw_{1}$, and transforms $Z$ in $Z+X_{\zf}$ and $G$ in
$G+\zx$ where $\zx=\zsu_{j=1}^{r}X_{j}\zte\za_{j}$ and
$X_{1},...,X_{r}\zpe{\mathcal A}$. As $\zx$ is irrelevant, that is
the problem is the same for $G$ and $G+\zx$, the equivalence is
established.

Coming back to the main question, suppose $m_{i} \zmai 2$ when
$i=1,...,s'$ and $m_{i}=1$ otherwise. Set $Z=Z_{1}+Z_{2}$ where
$Z_{1}$ corresponds to the long boxes and $Z_{2}$ to the sort
ones. Remark that $Z_{1}$ is $H$-generic since $Z_{2}\zpe KerH$.
On the other hand set \vskip .3truecm

\centerline {$\zgran {\widetilde Z}=\zsu_{i=1}^{s} \zcizq
c\zsu_{k=1}^{m_{i}} \zpizq -kz_{2k-1}^{i}{\frac {\zpar} {\zpar
z_{2k-1}^{i}}}+ kz_{2k}^{i}{\frac {\zpar} {\zpar
z_{2k}^{i}}}\zpder  +{\frac {c'}
{2}}\zsu_{k=1}^{2m_{i}}z_{k}^{i}{\frac {\zpar} {\zpar
z_{k}^{i}}}\zcder$.} \vskip .3truecm

Then $L_{\widetilde Z}\zw=c'\zw$ and $L_{\widetilde Z}G=cH$ so
$(L_{\widetilde Z}G-cH)\zex\za_{1}\zex...\zex\za_{r}=0$. Thus,
after shrinking $A$, $A'$ and $B$ if necessary, there exists a
function $h=h_{1}+h_{2}$ such that ${\widetilde Z}=Z+X_{h}$,
$\za_{1}\zex...\zex\za_{r}\zex d(dh_{1}\zci G)$ and
$\za_{1}\zex...\zex\za_{r}\zex d(dh_{2}\zci G)$ are semi-basic for
${\mathcal A}$, $h_{1}$ does not depend on the short boxes and
$Dh_{2}\zci H=0$ (indeed, $L_{\widetilde Z}\zw=c'\zw=L_{Z}\zw$
implies that ${\widetilde Z}=Z+X_{h}$ for some $h$; now decompose
this function into a sum $h=h_{1}+h_{2}$ as it was showed before).

Since the components of $X_{h_{1}}$ in the short boxes vanish,
${\widetilde Z}=Z+X_{h_{1}}+X_{h_{2}}$ and the vector fields $Z$,
$Z+X_{h_{2}}$ are equivalent, we may assume  \vskip .3truecm

\centerline {$\zgran Z_{2}=\zsu_{i=s'+1}^{s} \zcizq
c\zsu_{k=1}^{m_{i}} \zpizq -kz_{2k-1}^{i}{\frac {\zpar} {\zpar
z_{2k-1}^{i}}}+ kz_{2k}^{i}{\frac {\zpar} {\zpar
z_{2k}^{i}}}\zpder  +{\frac {c'}
{2}}\zsu_{k=1}^{2m_{i}}z_{k}^{i}{\frac {\zpar} {\zpar
z_{k}^{i}}}\zcder$,} \vskip .3truecm

\noindent which implies that the coefficients functions of $Z_{1}$
do not depend on the short boxes, otherwise $L_{Z}\zw\znoi c'\zw$.

Decompose $g$ into a sum $g=g_{1}+g_{2}$ in such a way that
$g_{1}$ does not depend on the short boxes, $Dg_{2}\zci H=0$ and
$\za_{1}\zex...\zex\za_{r}\zex d(dg_{1}\zci G)$,
$\za_{1}\zex...\zex\za_{r}\zex d(dg_{2}\zci G)$ are semi-basic for
${\mathcal A}$. As proposition 4.1 was already proved for $g_{2}$
since $Dg_{2}\zci H=0$, it suffices to show it for $g_{1}$. But
$Z_{1}$ is $H$-generic and its components do not depend on the
short boxes, therefore it is enough considering the problem on the
long boxes only.

In other words, we may suppose that there is no short box. Note
that in this case $KerH\zco ImH$, therefore if $\zt$ is a $1$-form
such that $(\zt\zci H)(KerH^{2})=0$ then $\zt(KerH)=0$.

After shrinking $A$, $A'$ and $B$ if necessary, consider a
function ${\widetilde g}$ such that $D{\widetilde g}=Dg\zci H$;
then $\za_{1}\zex...\zex\za_{r}\zex d(d{\widetilde g}\zci G)$ is
semi-basic for ${\mathcal A}$ (it suffices reasoning as in the
case $D{\widetilde h}=Dh\zci H$). On the quotient $A\zpor
(A'/KerH)\zpor B$ one may project ${\widetilde g}$, $Z$, $G$,
${\mathcal A}$, $H$ and the $2$-forms $\zw_{1}$, which becomes
symplectic, and $\zw_{1}(G,\quad)$. Then the order of nilpotency
of the projection ${\widetilde H}$ of $H$ equals ${\widetilde
k}-1$ and by the induction hypothesis there is a solution of the
problem for the scalars $a+c$, $c$ and $c+c'$ [now
$L_{Z}\zw_{1}=L_{Z}(\zw(G,\quad))=(c+c')\zw_{1}$ and the same
equality holds on the quotient $A\zpor (A'/KerH)\zpor B$] and the
projection of ${\widetilde g}$. Pulling-back this solution yields
a function ${\widetilde f}:{\widetilde U}\zpor{\widetilde
U}'\zpor{\widetilde V}\zfl{\mathbb K}$ such that $D{\widetilde
f}(KerH)=0$, $D(Z{\widetilde f}-[a+c]{\widetilde f}-{\widetilde
g})(KerH^{2})=0$ since the pull-back of $Ker{\widetilde H}$ is
$KerH^{2}$, $\za_{1}\zex...\zex\za_{r}\zex d(d{\widetilde f}\zci
G)$ is semi-basic for ${\mathcal A}$ and $D{\widetilde f}=0$ on
$\{(p,q)\}\zpor {\widetilde V}$.

Let $f:U\zpor U'\zpor V\zfl{\mathbb K}$ a function given by lemma
4.2 applied to ${\widetilde f}$. Then $Df=0$ on $\{(p,q)\}\zpor V$
since $D{\widetilde f}=0$ on this set; besides $D(Zf-af-g)\zci
H=D(Z{\widetilde f}-[a+c]{\widetilde f}-{\widetilde g})$ because
$(L_{Z}Df)\zci H=L_{Z}(Df\zci H)-cDf\zci H$. But $D(Z{\widetilde
f}-[a+c]{\widetilde f}-{\widetilde g})(KerH^{2})=0$ therefore
$D(Zf-af-g)(KerH)=0$, {\it which finishes the proof of proposition
4.1}.

From now on and  until the end of this section, consider an open
set $A\zco{\mathbb K}^{m}$ endowed with coordinates
$z=(z_{1},...,z_{m})$, a manifold $B$ whose points are denoted by
$u$ and on $A\zpor B$ the following objects:

\noindent ${\mathcal G}$: foliation defined by setting $u$
constant,

\noindent $d$: exterior derivative along ${\mathcal G}$,

\noindent $Z$: vector field tangent to ${\mathcal G}$,

\noindent $H$: $(1,1)$-tensor field along ${\mathcal G}$.

Suppose that $H$ is nilpotent and written with constant
coefficients with respect to $(\zpar/\zpar z_{j})\zte dz_{k}$,
$j,k=1,...,m$, where $(z_{1},...,z_{m})$ are regarded as
coordinates along ${\mathcal G}$.
\bigskip

{\bf Proposition 4.2.} {\it Given an integer $k\zmai 0$, a point
$p\zpe A$, a compact set $K\zco B$, two scalars $a,c$ and a
function $g:A\zpor B\zfl{\mathbb K}$ such that:

\noindent (1) $L_{Z}H=cH$,

\noindent (2) $Z$ is $H$-generic on $\{p\}\zpor K$,

\noindent (3) $d(dg\zci H)=0$, $dg(KerH^{k})=0$ and $g(\{p\}\zpor
B)=0$,

\noindent then there exist a product open set $U\zpor V\zco A\zpor
B$, which contains $\{p\}\zpor K$, such that:

\noindent (I) $Zf=af+g$,

\noindent (II) $d(df\zci H)=0$ and $df(KerH^{k})=0$,

\noindent (III) $df=0$ on $\{p\}\zpor V$.} \bigskip

Let us proof proposition 4.2. Reasoning as in the proof of
proposition 4.1 reduces the problem to the case $k=0$. On the
other hand, if $H=0$ one has just a ordinary differential equation
and it suffices considering a solution $f$ that vanishes on a
suitable transverse section of $Z$ containing $\{p\}\zpor V$,
where $V$ is an open neighbourhood of $K$ on $B$ (note that
$(Zf)(\{p\}\zpor B)=0$ which implies $df=0$ on $\{p\}\zpor V$).

Now suppose that proposition 4.2 holds up to the dimension $m-1$
for any scalars $a,c$.
\bigskip

{\bf Lemma 4.3.} {\it Given a function $h:A\zpor B\zfl{\mathbb K}$
such that $dh(KerH)=0$ and $d(dh\zci H)=0$, then there exist a
product open set $U\zpor V\zco A\zpor B$, which contains
$\{p\}\zpor K$, and a function $\zf:U\zpor V\zfl{\mathbb K}$ such
that:

\noindent (I) $d\zf\zci H=dh$ and $d(d\zf\zci H)=0$,

\noindent (II) $\zf(\{p\}\zpor V)=0$ and $d\zf(p,u)=0$ for every
$u\zpe V$ such that $dh(p,u)=0$.}
\bigskip

{\bf Proof.} Consider coordinates $(z_{j}^{i})$, $j=1,...,m_{i}$,
$i=1,...,s$, on $A$, constant linear combination of
$(z_{1},...,z_{m})$, such that
$H=\zsu_{i=1}^{s}\zsu_{k=1}^{m_{i}-1}(\zpar/\zpar z_{k+1}^{i})\zte
dz_{k}^{i}$; that is $dz_{k}^{i}\zci H=dz_{k-1}^{i}$ if $k\zmai 2$
and $dz_{1}^{i}\zci H=0$. Therefore $\zpar h/\zpar
z_{m_{i}}^{i}=0$, $i=1,...,s$, and $\zb\zci H=dh$ where
$\zb=\zsu_{i=1}^{s}\zsu_{k=1}^{m_{i}-1}(\zpar h/\zpar
z_{k}^{i})dz_{k+1}^{i}$. On the other hand, after shrinking $A$,
we may suppose that in these coordinates $A$ is a polycylinder.

Now from lemma 1.1 follows that $\zb_{\zbv ImH}$ is closed;
therefore there exists a function $\zq:A\zpor B\zfl{\mathbb K}$
such that $(d\zq-\zb)_{\zbv ImH}=0$, whence $d\zq\zci H=dh$. Let
${\widetilde\zq}$ be the function given by ${\widetilde\zq}(z,u)=
\zsu_{i=1}^{s}(z_{1}^{i}-z_{1}^{i}(p))(\zpar\zq/\zpar
z_{1}^{i})(p,u)+\zq(p,u)$; then $\zf=\zq-{\widetilde\zq}$
satisfies (I) and (II). $\square$

Let us resume the proof of proposition 4.2. Since $d(dg\zci H)=0$,
after shrinking $A$ if necessary, there is a function ${\widetilde
g}:A\zpor B\zfl{\mathbb K}$ vanishing on $\{p\}\zpor B$ such that
$d{\widetilde g}=dg\zci H$; moreover by lemma 1.1 $d(d{\widetilde
g}\zci H)=0$. The equation $Z{\widetilde f}=(a+c){\widetilde
f}+{\widetilde g}$ has some solution satisfying (II) and (III) and
such that $d{\widetilde f}(KerH)=0$. Indeed, project $Z$, $H$ and
${\widetilde g}$ on the quotient $(A/KerH)\zpor B$, apply the
induction hypothesis and then pull-back a suitable solution. Note
that ${\widetilde f}$ is defined on a product open set containing
$\{p\}\zpor K$ and which we will call $A\zpor B$ for simplifying.

Let $\zf:U\zpor V\zfl{\mathbb K}$ be a function given by lemma 4.3
applied to ${\widetilde f}$. Set $g_{0}=g+a\zf-Z\zf$; then
$g_{0}(\{p\}\zpor V)=0$ and $dg_{0}\zci H=0$ since
$(L_{Z}d\zf)\zci H=L_{Z}(d\zf\zci H)-cd\zf\zci H$. In turn, the
equation $Zf_{0}=af_{0}+g_{0}$ has some solution satisfying (II)
and (III). Indeed, project $Z$, $H$ and $g_{0}$ on the quotient
$(U/ImH)\zpor V$ and reason as before. Now it suffices to set
$f=f_{0}+\zf$ {\it for finishing the proof of proposition 4.2}.
\bigskip

{\bf 5. The non-real eigenvalue case}

In this section ${\mathbb K}={\mathbb R}$ and the manifolds
considered will be real unless another thing is stated. Let
$({\mathcal F},\zlma,\zw,\zw_{1})$ be a Veronese flag on a
manifold $P$ or at some point of $P$, ${\mathcal A}$ the foliation
of the largest $\zlma$-invariant vector subspace (see section 1)
and $\zp:P\zfl N$ a local quotient of $P$ by ${\mathcal A}$. Set
$codim{\mathcal F}=r$, $dim{\mathcal A}=2m$ and $dimN=n$. Recall
that $N$ is endowed with a $r$-codimensional Veronese web whose
limit when $t\zfl\zinf$ equals the quotient foliation ${\mathcal
F}'={\mathcal F}/{\mathcal A}$ and $\zlma$ projects in the
morphism $\zlma'$ associated to this Veronese web.

Suppose that the characteristic polynomial $\zf$ of
$\zlma_{\zbv\mathcal A}$ equals $(t^{2}+ft+g)^{m}$ where
$f^{2}<4g$, that is $\zf$ has no real roots. Set
$g_{k}=trace((\zlma_{\zbv\mathcal A})^{k})$; by lemma 1.2 one has
$kdg_{k+1}=(k+1)dg_{k}\zci\zlma$ on ${\mathcal F}$.

When $df_{\zbv{\mathcal A}(p)}\znoi 0$ and the algebraic type of
$\zlma_{\zbv\mathcal A}$ is constant about $p$, one may construct
the symplectic reduction of the Veronese flag as follows. First
observe that each $g_{k}$ is function of $g_{1},g_{2}$ since $f,g$
are the only significant coefficients of the elementary divisors,
$g_{1}=-mf$, $g_{2}=m(f^{2}-2g)$ and $dg_{2}=2dg_{1}\zci\zlma$ on
${\mathcal F}$. Therefore $X_{g_{2}}=2\zlma X_{g_{1}}$ and
$(dg_{1}\zex dg_{2})_{\zbv {\mathcal A}(p)}\znoi 0$, otherwise
$(X_{g_{1}}\zex  X_{g_{2}})(p)=0$ and $\zlma$ has an eigenvalue on
${\mathcal A}(p)-\{0\}$. Thus $X_{g_{1}},X_{g_{2}},X_{g_{3}},...$
give rise to a $\zlma$-invariant vector sub-bundle $E$ of
dimension two.

On the other hand $\zw(X_{g_{1}},X_{g_{2}})=2\zw(X_{g_{1}},\zlma
X_{g_{1}})=\zw_{1}(X_{g_{1}},X_{g_{1}})=0$ and
$\zw_{1}(X_{g_{1}},X_{g_{2}})=2\zw(\zlma X_{g_{1}},\zlma
X_{g_{1}})=0$. Hence
$X_{g_{1}}g_{1}=X_{g_{1}}g_{2}=X_{g_{2}}g_{1}=X_{g_{2}}g_{2}=0$
and $[X_{g_{1}},X_{g_{2}}]=0$; in particular $E$ is a foliation.
Besides $L_{X_{g_{1}}}\zlma=L_{X_{g_{2}}}\zlma=0$ since
$kdg_{k+1}=(k+1)dg_{k}\zci\zlma$ on ${\mathcal F}$.

Denoted by ${\bar P}$ and ${\bar\zp}:P\zfl {\bar P}$,
respectively, the local quotient of $P$ by $E$ and its canonical
projection. Consider coordinates
$(y,z)=(y_{1},...,y_{n},z_{1},...,z_{2m})$ around $p$ such that
$dy_{1}=...=dy_{r}=0$ defines ${\mathcal F}$,
$dy_{1}=...=dy_{n}=0$ the foliation ${\mathcal A}$,
$g_{1}=z_{2m-1}$, $g_{2}=z_{2m}$, $X_{g_{1}}=-\zpar/\zpar
z_{2m-3}$ and $X_{g_{2}}=-\zpar/\zpar z_{2m-2}$. Thus
$(y_{1},...,y_{n})$ can be regarded as coordinates on $N$,
$(y_{1},...,y_{n},z_{1},...,z_{2m-4},z_{2m-1},z_{2m})$ as
coordinates on $\bar P$ and $g_{1},g_{2}$ as functions on this
last manifold. Now it is obvious that $Ker(dg_{1}\zex dg_{2})$ and
${\mathcal F}\zin Ker(dg_{1}\zex dg_{2})$ project in two
foliations ${\bar {\mathcal F}}_{1}$ and ${\bar {\mathcal F}}$ on
$\bar P$, respectively, and $\zlma_{\zbv{\mathcal F}\zin
Ker(dg_{1}\zex dg_{2})}$ does in a morphism ${\bar\zlma}:{\bar
{\mathcal F}}\zfl {\bar {\mathcal F}}_{1}$; moreover $({\bar
{\mathcal F}},{\bar\zlma})$ is a weak Veronese flag along ${\bar
{\mathcal F}}_{1}$ (locally any extension of $\bar\zlma$ can be
lifted to an extension of $\zlma$), whose foliation ${\bar
{\mathcal A}}$ of the largest $\bar\zlma$-invariant vector
subspaces equals the projection of ${\mathcal A}\zin
Ker(dg_{1}\zex dg_{2})$, ${\bar P}/{\bar{\mathcal A}}$ is
identified to $N\zpor B$, where $B$ is an open neighbourhood of
$(g_{1}(p),g_{2}(p))$ on ${\mathbb R}^{2}$, and ${\bar {\mathcal
F}}_{1}$ projects in the foliation of $N\zpor B$ by the first
factor. Besides, the Veronese web induced by $({\bar{\mathcal
F}},{\bar\zlma})$ on each leaf $N\zpor \{b\}$ of this last
foliation equals the pull-back, by the first projection
$\zp_{1}:N\zpor B\zfl N$, of that induced by $({\mathcal
F},\zlma)$.

On the other hand, since $i_{X_{g_{1}}}\zw$, $i_{X_{g_{2}}}\zw$,
$i_{X_{g_{1}}}\zw_{1}$ and $i_{X_{g_{2}}}\zw_{1}$ are functional
combination of ${dg_{1}}_{\zbv\mathcal A}$,
${dg_{2}}_{\zbv\mathcal A}$ (recall that every $g_{k}$ is function
of $g_{1},g_{2}$ and $kdg_{k+1}=(k+1)dg_{k}\zci\zlma$ on
${\mathcal F}$) one has $Ker(\zw_{{\zbv}\mathcal A\zin
Ker(dg_{1}\zex dg_{2})})$

\noindent $=Ker({\zw_{1}}_{{\zbv}\mathcal A\zin
Ker(dg_{1}\zex dg_{2})})=E$. Therefore $\zw_{{\zbv}\mathcal A\zin
Ker(dg_{1}\zex dg_{2})}$ and ${\zw_{1}}_{{\zbv}\mathcal A\zin
Ker(dg_{1}\zex dg_{2})}$ project in two symplectic forms
${\bar\zw}$, ${\bar\zw}_{1}$ along $\bar \mathcal A$; moreover
${\bar\zw}_{1}={\bar\zw}({\bar\zlma},\quad)$. The family
$({\bar{\mathcal F}},{\bar\zlma},\bar\zw,{\bar\zw}_{1})$ will be
called {\it the symplectic reduction (near $p$) of} $({\mathcal
F},\zlma,\zw,\zw_{1})$. As in section 3, for proving that this
family is a Veronese flag it suffices to check the third condition
of its definition.

On $N$ consider coordinates $(x_{1},...,x_{n})$ and a
$(1,1)$-tensor field

\noindent $J=\zsu_{j=1}^{n}a_{j}(\zpar/\zpar x_{j})\zte
dx_{j}$ where $a_{1},...,a_{n}$ are real numbers.
\bigskip

{\bf Theorem 5.1.} {\it In the real analytic category consider a
$(1,1)$-tensor field $G$, which extends $\zlma$ and projects in
$J$, defined around a point $p$ of $P$ such that $({\mathcal
F},\zlma,\zw,\zw_{1})$ is a Veronese flag at this point. Assume
that:

\noindent (a) the characteristic polynomial of $\zlma_{\zbv
\mathcal A}$ equals $(t^{2}+ft+g)^{m}$ where $f^{2}<4g$,

\noindent (b) $p$ is a regular point of $\zlma_{\zbv \mathcal A}$,

\noindent (c) if $df_{\zbv {\mathcal A}(p)}=0$ then $f$ is
constant close to $p$,

\noindent (d) if $df_{\zbv {\mathcal A}(p)}\znoi 0$ then the
symplectic reduction of $({\mathcal F},\zlma,\zw,\zw_{1})$ is a
Veronese flag at ${\bar\zp}(p)$,

\noindent then around $p$ there exist a $(1,1)$-tensor field $G'$
extending $\zlma$ and projecting in $J$ and functions
$z_{1},...,z_{2m}$ such that
$(x,z)=(x_{1},...,x_{n},z_{1},...,z_{2m})$ is a system of
coordinates,

\centerline{$G'=\zsu_{j=1}^{n}a_{j}(\zpar/\zpar x_{j})\zte dx_{j}+
\zsu_{j,k=1}^{2m}h_{jk}(z)(\zpar/\zpar z_{j})\zte dz_{k}$}

\noindent and $\zw,\zw_{1}$ are expressed relative to ${dz_{1}}
_{\zbv \mathcal A},...,{dz_{2m}} _{\zbv \mathcal A}$ with
coefficient functions only depending on $z$.}
\bigskip

Before proving this result, let us recall a few facts on the
relationship between complex and real manifolds. Let $Q$ be a real
manifold of dimension $2k$ endowed with a complex structure $H$,
which allows us to regard $Q$ as a complex manifold of dimension
$k$. A real tangent vector field at $q\zpe Q$ is a linear
derivation of the algebra of germs at this point of differentiable
functions; therefore it acts too on the germs at $q$ of
holomorphic functions and it can be regarded as a complex tangent
vector at this point. In other words, the real and the complex
tangent vector space at the same point may be identified in a
canonical way.

In turn, if $X$ is a real vector field then $L_{X}H=0$ if and only
if the (infinitesimal) action of $X$ sends holomorphic functions
into holomorphic functions. In terms of complex coordinates
${\mathbf z}_{1}={\mathbf x}_{1}+\zima{\mathbf y}_{1},...,{\mathbf
z}_{k}={\mathbf x}_{k}+\zima{\mathbf y}_{k}$ and the associated real
coordinates $({\mathbf x}_{1},{\mathbf y}_{1},...,{\mathbf
x}_{k},{\mathbf y}_{k})$, if
$X=\zsu_{j=1}^{k}(\zf_{j}\zpar/\zpar{\mathbf x}_{j}
+\zq_{j}\zpar/\zpar{\mathbf y}_{j})$ then $L_{X}H=0$ if and only
if $\zf_{1}+\zima\zq_{1},...,\zf_{k}+\zima\zq_{k}$ are holomorphic
functions of ${\mathbf z}=({\mathbf z}_{1},...,{\mathbf z}_{k})$;
in this case from the complex viewpoint
$X=\zsu_{j=1}^{k}(\zf_{j}+\zima\zq_{j})\zpar/\zpar{\mathbf z}_{j}$
(warning this identification only works for the action of $X$ on
holomorphic functions but not for any complex-valued function).
This kind of vector fields are named {\it holomorphic}.

A complex ${\widetilde k}$-form [that is of type $({\widetilde
k},0)$] $\zb$ decompose into a sum $\zb=\zb_{1}+\zima\zb_{2}$, where
$\zb_{1},\zb_{2}$ are real ${\widetilde k}$-forms such that
$\zb_{1}(H,\quad,...,\quad)=-\zb_{2}$ and
$\zb_{2}(H,\quad,...,\quad)=\zb_{1}$ [which implies
$\zb_{j}(H,\quad,...,\quad)=\zb_{j}(\quad,H,...,\quad)$

\noindent $=\zb_{j}(\quad,\quad,...,H)$, $j=1,2$], and conversely. Besides
$\zb$ is holomorphic if and only if regarded from the real
viewpoint $\zb_{1}(X_{1},...,X_{\widetilde k})
+\zima\zb_{2}(X_{1},...,X_{\widetilde k})$ is a holomorphic function
whatever $X_{1},...,X_{\widetilde k}$ are holomorphic vector
fields. In particular if $\zb$ is a complex ${\widetilde k}$-form
and its real part is closed then $\zb$ is holomorphic and closed.

Finally, a holomorphic $(1,1)$-tensor field regarded from the real
point of view is just a real $(1,1)$-tensor field that commutes
with $H$ and transforms holomorphic vector fields into holomorphic
vector fields.

Until the end of this section one works in the analytic category,
in which theorem 5.1 will be deduced from the complex case of
theorems 2.1 and 3.1. We start constructing a complex structure
along ${\mathcal A}$. Shrinking $P$ if necessary, one may suppose
that the algebraic type of $\zlma_{\zbv\mathcal A}$ and that of
$\zlma_{\zbv{\mathcal A}\zin Ker(df\zex dg)}$ are constant. Set
$H_{0}=(4g-f^{2})^{-1/2}(2\zlma_{\zbv\mathcal A}+fI)$; then
$(H_{0}^{2}+I)^{m}=0$. Therefore $H_{0}$ is $0$-deformable. Let
$H$ be its semi-simple part; by construction $H^{2}=-I$ and there
is a polynomial $\zq(t)$ with real coefficients such that
$H=\zq(H_{0})$, so $H={\tilde\zq}(\zlma_{\zbv{\mathcal A}})$ for some
polynomial ${\tilde\zq}\zpe {\mathbb R}_{P}[t]$. From section 6
of \cite{TU} applied on each leaf of ${\mathcal A}$ follows $N_{H}=0$;
in other words $H$ is a complex structure along ${\mathcal A}$. Moreover
$\zw(H,\quad)=\zw(\quad,H)$ and
$\zw_{1}(H,\quad)=\zw_{1}(\quad,H)$ since
$H={\tilde\zq}(\zlma_{\zbv{\mathcal A}})$.

Thus the $2$-forms $\zW=\zw+\zima{\widetilde \zw}$ and
$\zW_{1}=\zw_{1}+\zima{\widetilde \zw}_{1}$, where ${\widetilde
\zw}(X,Y)=-\zw(HX,Y)$ and ${\widetilde
\zw}_{1}(X,Y)=-\zw_{1}(HX,Y)$, are holomorphic and closed because
$d\zw=d\zw_{1}=0$. Observe that
$\zW$ and $\zW_{1}$ are symplectic,
$\zlma_{\zbv\mathcal A}\zci
H=H\zci \zlma_{\zbv\mathcal A}$ and
$\zW_{1}=\zW(\zlma_{\zbv\mathcal A},\quad)$,
so $m$ is even and
$\zlma_{\zbv\mathcal A}$ is holomorphic along ${\mathcal A}$.

As $H$ is the semi-simple part of $H_0$ the tensor field $H_{0}-H$
is nilpotent and commutes with $H$; therefore $(H_{0}-H)^{m}=0$.
Hence $(\zlma_{\zbv\mathcal A}
+{\frac{1}{2}}[fI-(4g-f^{2})^{\frac{1}{2}}H])^{m}=0$. In other
words the complex polynomial $(t+h)^{m}$ where
$h={\frac{1}{2}}[f-\zima (4g-f^{2})^{\frac{1}{2}}]$ annuls
$\zlma_{\zbv\mathcal A}$, which implies that $(t+h)^{m}$ is the
complex characteristic polynomial of the holomorphic tensor
$\zlma_{\zbv\mathcal A}$. In particular $h$ is holomorphic along
$\mathcal A$ and $Kerdh=Kerdf\zin Kerdg=Kerdg_{1}\zin Kerdg_{2}$.

Shrinking $P$ allows us to suppose that the elementary divisors of
$\zlma_{\zbv\mathcal A}$ on this manifold are
$(t^{2}+ft+g)^{a_{1}}$,..., $(t^{2}+ft+g)^{a_{\widetilde k}}$.
Consider any point $q\zpe P$ and a cyclic decomposition ${\mathcal
A}(q)={\mathcal B}\zdi...\zdi{\mathcal B}_{\widetilde k}$
associated to these elementary divisors. Then $H{\mathcal
B}_{j}={\mathcal B}_{j}$ since $\zlma{\mathcal B}_{j}\zco{\mathcal
B}_{j}$; that is each ${\mathcal B}_{j}$ is a cyclic complex
vector subspace. Now reasoning as before on every ${\mathcal
B}_{j}$ at each $q\zpe P$ shows that $(t+h)^{a_{1}}$,...,
$(t+h)^{a_{\widetilde k}}$ are the complex elementary divisors of
$\zlma_{\zbv\mathcal A}$. By the same reason if
$(t^{2}+ft+g)^{b_{1}}$,..., $(t^{2}+ft+g)^{b_{k'}}$ are the
elementary divisors of $\zlma_{\zbv{\mathcal A}\zin Kerdh}$ then
$(t+h)^{b_{1}}$,..., $(t+h)^{b_{k'}}$ are the complex elementary
divisors of $\zlma_{\zbv{\mathcal A}\zin Kerdh}$.

The analytic complex Frobenius theorem yields functions
$z_{1},...,z_{2m}$ such that $(x_{1},...,x_{n},z_{1},...,z_{2m})$
is a system of coordinates around $p$ and

\noindent $H=\zsu_{k=1}^{2m}((\zpar/\zpar z_{2k})\zte dz_{2k-1}-(\zpar/\zpar
z_{2k-1})\zte dz_{2k})_{\zbv\mathcal A}$. Now we may consider the
complex coordinates $v_{1}=z_{1}+\zima z_{2}$,...,
$v_{m}=z_{2m-1}+\zima z_{2m}$ and, after shrinking, identify $P$
to a product open set $A\zpor B\zco{\mathbb R}^{n}\zpor{\mathbb
C}^{m}$. In complex notation $H$ equals $\zima I$ on ${\mathcal
A}$.

Moreover functions $z_{1},...,z_{2m}$ can be choose in such a way
that $\zW=(dv_{1}\zex dv_{2}+...+dv_{m-1}\zex
dv_{m})_{\zbv\mathcal A}$, and $h=v_{m}$ if $df_{\zbv\mathcal
A(p)}\znoi 0$. Indeed, consider complex variables
$u_{1}=x_{1}+\zima y_{1}$,..., $u_{n}=x_{n}+\zima y_{n}$ on an
open set $A'\zco {\mathbb C}^{n}$ such that $A'\zin{\mathbb
R}^{n}=A$ and by means of the analyticity extend $\mathcal A$,
$\zW$ and $h$, in the obvious way, to an open neighbourhood of
$(p,0)$ on $A'\zpor B$. Then apply the Darboux theorem for
obtaining suitable coordinates $(u_{1},...,u_{n},{\widetilde
v}_{1},...,{\widetilde v}_{m})$ and, finally, restraint functions
${\widetilde v}_{1},...,{\widetilde v}_{m}$ to $P=A\zpor B$.

On the other hand $G=\zsu_{j=1}^{n}a_{j}(\zpar/\zpar x_{j})\zte
dx_{j} +\zsu_{j,k=1}^{2m}g_{jk}(\zpar/\zpar z_{j})\zte dz_{k}
+\zsu_{j=1}^{n}X_{j}\zte dx_{j}$ where
$X_{1},...,X_{n}\zpe{\mathcal A}$.

After changing the other of variables $x_{1},...,x_{n}$ and
shrinking $A$, we may assume that $dx_{1}\zex...\zex
dx_{n-r}\zex\za_{1}\zex...\zex\za_{r}$ has no zeros. Set
$T=\zsu_{j=1}^{n}X_{j}\zte dx_{j}$. Since it is enough proving
theorem 5.1 for some $G+\zsu_{k=1}^{r}Y_{k}\zte \za_{k}$ where
$Y_{1},...,Y_{r}\zpe{\mathcal A}$, one can suppose
$X_{n-r+1}=...=X_{n}=0$. Thus there exist vector fields
$X'_{1},...,X'_{n}\zpe{\mathcal F}$ functional combinations of
$\zpar/\zpar x_{1},...,\zpar/\zpar x_{n}$ whose coefficients do
not depend on $z$ such that $TX'_{j}=X_{j}$, $j=1,...,n$. Observe
that $L_{X'_{j}}H=0$.

But $N_{G}({\mathcal F},{\mathcal F})=0$ so
$N_{G}(X'_{j},{\mathcal A})=0$, whence $G\zci
L_{X'_{j}}(G_{\zbv\mathcal A})-
L_{(JX'_{j}+X_{j})}(G_{\zbv\mathcal A})=0$ that is
$L_{X_{j}}(G_{\zbv\mathcal A})=G\zci L_{X'_{j}}(G_{\zbv\mathcal
A})-L_{JX'_{j}}(G_{\zbv\mathcal A})$. As $JX'_{j}$ is a functional
combination of $\zpar/\zpar x_{1},...,\zpar/\zpar x_{n}$ with
coefficients only depending on $x$ one has $L_{JX'_{j}}H=0$. By
construction $G_{\zbv\mathcal A}$ and $H$ commute, therefore
$L_{X'_{j}}(G_{\zbv\mathcal A})$ and $L_{JX'_{j}}(G_{\zbv\mathcal
A})$ commute with $H$ and from the expression above follows that
$L_{X_{j}}(G_{\zbv\mathcal A})$ does too.

On the other hand $H={\widetilde \zq}(G_{\zbv\mathcal A})$ so
$L_{X_{j}}H$ equals a polynomial in $G_{\zbv\mathcal A}$ and
$L_{X_{j}}(G_{\zbv\mathcal A})$, which implies that $H\zci
L_{X_{j}}H=(L_{X_{j}}H)\zci H$. In turn from $H^{2}=-I$ follows
$H\zci L_{X_{j}}H=-(L_{X_{j}}H)\zci H$, therefore $ L_{X_{j}}H=0$
since $H$ is invertible. In short, we may assume $X_{1},...,X_{n}$
holomorphic without loss of generality. Now in complex notation
one has $G=\zsu_{j=1}^{n}a_{j}(\zpar/\zpar x_{j})\zte dx_{j}
+\zsu_{j,k=1}^{m}h_{jk}(\zpar/\zpar v_{j})\zte dv_{k}
+\zsu_{j=1}^{n}(\zsu_{k=1}^{m}h'_{jk}\zpar/\zpar v_{k})\zte
dx_{j}$ where $h_{jk}$ and $h'_{jk}$ are holomorphic along
${\mathcal A}$.

If as before we consider complex variables $u_{1}=x_{1}+\zima
y_{1}$,..., $u_{n}=x_{n}+\zima y_{n}$ on an open set $A'\zco
{\mathbb C}^{n}$ such that $A'\zin{\mathbb R}^{n}=A$, then through
the analyticity $\za_{1},...,\za_{r}$, ${\mathcal F}$, $\zlma$,
${\mathcal A}$, $\zW$, $\zW_{1}$ and $G$ may be extended  to
similar holomorphic objects
$\widetilde\za_{1},...,\widetilde\za_{r}$, ${\widetilde{\mathcal
F}}$, $\widetilde\zlma$, ${\widetilde{\mathcal A}}$,
$\widetilde\zW$, ${\widetilde\zW}_{1}$ and ${\widetilde G}$, which
are defined on an open neighbourhood $\widetilde P$ of
${\widetilde p}=(p,0)$ on $A'\zpor B$. In particular
$\widetilde\za_{1},...,\widetilde\za_{r}$ defines
${\widetilde{\mathcal F}}$, $du_{1}=...=du_{n}=0$ defines
${\widetilde{\mathcal A}}$ and ${\widetilde\zW}=(dv_{1}\zex
dv_{2}+...+dv_{m-1}\zex dv_{m})_{\zbv{\widetilde{\mathcal A}}}$.
After shrinking $\widetilde P$ if necessary, it is easily checked
that $({\widetilde{\mathcal
F}},{\widetilde\zlma},{\widetilde\zW},{\widetilde\zW}_{1})$
verifies the two first conditions of Veronese flag.

Observe that given a holomorphic function
$\zm=\zm_{1}+\zima\zm_{2}$ then its ${\widetilde\zW}$-hamiltonian
$X_{\zm}$ equals the $\zw_{\mathbb R}$-hamiltonian of $\zm_{1}$
where $\zw_{\mathbb R}$ is the real part of ${\widetilde\zW}$
(note that $\zw_{\mathbb R}=\zw$ on $P$). On the other hand if
$\zm$ is defined around $\widetilde p$ and
${\widetilde\zlma}^{*}d\zm$ is closed on ${\widetilde{\mathcal
F}}$, that is
${\widetilde\za}_{1}\zex...\zex{\widetilde\za}_{r}\zex d(d\zm\zci
{\widetilde G})=0$, then $\za_{1}\zex...\zex\za_{r}\zex
d(d\zm_{1}\zci G)=0$ on $P$ which implies
$\za_{1}\zex...\zex\za_{r}\zex L_{X_{\zm}}G=0$ around $p$ on $P$;
therefore ${\widetilde\za}_{1}\zex...\zex{\widetilde\za}_{r}\zex
L_{X_{\zm}}{\widetilde G}=0$ since this last tensor field is the
extension of $\za_{1}\zex...\zex\za_{r}\zex L_{X_{\zm}}G$. Hence
$L_{X_{\zm}}{\widetilde \zlma}=0$ near $\widetilde p$; in other
words $({\widetilde{\mathcal
F}},{\widetilde\zlma},{\widetilde\zW},{\widetilde\zW}_{1})$ is a
Veronese flag at $\widetilde p$.

Set ${\widetilde G}=\zsu_{j=1}^{n}a_{j}(\zpar/\zpar u_{j})\zte
du_{j} +\zsu_{j,k=1}^{m}{\widetilde h}_{jk}(\zpar/\zpar v_{j})\zte
dv_{k} +\zsu_{j=1}^{n}{\widetilde X}_{j}\zte du_{j}$. Then
${\widetilde h}_{jk}$ is the prolongation to $\widetilde P$ of
$h_{jk}$ and ${\widetilde X}_{j}$ that of $X_{j}$. Consider the
$m\zpor m$ matrix  ${\mathcal M}=(h_{jk})+hI$ and its prolongation
${\widetilde{\mathcal M}}=({\widetilde h}_{jk})+{\widetilde h}I$
where ${\widetilde h}$ is the prolongation of $h$, which equals
$v_{m}$ when $h$ is not constant. Recall that $(t+h)^{a_{1}}$,...,
$(t+h)^{a_{\widetilde k}}$ are the complex elementary divisors of
$\zlma_{\zbv\mathcal A}$ on $P$; since these elementary divisors
up to change of order are determined by
$dimKer(\zlma_{\zbv\mathcal A}+hI)^{a}$, $a=1,...,m$, these
dimensions have to be constant, that is to say each function
$rank{\mathcal M}^{a}$, $a=1,...,m$, is constant on $P$. This fact
implies that $rank{\widetilde{\mathcal M}}^{a}$, $a=1,...,m$, is
constant near $\widetilde p$ on $\widetilde P$ and equals
$rank{\mathcal M}^{a}$.

Indeed, as ${\widetilde{\mathcal M}}({\widetilde p})={\mathcal
M}(p)$ one has $rank{\widetilde{\mathcal M}}^{a}\zmai
rank{\mathcal M}^{a}$. Let $\widetilde\zr$ be a minor of
${\widetilde{\mathcal M}}^{a}$  and $\zr$ the similar one of
${\mathcal M}^{a}$. Then $\widetilde\zr$ is the prolongation of
$\zr$, so $\widetilde\zr$ vanishes on $\widetilde P$ if $\zr$ does
on $P$. Therefore $rank{\widetilde{\mathcal M}}^{a}\zmei
rank{\mathcal M}^{a}$. Thus
$dimKer({\widetilde\zlma}_{\zbv\widetilde{\mathcal A}}+{\widetilde
h}I)^{a}=dimKer(\zlma_{\zbv\mathcal A}+hI)^{a}$, $a=1,...,m$, and
consequently $(t+{\widetilde h})^{a_{1}}$,..., $(t+{\widetilde
h})^{a_{\widetilde k}}$ are the elementary divisors of
${\widetilde\zlma}_{\zbv\widetilde{\mathcal A}}$ closed to
$\widetilde p$.

A similar argument shows that $(t+{\widetilde h})^{b_{1}}$,...,
$(t+{\widetilde h})^{b_{k'}}$ are the elementary divisors of
${\widetilde\zlma}_{\zbv{\widetilde{\mathcal A}}\zin Kerdh}$
closed to $\widetilde p$. In short, the point $\widetilde p$ is
regular for ${\widetilde\zlma}_{\zbv\widetilde{\mathcal A}}$.

When $h$ is not constant, in coordinates $(u,v)$ one has
${\widetilde h}=v_{m}$ and ${\widetilde\zW}=(dv_{1}\zex
dv_{2}+...+dv_{m-1}\zex dv_{m})_{\zbv{\widetilde{\mathcal A}}}$;
in particular $\zw=\zsu_{j=1}^{m/2}(dz_{4j-3}\zex dz_{4j-1}
-dz_{4j-2}\zex dz_{4j})_{\zbv\mathcal A}$. Thus the symplectic
reduction of $({\widetilde{\mathcal
F}},{\widetilde\zlma},{\widetilde\zW},{\widetilde\zW}_{1})$ can be
identified to the extension by means of complex variables
$u_{1}=x_{1}+\zima y_{1}$,..., $u_{n}=x_{n}+\zima y_{n}$ of the
symplectic reduction of $({\mathcal F},\zlma,\zw,\zw_{1})$.
Indeed, $g_{1},g_{2}$ are function of $z_{2m-1},z_{2m}$ only so
$E$ is spanned by $\zpar/\zpar z_{2m-3},\zpar/\zpar z_{2m-2}$ and
from the complex viewpoint it is the foliation spanned by
$\zpar/\zpar v_{m-1}$. Since the symplectic reduction of
$({\mathcal F},\zlma,\zw,\zw_{1})$ is a Veronese flag at
${\bar\zp}(p)$, the symplectic reduction of $({\widetilde{\mathcal
F}},{\widetilde\zlma},{\widetilde\zW},{\widetilde\zW}_{1})$ is a
Veronese flag at $({\bar\zp}(p),0)$, which is the image of
${\widetilde p}=(p,0)$ by the canonical projection [just adapt the
argument showing that $({\widetilde{\mathcal
F}},{\widetilde\zlma},{\widetilde\zW},{\widetilde\zW}_{1})$
verifies condition 3') at ${\widetilde p}$].

Since ${\widetilde h}({\widetilde p})\znope{\mathbb R}$ from
theorems 2.1 and 3.1 follow the existence around ${\widetilde p}$
of a $(1,1)$-tensor field ${\widetilde G}'$ extending
$\widetilde\zlma$ and projecting in
$\zsu_{j=1}^{n}a_{j}(\zpar/\zpar u_{j})\zte du_{j}$ and functions
$w_{1},...,w_{m}$ such that
$(u,v)=(u_{1},...,u_{n},w_{1},...,w_{m})$ is a system of
coordinates,

\centerline{${\widetilde G}'=\zsu_{j=1}^{n}a_{j}(\zpar/\zpar
u_{j})\zte du_{j} +\zsu_{j,k=1}^{m}\zh_{jk}(w)(\zpar/\zpar
w_{j})\zte dw_{k}$}

\noindent and ${\widetilde\zW},{\widetilde\zW}_{1}$ are expressed
with coefficient functions only depending on $w$ (constant when
$h$ is constant).

Finally, set $w_{1}={\widetilde z}_{1}+\zima{\widetilde
z}_{2}$,..., $w_{m}={\widetilde z}_{2m-1}+\zima{\widetilde
z}_{2m}$ and observe that ${\widetilde G}'(TP)\zco TP$; therefore
the restriction to $P$ of ${\widetilde G}'$ defines a
$(1,1)$-tensor field $G'$ which projects in $J$ and extends
$\zlma$. Now for finishing the proof of theorem 5.1 it is enough
considering $G'$ and functions ${\widetilde z}_{1},...,{\widetilde
z}_{2m}$.
\bigskip

{\bf 6. The blocks of a Veronese flag}

The aim of this section is to reduce the local study of Veronese
flags to the case where their characteristic polynomial is a power
of an irreducible one. Let $({\mathcal F},\zlma,\zw,\zw_{1})$ be a
Veronese flag on a manifold $P$ or at some point of $P$,
${\mathcal A}$ the foliation of the largest vector subspaces and
$\zp:P\zfl N$ a local quotient of $P$ by ${\mathcal A}$. Set
$codim{\mathcal F}=r$, $dim{\mathcal A}=2m$ and $dimN=n$. On $N$
consider coordinates $(x_{1},...,x_{n})$, closed $1$-forms
$\za_{1},...,\za_{r}$ and a $(1,1)$-tensor field
$J=\zsu_{j=1}^{n}(\zpar/\zpar x_{j})\zte dx_{j}$, where
$a_{1},...,a_{n}\zpe{\mathbb K}$, such that the associated
Veronese web is given by $J,\za_{1},...,\za_{r}$ and
$dx_{1}\zex...\zex dx_{n-r}\zex\za_{1}\zex...\zex\za_{r}$ never
vanishes.

Assume that on an open neighbourhood of a regular point $p$ of
$\zlma_{\zbv\mathcal A}$ the characteristic polynomial $\zf$ of
$\zlma_{\zbv\mathcal A}$ is the product of two monic relatively
prime  polynomials $\zf_{1}$ and $\zf_{2}$. Then
$Im\zf_{1}(\zlma_{\zbv\mathcal A})=Ker\zf_{2}(\zlma_{\zbv\mathcal
A})$, $Im\zf_{2}(\zlma_{\zbv\mathcal
A})=Ker\zf_{1}(\zlma_{\zbv\mathcal A})$ and ${\mathcal
A}=Im\zf_{1}(\zlma_{\zbv\mathcal A})\zdi
Im\zf_{2}(\zlma_{\zbv\mathcal A})$; moreover
$Im\zf_{1}(\zlma_{\zbv\mathcal A})$ and
$Im\zf_{2}(\zlma_{\zbv\mathcal A})$ are foliations because
$N_{\zlma}=0$ (apply lemma 2 of \cite{TU}). Thus around $p$ there exist
coordinates $(x,z,{\widetilde
z})=(x_{1},...,x_{n},z_{1},...,z_{m'},{\widetilde
z}_{1},...,{\widetilde z}_{\widetilde m})$ such that
$Im\zf_{1}(\zlma_{\zbv\mathcal A})$ is spanned by $\zpar/\zpar
z_{1},...,\zpar/\zpar z_{m'}$ and $Im\zf_{2}(\zlma_{\zbv\mathcal
A})$ is spanned by $\zpar/\zpar {\widetilde z}_{1},...,\zpar/\zpar
{\widetilde z}_{\widetilde m}$.

On the other hand $\zw_{k}(Im\zf_{1}(\zlma_{\zbv\mathcal
A}),Im\zf_{2}(\zlma_{\zbv\mathcal
A}))=\zw_{k}(Im\zf(\zlma_{\zbv\mathcal A}),\quad)=0$, $k=0,1$,
where by definition $\zw_{0}=\zw$. Hence

\centerline{$\zw_{k}=\zsu_{1\zmei i<j\zmei
m'}f_{ijk}(x,z)dz_{i}\zex dz_{j} +\zsu_{1\zmei i<j\zmei
{\widetilde m}}{\widetilde f}_{ijk}(x,{\widetilde z})d{\widetilde
z}_{i}\zex d{\widetilde z}_{j}$}

\noindent because $d\zw_{k}=0$. In particular $m'$ and $\widetilde
m$ are even since $\zw_{0}$ is symplectic.

Therefore if $G$ is a $(1,1)$-tensor field extending $\zlma$ and
projecting in $J$ one has
$G=J+\zsu_{j,k=1}^{m'}h_{jk}(x,z)(\zpar/\zpar z_{j})\zte dz_{k}
+\zsu_{j,k=1}^{\widetilde m}{\widetilde h}_{jk}(x,{\widetilde
z})(\zpar/\zpar {\widetilde z}_{j})\zte d{\widetilde z}_{k}
+\zsu_{j=1}^{n}X_{j}\zte dx_{j}$ where
$X_{1},...,X_{n}\zpe{\mathcal A}$.

Let us see that $G$ may be chosen in such a way that
$X_{1},...,X_{n}$ are foliate both for $Im\zf_{1}(G_{\zbv\mathcal
A})$ and $Im\zf_{2}(G_{\zbv\mathcal A})$. Set
$T=\zsu_{j=1}^{n}X_{j}\zte dx_{j}$. By considering
$G+\zsu_{j=1}^{r}Y_{j}\zte \za_{j}$ instead of $G$ where
$Y_{1},...,Y_{r}$ are suitable vector fields tangent to ${\mathcal
A}$, we can suppose $X_{n-r+1}=...=X_{n}=0$ without loss of
generality. Thus locally there exist
$X'_{1},...,X'_{n}\zpe{\mathcal F}$ functional combinations of
$\zpar/\zpar x_{1},...,\zpar/\zpar x_{n}$ with coefficients only
depending on $x$ such that $TX'_{j}=X_{j}$, $j=1,...,n$.

As $Im\zf_{1}(G_{\zbv\mathcal A})$ is spanned by $\zpar/\zpar
z_{1},...,\zpar/\zpar z_{m'}$ and $Ker\zf_{1}(G_{\zbv\mathcal A})$
by $\zpar/\zpar {\widetilde z}_{1}$, ..., $\zpar/\zpar {\widetilde
z}_{\widetilde m}$, the morphism $\zf_{1}(G_{\zbv\mathcal
A}):Im\zf_{1}(G_{\zbv\mathcal A})\zfl Im\zf_{1}(G_{\zbv\mathcal
A})$ is in fact an isomorphism whose inverse equals
$\zq(\zf_{1}(G_{\zbv\mathcal A}))$ for some polynomial $\zq(t)$
[indeed, if $t^{m'}+\zsu_{j=0}^{m'-1}g_{j}t^{j}$ is the
characteristic polynomial of $\zf_{1}(G_{\zbv\mathcal A})$
restricted to $Im\zf_{1}(G_{\zbv\mathcal A})$ set
$\zq(t)=-g_{0}^{-1}(t^{m'-1}+\zsu_{j=1}^{m'-1}g_{j}t^{j-1})$].
Therefore $\zsu_{j=1}^{m'}(\zpar/\zpar z_{j})\zte
{dz_{j}}_{\zbv\mathcal A}=\zr(G_{\zbv\mathcal A})$ where
$\zr(t)=\zf_{1}(t)\zpu\zq(\zf_{1}(t))$.

Analogously there is a polynomial ${\widetilde\zr}(t)$ such that
$\zsu_{k=1}^{\widetilde m}(\zpar/\zpar {\widetilde z}_{k})\zte
{d{\widetilde z}_{k}}_{\zbv\mathcal
A}={\widetilde\zr}(G_{\zbv\mathcal A})$. Set
$H=\zsu_{j=1}^{m'}(\zpar/\zpar z_{j})\zte {dz_{j}}_{\zbv\mathcal
A} -\zsu_{k=1}^{\widetilde m}(\zpar/\zpar {\widetilde z}_{k})\zte
{d{\widetilde z}_{k}}_{\zbv\mathcal A}$; then
$H=\zq_{1}(G_{\zbv\mathcal A})$ where
$\zq_{1}(t)=\zr(t)-{\widetilde\zr}(t)$. Observe that $L_{X}H=0$
for any vector field $X$ such that
$X=\zsu_{j=1}^{n}f_{j}(x)\zpar/\zpar x_{j}$.

As in section 5, from $N_{G}({\mathcal F},{\mathcal F})=0$ follows
$L_{X_{j}}(G_{\zbv\mathcal A})=G\zci L_{X'_{j}}(G_{\zbv\mathcal
A})-L_{JX'_{j}}(G_{\zbv\mathcal A})$; therefore
$L_{X_{j}}(G_{\zbv\mathcal A})$ and $H$ commute since
$L_{X'_{j}}H=L_{JX'_{j}}H=0$ and $(G_{\zbv\mathcal A})\zci
H=H\zci(G_{\zbv\mathcal A})$ because $H=\zq_{1}(G_{\zbv\mathcal
A})$. In turn $L_{X_{j}}H=L_{X_{j}}(\zq_{1}(G_{\zbv\mathcal A}))$
is a polynomial in $G_{\zbv\mathcal A}$ and
$L_{X_{j}}(G_{\zbv\mathcal A})$, which implies that $H$ and
$L_{X_{j}}H$ commute. On the other hand since $H^{2}=I$ one has
$H\zci L_{X_{j}}H=-(L_{X_{j}}H)\zci H$, so $H\zci L_{X_{j}}H=0$
and finally $L_{X_{j}}H=0$. Therefore each $X_{j}$ is foliate for
$Im\zf_{1}(G_{\zbv\mathcal A})=Im(H+I)$ and
$Im\zf_{2}(G_{\zbv\mathcal A})=Im(H-I)$, that is
$X_{j}=\zsu_{i=1}^{m'}f_{ji}(x,z)\zpar/\zpar z_{i}
+\zsu_{k=1}^{\widetilde m}{\widetilde f}_{jk}(x,{\widetilde z})\zpar/\zpar
{\widetilde z}_{k}$.

Now in variables $(x,z)$ we can consider the foliation ${\mathcal
F}'$ defined by $\za_{1},...,\za_{r}$, the $(1,1)$-tensor field

$G'=J+\zsu_{j,k=1}^{m'}h_{jk}(x,z)(\zpar/\zpar z_{j})\zte dz_{k}
+\zsu_{j=1}^{n}(\zsu_{i=1}^{m'}f_{ji}(x,z)\zpar/\zpar z_{i})\zte
dx_{j}$, and its restriction $\zlma'$ to ${\mathcal F}'$, the
$2$-forms $\zw'=\zsu_{1\zmei i<j\zmei m'}f_{ij0}(x,z)dz_{i}\zex
dz_{j}$,  $\zw_{1}'=\zsu_{1\zmei i<j\zmei
m'}f_{ij1}(x,z)dz_{i}\zex dz_{j}$
(more exactly their restriction to $Im\zf_{1}(G_{\zbv\mathcal A})$
but we omit it for simplifying the notation)
and the point $p'$ corresponding
to $p$. It is easily checked that $({\mathcal
F}',\zlma',\zw',\zw_{1}')$ is a Veronese flag, respectively a
Veronese flag at $p'$, if that was the case of $({\mathcal
F},\zlma,\zw,\zw_{1})$. Similarly in variables $(x,{\widetilde z})$
one may consider the foliation $ {\widetilde {\mathcal F}}$
defined by  $\za_{1},...,\za_{r}$, the $(1,1)$-tensor field
${\widetilde G}=J+\zsu_{j,k=1}^{\widetilde m}{\widetilde
h}_{jk}(x,{\widetilde z})(\zpar/\zpar {\widetilde z}_{j})\zte
d{\widetilde z}_{k} +\zsu_{j=1}^{n}(\zsu_{k=1}^{\widetilde m
}f_{jk}(x,{\widetilde z})\zpar/\zpar {\widetilde z}_{k})\zte
dx_{j}$, and its restriction $\widetilde\zlma$ to ${\widetilde
{\mathcal F}}$, the $2$-forms ${\widetilde\zw}=\zsu_{1\zmei
i<j\zmei {\widetilde m}}{\widetilde f}_{ij0}(x,{\widetilde
z})d{\widetilde z}_{i}\zex d{\widetilde z}_{j}$, ${\widetilde
\zw}_{1}=\zsu_{1\zmei i<j\zmei {\widetilde m}}{\widetilde
f}_{ij1}(x,{\widetilde z})d{\widetilde z}_{i}\zex d{\widetilde
z}_{j}$ and the point $\widetilde p$ corresponding to $p$; then
$({\widetilde {\mathcal F}},{\widetilde\zlma},{\widetilde
\zw},{\widetilde \zw}_{1})$ is a Veronese flag or a Veronese flag
at $\widetilde p$ if that is case of $({\mathcal
F},\zlma,\zw,\zw_{1})$.

Moreover $p'$ is regular for $\zlma'_{\zbv\mathcal A'}$ and
$\widetilde p$ for  ${\widetilde\zlma}_{\zbv{\widetilde{\mathcal
A}}}$ since $p$ was regular for $\zlma_{\zbv\mathcal A}$,
$\zf_{2}$ is the characteristic polynomial of
$\zlma'_{\zbv\mathcal A'}$ and $\zf_{1}$ that of
${\widetilde\zlma}_{\zbv{\widetilde{\mathcal A}}}$. In a more
technical way we will say that, around $p$, $({\mathcal
F},\zlma,\zw,\zw_{1})$ is the fibered product over $N$, around $p'$ and
$\widetilde p$, of $({\mathcal F}',\zlma',\zw',\zw_{1}')$ and
$({\widetilde {\mathcal F}},{\widetilde\zlma},{\widetilde
\zw},{\widetilde \zw}_{1})$.

Obviously one may reiterate the process until the characteristic
polynomial of each factor is power of an irreducible one, which
thus becomes the only case to take into account.
\bigskip

{\bf 7. The local product theorem}

In this section is showed that, around every point of some dense
open set, an analytic bihamiltonian structure decomposes into a
product of a Kronecker bihamiltonian structure and a symplectic
one if a necessary condition stated later on holds (see \cite{TU11A}).

Consider a bihamiltonian structure $(\zL,\zL_{1})$ on a real or
complex manifold $M$ of dimension $m$. The set of all $p\zpe M$
such that $rank(\zL,\zL_{1})$ is constant about $p$ is  open
(obvious) and dense. Indeed, first recall that at any $q\zpe M$
$rank((1-t)\zL+t\zL_{1})(q)=rank(\zL,\zL_{1})(q)$ except for a
finite number of scalars $t$, which is $\zmei m/2$ (see section
1.2 of \cite{TUD}). Now choose non-equal scalars $b_{1},...,b_{k}$ with
$k\zmai (m/2)+2$; then the set of all $p\zpe M$ such that the rank
of each $rank((1-b_{j})\zL+b_{j}\zL_{1})$, $j=1,...,k$, is locally
constant at $p$ is dense, open and contained in the foregoing open
set.

For simplicity sake suppose $r=corank(\zL,\zL_{1})$ locally
constant. Since our problem is local, by considering
$rank((1-b_{j})\zL+b_{j}\zL_{1})$ and $rank((1-b_{\widetilde
j})\zL+b_{\widetilde j}\zL_{1})$ for suitable indices
$j,\widetilde j$ instead of $\zL,\zL_{1}$, we may assume maximal
$(\zL,\zL_{1})$, that is
$r=corank\zL=corank\zL_{1}=corank(\zL,\zL_{1})$, without loss of
generality.

As in sub-section 1.1, for each $p\zpe M$ let ${\mathcal
A}_{1}(p)$ be the intersection of all vector subspaces
$Im(\zL+t\zL_{1})(p)$, $t\zpe{\mathbb K}$, such that
$rank(\zL+t\zL_{1})(p)=m-r$. From the algebraic model follows that
the dimension of the symplectic factor at $p$ equals
$2dim{\mathcal A}_{1}(p)+r-m$. But if $c_{1},...,c_{m}$ are
different scalars such that $rank(\zL+c_{j}\zL_{1})(p)=m-r$,
$j=1,...,m$, then ${\mathcal
A}_{1}(p)=\zin_{j=1}^{m}Im(\zL+c_{j}\zL_{1})(p)$ [see section 1.2
of \cite{TUD} again]. By continuity ${\mathcal
A}_{1}(q)=\zin_{j=1}^{m}Im(\zL+c_{j}\zL_{1})(q)$ when $q$ is close
to $p$ therefore $dim{\mathcal A}_{1}$ is a locally decreasing
function, which implies that the dimension of ${\mathcal A}_{1}$
and that of the symplectic factor are locally constant on a dense
open set.

Observe that if $(\zL,\zL_{1})$ decomposes into a product near
$p$, then the dimension of the symplectic factor has to be
constant close to $p$.

In short, suppose that on an open set $M'\zco M$ the bihamiltonian
structure is maximal and its rank and the dimension of the
symplectic factor are constant. Then, following sub-section 1.1,
set $m=2m'+2n-r$ where $2m'$ is the dimension of the symplectic
factor and consider the Veronese flag $({\mathcal
F},\zlma,\zw,\zw_{1})$ on the local quotient $P$ of $M'$ by the
secondary axis ${\mathcal A}_{2}$.

Given a (linear) symplectic form $\zt$ and a $2$-form $\zt_{1}$ on
an even dimensional vector space $V$, let $K$ be the endomorphism
of $V$ defined by $\zt_{1}=\zt(K,\quad)$. By definition the
characteristic polynomial of $(\zt,\zt_{1})$ will be that of $K$.
Let ${\widetilde \zf}=t^{2m'}+\zsu_{j=0}^{2m'-1}{\widetilde
h}_{j}t^{j}$ be the characteristic polynomial of the symplectic
factor of $(\zL,\zL_{1})$ on $M'$, that is
$t^{2m'}+\zsu_{j=0}^{2m'-1}{\widetilde h}_{j}(p)t^{j}$, for each
$p\zpe M'$, is the characteristic polynomial of the symplectic
factor of $(\zL(p),\zL_{1}(p))$ when regarded as a couple of
(linear) symplectic forms.

On the other hand let $\zf=t^{2m'}+\zsu_{j=0}^{2m'-1} h_{j}t^{j}$
be the characteristic polynomial of $\zlma_{\zbv\mathcal A}$. By
means of the algebraic model of $(\zL(p),\zL_{1}(p))$ it is not
hard to see that the symplectic factor of $(\zL(p),\zL_{1}(p))$ is
isomorphic to $(\zw(\zp_{P}(p)),\zw_{1}(\zp_{P}(p)))$. Thus the
characteristic polynomial of $(\zlma_{\zbv\mathcal
A})(\zp_{P}(p))$ equals ${\widetilde \zf}(p)$, that is locally
${\widetilde h}_{j}=h_{j}\zci\zp_{P}$, $j=0,...,2m'-1$, which in
particular shows the differentiability of ${\widetilde
h}_{0},...,{\widetilde h}_{2m'-1}$.
\bigskip

{\bf Proposition 7.1.} {\it The functions ${\widetilde
h}_{0},...,{\widetilde h}_{2m'-1}$ are in involution both for
$\zL$ and $\zL_{1}$. Moreover $\{\zL(d{\widetilde
h}_{j},\quad)(p)\}_{j=0,...,2m'-1}$ and $\{\zL_{1}(d{\widetilde
h}_{j},\quad)(p)\}_{j=0,...,2m'-1}$ span the same vector subspace
of $T_{p}M'$ for any $p\zpe M'$.}
\bigskip

{\bf Proof.} Let $\{\quad,\quad\}_{\zw}$ be the Poisson structure
on $P$ defined by $({\mathcal A},\zw)$ and
$\{\quad,\quad\}_{\zw_{1}}$ that defined by $({\mathcal
A},\zw_{1})$. Recall that $\{\quad,\quad\}_{\zw}$ is the
projection of $\zL$ and $\{\quad,\quad\}_{\zw_{1}}$ that of
$\zL_{1}$. Thus for proving the involution of ${\widetilde
h}_{0},...,{\widetilde h}_{2m'-1}$ it is enough showing that
$h_{0},...,h_{2m'-1}$ are in involution with respect to
$\{\quad,\quad\}_{\zw}$ and $\{\quad,\quad\}_{\zw_{1}}$.

On the other hand by lemma 1.2, $kdg_{k+1}=(k+1)dg_{k}\zci \zlma$
on ${\mathcal F}$ where $g_{
k}=trace(\zlma_{\zbv\mathcal A})^{ k}$, $
k\zmai 0$, whence $(k+1)\zlma X_{g_{k}}=kX_{g_{k+1}}$. Therefore
if $1\zmei k\zmei{\widetilde k}$ one has $\zw(X_{g_{\widetilde
k}},X_{g_{k}})=C(k,{\widetilde k})\zpu\zw(\zlma^{{\widetilde
k}-k}X_{g_{k}},X_{g_{k}})=0$ since $\zw(\zlma^{{\widetilde
k}-k},\quad)$ is a $2$-form on ${\mathcal F}$; so $\{g_{\widetilde
k},g_{k}\}_{\zw}=0$. But $h_{0},...,h_{2m'-1}$ are function of
$g_{1},...,g_{k},...$ (see section 3) therefore
$\{h_{i},h_{j}\}_{\zw}=0$.

As it was pointed out in sub-section 1.1, from the algebraic model
follows that
$\zL_{1}(\zp_{P}^{*}\zlma^{*}\zb,\quad)=\zL(\zp_{P}^{*}\zb,\quad)$
for any $\zb\zpe T_{{\zp_{P}}(q)}^{*}P$ and $q\zpe M'$. Thus in
our case $kdg_{k+1}=(k+1)dg_{k}\zci \zlma$ on ${\mathcal F}$
implies $\zL(d(g_{k}\zci\zp_{P}),\quad)
=k(k+1)^{-1}\zL_{1}(d(g_{k+1}\zci\zp_{P}),\quad)$. Since
$h_{0},...,h_{2m'-1}$ are function of $g_{1},...,g_{k},...$ and
the traces are function of $h_{0},...,h_{2m'-1}$ (see section 3
again), the same thing happens with ${\widetilde
h}_{0},...,{\widetilde h}_{2m'-1}$ and
$g_{1}\zci\zp_{P},...,g_{k}\zci\zp_{P},...$ Therefore the vector
subspace spanned by $\{\zL(d{\widetilde
h}_{j},\quad)(p)\}_{j=0,...,2m'-1}$ is contained in that spanned
by $\{\zL_{1}(d{\widetilde h}_{j},\quad)(p)\}_{j=0,...,2m'-1}$.

For finishing the proof it is enough inverting the roles of $\zL$
and $\zL_{1}$ because the characteristic polynomial of the
symplectic factor of $(\zL_{1},\zL)$ equals
$t^{2m'}+\zsu_{j=1}^{2m'-1}{\widetilde h}_{2m'-j}{\widetilde
h}_{0}^{-1}t^{j}+{\widetilde h}_{0}^{-1}$. $\square$

Now assume that $(M',\zL,\zL_{1})$ is diffeomorphic to a product
of a Kronecker bihamiltonian structure and a symplectic one
$(M_{1},\zL',\zL'_{1})\zpor (M_{2},\zL'',\zL''_{1})$. Let
${\mathcal B}_{1}$ and ${\mathcal B}_{2}$ be the foliations given
by the first and second factor respectively. Then ${\mathcal
A}_{1}\zcco{\mathcal B}_{2}$ and ${\widetilde
h}_{0},...,{\widetilde h}_{2m'-1}$ are ${\mathcal B}_{1}$-foliate
functions; therefore the dimension of the vector subspace of
$T_{q}^{*}M'$ spanned by $d{\widetilde h}_{0}(q),...,d{\widetilde
h}_{2m'-1}(q)$ equals the dimension of the vector subspace of
${\mathcal A}_{1}^{*}(q)$ spanned by ${d{\widetilde
h}_{0}}_{\zbv{\mathcal A}_{1}(q)},...,{d{\widetilde
h}_{2m'-1}}_{\zbv{\mathcal A}_{1}(q)}$ whenever $q\zpe M'$.

Thus {\it the foregoing property is necessary} for the existence
of a local decomposition into a product of a Kronecker
bihamiltonian structure and a symplectic one.

A point $p$ of $M$ is called {\it regular} for $(\zL,\zL_{1})$ if
the three following conditions hold:

\noindent 1) The rank $(\zL,\zL_{1})$ is constant on an open
neighbourhood $M'$ of this point.

Observe that this first condition allows assuming maximal
$(\zL,\zL_{1})$ by replacing $(\zL,\zL_{1})$ by
$(1-b)\zL+b\zL_{1}$ and $(1-b')\zL+b'\zL_{1}$, for suitable
scalars $b,b'$, and shrinking $M'$. Then:

\noindent 2) The dimension of the symplectic factor is constant
near $p$, that is on $M'$ by shrinking this neighbourhood again if
necessary.

\noindent 3) The point $\zp_{P}(p)$ is regular for
$\zlma_{\zbv\mathcal A}$.

Obviously there are many choices of scalars $b,b'$ such that
$((1-b)\zL+b\zL_{1},(1-b')\zL+b'\zL_{1})$ is maximal around $p$,
but it is easily checked that conditions 2) and 3) do not depend
on them.

Since the set of regular points of $\zlma_{\zbv\mathcal A}$ is open
and dense and the projection $\zp_{P}$ is a submersion, the set of
regular points of $(\zL,\zL_{1})$ is dense and open on $M$; it
will be named {\it the regular open set}.
\bigskip

{\bf Theorem 7.1.} {\it Consider a real analytic or holomorphic
bihamiltonian structure $(\zL,\zL_{1})$ on $M$ and a regular point
$p$. Let ${\widetilde \zf}=t^{2m'}+\zsu_{j=0}^{2m'-1}{\widetilde
h}_{j}t^{j}$ be the characteristic polynomial of the symplectic
factor of $(\zL,\zL_{1})$ near $p$. Assume that when $q$ is close
to $p$ the vector subspace spanned by $d{\widetilde
h}_{0}(q),...,d{\widetilde h}_{2m'-1}(q)$ and that spanned by
${d{\widetilde h}_{0}}_{\zbv{\mathcal A}_{1}(q)},...,{d{\widetilde
h}_{2m'-1}}_{\zbv{\mathcal A}_{1}(q)}$ have the same dimension.
Then, around $p$,  $(\zL,\zL_{1})$ decomposes into a product of a
Kronecker bihamiltonian structure and a symplectic one.

Moreover, if $\zf(p)$ only has real roots then in the $C^{\zinf}$
category $(\zL,\zL_{1})$ locally decomposes into a product
Kronecker-symplectic.}
\bigskip

Let us prove theorem 7.1. Shrinking $M$ we may assume that the
hypothesis of theorem hold for every point of this manifold.
Around $\zp_{N}(p)$ consider coordinates $(x_{1},...,x_{n})$,
scalars $a_{1},...,a_{n}$, the tensor field
$J=\zsu_{j=1}^{n}a_{j}(\zpar/\zpar x_{j})\zte dx_{j}$ and closed
$1$-forms  $\za_{1},...,\za_{n}$ such that  $a_{1},...,a_{n}$ are
not roots of the characteristic polynomial $\zf(\zp_{P}(p))$ of
$(\zlma_{\zbv\mathcal A})(\zp_{P}(p))$ and $\za_{1},...,\za_{n},J$
define the Veronese web associated to the Veronese flag
$({\mathcal F},\zlma,\zw,\zw_{1})$ on $P$ induced in turn by
$(\zL,\zL_{1})$. Now shrinking $P$ allows supposing that
$a_{1},...,a_{n}$ never are roots of the characteristic polynomial
$\zf$ of $\zlma_{\zbv\mathcal A}$.

The next aim will be to show the existence near $\zp_{P}(p)$ of
functions $z_{1},...,z_{2m'}$ and a $(1,1)$-tensor field $G$
extending $\zlma$ such that
$(x,z)=(x_{1},...,x_{n},z_{1},...,z_{2m'})$ is a system of
coordinates, $G=\zsu_{j=1}^{n}a_{j}(\zpar/\zpar x_{j})\zte
dx_{j}+\zsu_{j,k=1}^{n}h_{jk}(z)(\zpar/\zpar z_{j})\zte dz_{k}$
and $\zw,\zw_{1}$ are expressed  relative to
${dz_{1}}_{\zbv\mathcal A},...,{dz_{2m'}}_{\zbv\mathcal A}$ with
coefficient functions only depending on $z$.

First, around $\zp_{P}(p)$, consider a $(1,1)$-tensor field
$G_{0}$ extending $\zlma$ and projecting in $J$. Taking into
account section 6 one may suppose $G_{0}$ adapted to the blocks of
$({\mathcal F},\zlma,\zw,\zw_{1})$, where each of them has a
characteristic polynomial power of an irreducible one. Therefore
it suffices dealing with the problem in every block $({\mathcal
F}',\zlma',\zw',\zw'_{1})$. Observe that the corresponding point
$p'$ is regular for $\zlma'_{\zbv\mathcal A'}$.

On the other hand each ${\widetilde h}_{j}=h_{j}\zci\zp_{P}$ and
the foliation ${\mathcal A}_{1}$ projects in ${\mathcal A}$,
therefore at every point the vector subspace spanned by
$dh_{0},...,dh_{2m'-1}$ and that spanned by
${dh_{0}}_{\zbv\mathcal A},...,{dh_{2m'-1}}_{\zbv\mathcal A}$ have
the same dimension. Thus since $\zf$ is the product of the
characteristic polynomial of the blocks one has the following two
cases:

\noindent (I) If $(t-f)^{2m''}$ is the characteristic polynomial
of $\zlma'_{\zbv\mathcal A'}$ then $f$ is either constant or
$(df_{\zbv\mathcal A'}) (p')\znoi 0$; besides $f$ never takes the
values $a_{1},...,a_{n}$.

\noindent (II) If $(t^{2}+ft+g)^{m''}$, where ${\mathbb
K}={\mathbb R}$ and  $f^{2}<4g$, is the the characteristic
polynomial of $\zlma'_{\zbv\mathcal A'}$ then $f$ is either
constant or $(df_{\zbv\mathcal A'}) (p')\znoi 0$.

When $f$ is constant theorems 2.1 and 5.1 give us the required
coordinates and the $(1,1)$-tensor field. If $(df_{\zbv\mathcal
A'}) (p')\znoi 0$ these objects are given by theorems 3.1 and 5.1,
provided that we are able to show that the symplectic reduction is
a Veronese flag or, more exactly, to check the third condition of
this notion. Let $({\bar{\mathcal
F}'},{\bar\zlma'},{\bar\zw'},{\bar\zw'}_{1})$ be the symplectic
reduction of $({\mathcal F}',\zlma',\zw',\zw'_{1})$ and $\zp'$ its
corresponding canonical projection. Consider a function $h$ on an
open set of the symplectic reduction such that
$({\bar\zlma'})^{*}dh$ is closed along ${\bar{\mathcal F}'}$. Then
regarded as a function on an open set of $P$ in the obvious way
(that is first compose with $\zp'$ and then extend from the block
to $P$) $\zlma^{*}dh$ is closed along the foliation ${\mathcal
F}\zin Kerdf$ or ${\mathcal F}\zin Kerdf\zin Kerdg$. By lemma 1.6,
at each point $L_{X_{f}}\zlma$ sends ${\mathcal F}\zin Kerdf$,
respectively ${\mathcal F}\zin Kerdf\zin Kerdg$, into the vector
space spanned by $X_{f}$, respectively $X_{f},X_{g}$.

But $X_{h}$ is tangent to the block corresponding to $({\mathcal
F}',\zlma',\zw',\zw'_{1})$ since $h$ only depends on the variables
of this block. Therefore $L_{X'_{f}}\zlma'$ sends ${\mathcal
F'}\zin Kerdf$ or ${\mathcal F'}\zin Kerdf\zin Kerdg$ into the
vector space spanned by $X'_{f}$ or $X'_{f},X'_{g}$, where
$X'_{h},X'_{f},X'_{g}$ are the $\zw'$-hamiltonians of $h,f,g$
respectively.

On the other hand $X'_{f}h=X'_{g}h=0$, therefore
$X'_{h}f=X'_{h}g=0$; that is to say $X'_{h}$ is tangent to
${\mathcal A'}\zin Kerdf$ or to ${\mathcal A'}\zin Kerdf\zin
Kerdg$. Moreover by $\zp'$ the vector field $X'_{h}$ projects in
the ${\bar\zw'}$-hamiltonian ${\bar X'}_{h}$ of $h$, whereas
${\zlma'}_{\zbv{\mathcal F'}\zin Kerdf}$ or
${\zlma'}_{\zbv{\mathcal F'}\zin Kerdf\zin Kerdg}$ do in
${\bar\zlma'}$. Thus $L_{X'_{h}}\zlma'$ restricted to ${\mathcal
F'}\zin Kerdf$ or to ${\mathcal F'}\zin Kerdf\zin Kerdg$ projects
in $L_{{\bar X'}_{h}}{\bar\zlma'}$, whence $L_{{\bar
X'}_{h}}{\bar\zlma'}=0$. {\it In short, the symplectic reduction
is a Veronese flag.}

We need the following lemma whose proof is an exercise on Poisson
structures (see \cite{WE}).
\bigskip

{\bf Lemma 7.1.} {\it On a manifold ${\widetilde M}$ consider a
Poisson structure $\widetilde\zL$ and a $2\widetilde
m$-codimensional foliation ${\mathcal G}$. Assume that:

\noindent (a) The bracket of any two foliate functions is a
foliate function.

\noindent (b) The hamiltonians of the foliate functions give rise
to a $2{\widetilde m}$-dimensional vector sub-bundle
${\widetilde{\mathcal G}}$ of $T{\widetilde M}$.

Then ${\widetilde{\mathcal G}}$ is a foliation, $TM={\mathcal
G}\zdi{\widetilde{\mathcal G}}$ and, in coordinates
$(u,v)=(u_{1},...,u_{k},v_{1},...,v_{2{\widetilde m}})$ such that
${\mathcal G}$ and ${\widetilde{\mathcal G}}$ are defined by
$dv_{1}=...=dv_{2{\widetilde m}}=0$ and $du_{1}=...=du_{k}=0$
respectively, one has

\centerline{${\widetilde\zL}=\zsu_{1\zmei i<j\zmei
k}{\zh}_{ij}(u)(\zpar/\zpar u_{i})\zex (\zpar/\zpar u_{j})
+\zsu_{1\zmei i<j\zmei 2{\widetilde m} }{\widetilde
\zh}_{ij}(v)(\zpar/\zpar v_{i})\zex (\zpar/\zpar v_{j})$.}

Moreover $\zsu_{1\zmei i<j\zmei 2{\widetilde m} }{\widetilde
\zh}_{ij}(v)(\zpar/\zpar v_{i})\zex (\zpar/\zpar v_{j})$ is a
symplectic Poisson structure in variables
$(v_{1},...,v_{2{\widetilde m}})$.}
\bigskip

By means of $\zp_{P}$ functions $z_{1},...,z_{2m'}$ may be
regarded as functions defined around $p$ on $M$; since $\zp_{P}$
is a submersion ${\mathcal G}_{0}=Ker(dz_{1}\zex...\zex dz_{2m'})$
is a $2m'$-codimensional foliation about $p$. On the other hand,
since $\{z_{i},z_{j}\}_{\zw}$ and $\{z_{i},z_{j}\}_{\zw_{1}}$ are
only function of $z$ and $\zL$, $\zL_{1}$ project in the bivectors
associated to $({\mathcal A},\zw)$ and $({\mathcal A},\zw_{1})$
respectively, the functions $\zL(dh_{1},dh_{2})$ and
$\zL_{1}(dh_{1},dh_{2})$ are ${\mathcal G}_{0}$-foliate whenever
$h_{1},h_{2}$ are ${\mathcal G}_{0}$-foliate. Besides, near $p$,
the $\zL$-hamiltonians of the ${\mathcal G}_{0}$-foliate functions
give rise to a vector sub-bundle ${\mathcal G}_{1}$ of dimension
$2m'$ because $\zw$ is symplectic on ${\mathcal A}$. In the same
way, the $\zL_{1}$-hamiltonians of the ${\mathcal G}_{0}$-foliate
functions give rise to a vector sub-bundle ${\mathcal G'}_{1}$ of
dimension $2m'$. But
$\zL_{1}(\zp_{P}^{*}\zlma^{*}\zb,\quad)=\zL(\zp_{P}^{*}\zb,\quad)$,
$dz_{j}\zci G=\zsu_{k=1}^{2m'}h_{jk}(z)dz_{k}$, $j=1,...,2m'$, and
$\zlma_{\zbv\mathcal A}=G_{\zbv\mathcal A}$ is invertible;
therefore ${\mathcal G'}_{1}={\mathcal G}_{1}$.

By lemma 7.1 applied to $\zL,\zL_{1}$ and ${\mathcal G}_{0}$, the
vector sub-bundle ${\mathcal G}_{1}$ is a foliation and locally
$TM={\mathcal G}_{0}\zdi{\mathcal G}_{1}$. Thus around $p$ there
exist functions $u_{1},...,u_{m-2m'}$ such that
$(u,z)=(u_{1},...,u_{m-2m'},z_{1},...,z_{2m'})$ is a system of
coordinates, ${\mathcal G}_{0}$ is defined by
$dz_{1}=...=dz_{2m'}=0$ and  ${\mathcal G}_{1}$ by
$du_{1}=...=du_{m-2m'}=0$. Now from lemma 7.1 follows that

$\zL=\zsu_{1\zmei i<j\zmei m-2m'}{\zh}_{ij}(u)(\zpar/\zpar
u_{i})\zex (\zpar/\zpar u_{j})$

\hskip 4.5truecm $+\zsu_{1\zmei i<j\zmei 2m'}{\widetilde
\zh}_{ij}(z)(\zpar/\zpar z_{i})\zex (\zpar/\zpar z_{j})$

$\zL_{1}=\zsu_{1\zmei i<j\zmei m-2m'}{\zh}_{1ij}(u)(\zpar/\zpar
u_{i})\zex (\zpar/\zpar u_{j})$

\hskip 4.5truecm $+\zsu_{1\zmei i<j\zmei 2m'}{\widetilde
\zh}_{1ij}(z)(\zpar/\zpar z_{i})\zex (\zpar/\zpar z_{j})$

\noindent which decomposes $(\zL,\zL_{1})$ into a product of a
Kronecker bihamiltonian structure [variables
$(u_{1},...,u_{m-2m'})$] and a symplectic one [variables
$(z_{1},...,z_{2m'})$] and {\it finishes the proof of theorem
7.1}.
\bigskip

{\bf 8. A counter-example}

In this section one will give an example, in the $C^{\zinf}$
category, of a bihamiltonian structure for which theorem 7.1
fails (see \cite{TU11B}); more exactly one will show that the partial tensor field
$\zlma$, of the associated Veronese web, cannot be extended to a
$(1,1)$-tensor field with no Nijenhuis torsion. In our example the
bihamiltonian structured considered defines a $G$-structure and
the Lewy's result \cite{LW}  prevents us to find an extension of $\zlma$
with vanishing Nijenhuis torsion, which clearly contradicts
theorem 7.1 (the reader interested in a classic example of
non-equivalent $G$-structures may see \cite{GS}).

First let us establish some auxiliary results. Consider on a
manifold $P$ endowed with coordinates
$(x,y)=(x_{1},...,x_{n},y_{1},...,y_{m})$, for example on an open
set of ${\mathbb K}^{n+m}$, the foliation ${\mathcal A}$ given by
$dx_{1}=...=dx_{n}=0$, and in coordinates $x=(x_{1},...,x_{n})$,
that is on the quotient of $P$ by ${\mathcal A}$, a Veronese web
defined by $J=\zsu_{j=1}^{n}a_{j}(\zpar/\zpar x_{j})\zte dx_{j}$
where $a_{1},...,a_{n}\zpe{\mathbb K}-\{0\}$ and the closed
$1$-forms $\za_{1},...,\za_{r}$. Recall that in this case
$\za_{1}\zex...\zex\za_{r}\zex d(\za_{j}\zci J)=0$, $j=1,...,r$,
and $\za_{1},...,\za_{r},J^{*}$ span, at each point, the same
vector space that $dx_{1},...,dx_{n}$. In the obvious way
$J,\za_{1},...,\za_{r}$ will be regarded as objects on $P$ too. On
the other hand, assume that the $n$-form $dx_{1}\zex...\zex
dx_{n-r}\zex\za_{1}\zex...\zex\za_{r}$ never vanishes.

On $P$ let $G=J+H+\zsu_{j=1}^{n-r}X_{j}\zte dx_{j}$ where
$X_{1},...,X_{n-r}\zpe{\mathcal A}$,
$H=\zsu_{j,k=1}^{m}a_{jk}(y)(\zpar/\zpar y_{j})\zte dy_{k}$ and
the Nijenhuis torsion of $H_{\zbv \mathcal A}$ vanishes.
\bigskip

{\bf Lemma 8.1.} {\it One has:

\noindent (a) If $N_{G}\zex\za_{1}\zex...\zex\za_{r}=0$ then
${(L_{X_{j}}H)}_{\zbv\mathcal A}=0$, $j=1,...,n-r$.

\noindent (b) If ${(L_{X_{j}}H)}_{\zbv\mathcal A}=0$,
$j=1,...,n-r$, then $N_{G}({\mathcal A},\quad)=0$.}
\bigskip

{\bf Proof.} (a) From the formula $N_{G}(X,\quad)=L_{GX}G-GL_{X}G$
straightforward follows $N_{G}(\zpar/\zpar x_{i},{\mathcal A})=0$,
$i=n-r+1,...,n$. But $N_{G}(\zpar/\zpar
y_{k},\quad)\zex\za_{1}\zex...\zex\za_{r}=0$, $k=1,...,m$, so
$N_{G}(\zpar/\zpar y_{k},\quad)=0$ since the $1$-forms
$\za_{1},...,\za_{r}$ restricted to $dx_{1}=...=dx_{n-r}=0$ are
linearly independent everywhere. Thus $N_{G}(\zpar/\zpar
x_{j},{\mathcal A})=0$, $j=1,...,n-r$, which implies
${(L_{X_{j}}H)}_{\zbv\mathcal A}=0$.

(b) Clearly $N_{G}({\mathcal A},{\mathcal A})=0$ and
$N_{G}(\zpar/\zpar x_{i},{\mathcal A})=0$ when $i$ runs from
$n-r+1$ to $n$. On the other hand $N_{G}(\zpar/\zpar
x_{j},{\mathcal A})=(L_{X_{j}}H)({\mathcal A})$ if $j=1,...,n-r$.
$\square$
\bigskip

{\bf Lemma 8.2.} {\it Consider a tensor field
$G''=J+H+\zsu_{j=1}^{n}X_{j}\zte dx_{j}$ where
$X_{1},...,X_{n}\zpe{\mathcal A}$. If $N_{G''}=0$ then
${(L_{X_{j}}H)}_{\zbv\mathcal A}=0$, $j=1,...,n$.}
\bigskip

{\bf Proof.} Now $(L_{\zpar/\zpar x_{j}}G'')({\mathcal A})=0$
whereas $(L_{X_{j}}H)({\mathcal A})=(L_{G''(\zpar/\zpar
x_{j})}G'')({\mathcal A})=0$. $\square$

Hereafter $n=3$ and
$J=\zsu_{j=1}^{3}a_{j}(\zpar/\zpar x_{j})\zte dx_{j}$ where
$a_{1},a_{2},a_{3}$ are non-equal and non-vanishing real numbers.
Besides one will replace $m$ by $4m$, that is we will consider
coordinates $(x_{1},x_{2},x_{3},y_{1},...,y_{4m})$, and $P$ will
be an open set of ${\mathbb R}^{4m+3}$.
On the other
hand one will set $r=2$, $\za_{1}=dx_{1}-dx_{2}$ and
$\za_{2}=x_{2}dx_{2}-dx_{3}$. Then $\za_{1}$, $\za_{2}$,
$\za_{1}\zci J$ and $\za_{2}\zci J$ are closed. It is easily
checked that $\za_{1},\za_{2},J$ define a Veronese web of
codimension two in variables $x$ and
$dx_{1}\zex\za_{1}\zex\za_{2}=dx_{1}\zex dx_{2}\zex dx_{3}$.

For making calculations easy, we introduce a complex structure along
${\mathcal A}$ by means of the complex variables
$(z,u)=(z_{1},...,z_{m},u_{1},...,u_{m})$, where
$z_{1}=y_{1}+\zima y_{2}$, $u_{1}=y_{3}+\zima y_{4}$,...,
$z_{m}=y_{4m-3}+\zima y_{4m-2}$, $u_{m}=y_{4m-1}+\zima y_{4m}$.
Set $H=\zima I_{(z,u)}+\zsu_{j=1}^{m}(\zpar/\zpar z_{j})\zte
du_{j}$ where $I_{(z,u)}=\zsu_{j=1}^{m}[(\zpar/\zpar z_{j})\zte
dz_{j}+(\zpar/\zpar u_{j})\zte du_{j}]$.

From the real viewpoint $H$ is a $(1,1)$-tensor field with
constant coefficients and minimal polynomial $t(t^{2}+1)^{2}$,
whose semi-simple and nilpotent parts equal $\zima I_{(z,u)}$ and
$\zsu_{j=1}^{m}(\zpar/\zpar z_{j})\zte du_{j}$ respectively.
\bigskip

{\bf Lemma 8.3.} {\it Consider the $(1,1)$-tensor field $G'=J+H$
and a complex valued function $f(x,u)$ holomorphic along
${\mathcal A}$. If $d(df\zci G')\zex\za_{1}\zex\za_{2}=0$ then,
locally, there exists a complex valued function $g(x,z,u)$
holomorphic along ${\mathcal A}$ such that $d(dg\zci
G')\zex\za_{1}\zex\za_{2}=0$ and ${(dg\zci (H-\zima
I_{(z,u)}))}_{\zbv \mathcal A}={df}_{\zbv \mathcal A}$.}
\bigskip

Let us prove this result. First consider the basis of the
cotangent bundle, with respect to variables $x$,
$\{dx_{1},\za_{1},\za_{2}\}$ and its dual basis $X=\zpar/\zpar
x_{1}+\zpar/\zpar x_{2}+x_{2}\zpar/\zpar x_{3}$,
$X_{1}=-\zpar/\zpar x_{2}-x_{2}\zpar/\zpar x_{3}$,
$X_{2}=-\zpar/\zpar x_{3}$. Taking into account that
$dh=X(h)dx_{1}+X_{1}(h)\za_{1}+X_{2}(h)\za_{2}+\zsu_{j=1}^{m}[(\zpar
h/\zpar z_{j})dz_{j}+(\zpar h/\zpar u_{j})du_{j}]$ when $h$ is
holomorphic along ${\mathcal A}$, a calculation shows that
$d(df\zci G')\zex\za_{1}\zex\za_{2}=0$ if and only if $(JX-\zima
X)\zpu (\zpar f/\zpar u_{j})=0$, $j=1,...,m$. On the other hand
${(dg\zci (H-\zima I_{(z,u)}))}_{\zbv \mathcal A}={df}_{\zbv
\mathcal A}$ means that $g=\zsu_{j=1}^{m}z_{j}\zpar f/\zpar
u_{j}+\zf(x,u)$. Therefore we have to find a function $\zf(x,u)$
holomorphic along ${\mathcal A}$ in such a way that $d(dg\zci
G')\zex\za_{1}\zex\za_{2}=0$. But again a calculation shows that
this last condition is equivalent to the equation $(JX-\zima
X)\zpu (\zpar \zf/\zpar u_{j})=X\zpu (\zpar f/\zpar u_{j})$,
$j=1,...,m$ [observe that
$d(d(\zpar f/\zpar u_{k})\zci G')\zex\za_{1}\zex\za_{2}=0$,
$k=1,...,m$, since $d(df\zci G')\zex\za_{1}\zex\za_{2}=0$].

Now consider a function $\zq(x,u)$ such that $(JX-\zima
X)\zpu\zq=0$ and set $Y=(a_{1}-\zima)^{-1}x_{1}\zpar/\zpar
x_{1}+(a_{2}-\zima)^{-1}x_{2}\zpar/\zpar x_{2}+
((a_{2}-\zima)^{-1}+(a_{3}-\zima)^{-1})x_{3}\zpar/\zpar x_{3}$.
Then $[JX-\zima X,Y]=X$, which implies $(JX-\zima X)\zpu
h=X\zpu\zq$ where $h=Y\zpu\zq$.

Since $Y$ commutes with $\zpar/\zpar u_{j},\zpar/\zpar {\bar
z}_{j},\zpar/\zpar {\bar u}_{j}$, $j=1,...,m$, it suffices to set
$\zf=Y\zpu f$ {\it for finishing the proof of lemma 8.3}.

Now let $G=J+H+Z\zte dx_{1}$ where
$Z=\zsu_{j=1}^{m}f_{j}(x,u)\zpar/\zpar z_{j}$ and each $f_{j}$ is
holomorphic along ${\mathcal A}$. Then ${(L_{Z}H)}_{\zbv\mathcal
A}=0$ and, by lemma 8.1, one has $N_{G}({\mathcal A},\quad)=0$;
thus $N_{G}\zex\za_{1}\zex\za_{2}=0$ since there are only three
variables $x$.
\bigskip

{\bf Theorem 8.1.} {\it There exists a complex valued function $f(x)$,
$x\zpe{\mathbb R}^{3}$, such that, if one sets $f_{1}=u_{1}f$, then
the Nijenhuis torsion of the $(1,1)$-tensor field
${\widetilde G}=G+Z_{1}\zte\za_{1}+Z_{2}\zte\za_{2}$ never
vanishes around any point whatever $Z_{1},Z_{2}\zpe{\mathcal A}$.}
\bigskip

{\bf Proof.} The real characteristic polynomial $\zq$ of
${\widetilde G}$ equals $\zq_{1}\zpu\zq_{2}$
where $\zq_{1}=(t-a_{1})(t-a_{2})(t-a_{3})$ and
$\zq_{2}=(t^{2}+1)^{2m}$. Assume $N_{\widetilde G}=0$ for some
$Z_{1},Z_{2}\zpe{\mathcal A}$. Then locally ${\widetilde G}$
decompose into a product of two manifolds endowed each of them
with a $(1,1)$-tensor field whose real characteristic polynomials
are $\zq_{1}$ and $\zq_{2}$ respectively. Both factor tensor
fields are flat because the first one can be identified to $J$ and
the second one to the restriction of ${\widetilde G}$ to a leaf of
the foliation ${\mathcal A}$ and, obviously, this restriction is
flat (in fact the foliation by the second factor
equals $\mathcal A$).

Note that the complex structure along ${\mathcal A}$ is given by
the semi-simple part of ${\widetilde G}_{\zbv\mathcal A}$. So
there exist coordinates $({\widetilde x},{\widetilde
z},{\widetilde u})=({\widetilde x}_{1},{\widetilde
x}_{2},{\widetilde x}_{3},{\widetilde z}_{1},...,{\widetilde
z}_{m},{\widetilde u}_{1},...,{\widetilde u}_{m})$, where
${\tilde x}_{1}=x_{1}$, ${\tilde x}_{2}=x_{2}$,
${\tilde x}_{3}=x_{3}$, such that
$d{\widetilde x}_{1}=d{\widetilde x}_{2}=d{\widetilde x}_{3}=0$
defines ${\mathcal A}$,
$\za_{1}=d{\tilde x}_{1}-d{\tilde x}_{2}$,
$\za_{2}={\tilde x}_{2}d{\tilde x}_{2}-d{\tilde x}_{3}$,
${\widetilde z}_{1},...,{\widetilde
z}_{m},{\widetilde u}_{1},...,{\widetilde u}_{m}$ are holomorphic
along ${\mathcal A}$ and ${\widetilde G}={\widetilde
J}+{\widetilde H}$ where ${\widetilde
J}=\zsu_{j=1}^{3}a_{j}(\zpar/\zpar {\widetilde x}_{j})\zte
d{\widetilde x}_{j}$ and ${\widetilde H}=\zima I_{({\widetilde
z},{\widetilde u})}+ \zsu_{j=1}^{m}(\zpar/\zpar {\widetilde
z}_{j})\zte d{\widetilde u}_{j}$.

Clearly $d(du_{1}\zci{\widetilde G})\zex\za_{1}\zex\za_{2}=0$
[calculate it in coordinates $(x,z,u)$]. Besides in coordinates
$({\widetilde x},{\widetilde z},{\widetilde u})$  function $u_1$
does not depend on ${\widetilde z}$ since it is foliate with
respect to the foliation $Ker({({\widetilde G}-\zima I)}_{\zbv\mathcal
A})$. Therefore from lemma 8.3 applied, in coordinates
$({\widetilde x},{\widetilde z},{\widetilde u})$, to $u_1$ and
${\widetilde G}$ follows the local existence of a function $g$
holomorphic along ${\mathcal A}$ such that $d(dg\zci {\widetilde
G})\zex\za_{1}\zex\za_{2}=0$ and ${du_{1}}_{\zbv\mathcal
A}={(dg\zci({\widetilde H}-\zima I_{({\widetilde z},{\widetilde
u})}))}_{\zbv\mathcal A}$. But ${({\widetilde H}-\zima
I_{({\widetilde z},{\widetilde u})})}_{\zbv\mathcal A}={(H-\zima
I_{(z,u)})}_{\zbv\mathcal A}$ since this object is the nilpotent part
of ${\widetilde G}_{\zbv\mathcal A}$,  so ${du_{1}}_{\zbv\mathcal
A}={(dg\zci (H-\zima I_{(z,u)}))}_{\zbv\mathcal A}$; moreover
$d(dg\zci G)\zex\za_{1}\zex\za_{2}=0$ because $({\widetilde
G}-G)\zex\za_{1}\zex\za_{2}=0$. The first condition implies that
$g=z_{1}+\zr(x,u)$ where $\zr$ is holomorphic along ${\mathcal
A}$.

Now take $f_{1}=u_{1}f(x)$;
then $0=(d(dg\zci G)\zex\za_{1}\zex\za_{2})(X,X_{1},X_{2},
\zpar/\zpar u_{1})=d(dg\zci G)(X,\zpar/\zpar u_{1})=
X(dg(G\zpar/\zpar u_{1}))-\zpar/\zpar u_{1}(dg(GX))$, which
yields the equation

\centerline{(*)\hskip .5truecm  $(JX-\zima X)\zpu (\zpar\zr/\zpar
u_{1})+f=0$.}

Let ${\widetilde X}=[JX,-X]$. Then ${\widetilde
X}=(a_{3}-a_{2})\zpar/\zpar x_{3}\znoi 0$ since $a_{2}\znoi
a_{3}$. Regarded on ${\mathbb R}^{3}$, the vector fields
$JX,-X,{\widetilde X}$, which are
linearly independent everywhere, define a $3$-dimensional Lie
algebra whose center is spanned by ${\widetilde X}$.
Moreover $b_{1}JX-b_{2}X+b_{3}{\widetilde X}$ is complete for any
$b_{1},b_{2},b_{3}\zpe \mathbb R$.

On ${\mathbb R}^{3}$ endowed with coordinates
$y=(y_{1},y_{2},y_{3})$ set $Y_{1}=-\zpar/\zpar
y_{1}-2y_{2}\zpar/\zpar y_{3}$, $Y_{2}=-\zpar/\zpar
y_{2}+2y_{1}\zpar/\zpar y_{3}$ and ${\widetilde
Y}=[Y_{1},Y_{2}]=-4\zpar/\zpar y_{3}$; note that
$Y_{1},Y_{2},{\widetilde Y}$ are linearly independent
everywhere and define a $3$-dimensional Lie algebra whose center
is spanned by ${\widetilde Y}$;
moreover $b_{1}Y_{1}+b_{2}Y_{2}+b_{3}{\widetilde Y}$ is complete
for any $b_{1},b_{2},b_{3}\zpe \mathbb R$.
As ${\mathbb R}^{3}$ is simply
connected there is a diffeomorphism  of this space transforming
$JX,-X,{\widetilde X}$ in $Y_{1},Y_{2},{\widetilde Y}$
respectively.

From the Lewy's example (see (5) of page 156 of \cite{LW}) follows the
existence of a $C^{\zinf}$ function $F:{\mathbb R}^{3}\zfl{\mathbb
C}$ such that the equation $(Y_{1}+\zima Y_{2}){\widetilde F}=F$
has no solution in any neighbourhood of any point of ${\mathbb
R}^{3}$. Pulling-back $-F$ gives a function $f$ for which equation
(*) has no solution at all (regard $\zpar\zr/\zpar u_{1}$ as a
function of $x$ and $u_{1},...,u_{m}$ as parameters);
in other words if one takes $f_{1}=u_{1}f$ then
$N_{\widetilde G}$ never vanishes around any point. $\square$

The next step will be to apply the construction of sub-section 1.2 to
a foliation and a particular $(1,1)$-tensor field on ${\mathbb R}^{7}$.
More exactly, set $m=1$ and
$S=J+H+u_{1}f(x)(\zpar/\zpar z_{1})\zte dx_{1}$
or in real notation
$S=\zsu_{j=1}^{3}a_{j}(\zpar/\zpar x_{j})\zte dx_{j}+
\zsu_{j=1}^{2}[(\zpar/\zpar y_{2j})\zte dy_{2j-1}-
(\zpar/\zpar y_{2j-1})\zte dy_{2j}]+
(\zpar/\zpar y_{1})\zte dy_{3}+(\zpar/\zpar y_{2})\zte dy_{4}+
[(y_{3}g_{1}-y_{4}g_{2})(\zpar/\zpar y_{1})+
(y_{3}g_{2}+y_{4}g_{1})(\zpar/\zpar y_{2})]\zte dx_{1}$ where
$f=g_{1}+\zima g_{2}$.

Note that, as it was pointed out before, $N_{S}\zex\za_{1}\zex\za_{2}=0$.
Moreover, if $\za$ is a closed $1$-form such that
$Ker\za\zcco Ker(\za_{1}\zex\za_{2})$ then
$\za_{1}\zex\za_{2}\zex d(\za\zci S)=0$ since
$\za=\zsu_{j=1}^{3}h_{k}(x)dx_{k}$. In other words the construction of
sub-section 1.2 applies to $S$ and ${\mathcal G}= Ker(\za_{1}\zex\za_{2})$.

Let $(x,y,{\tilde x},{\tilde y})=(x_{1},x_{2},x_{3},y_{1},...,y_{4},
{\tilde x}_{1},{\tilde x}_{2},{\tilde x}_{3},{\tilde y}_{1},
...,{\tilde y}_{4})$ be the coordinates of $T^{*}{\mathbb R}^{7}$
associated to $(x,y)$. Then $\zw=\zsu_{j=1}^{3}d{\tilde x}_{j}\zex dx_{j}
+\zsu_{j=1}^{4}d{\tilde y}_{j}\zex dy_{j}$, $\zw_{1}=\zw(S^{*},\quad)$
and

\noindent $S^{*}=\zsu_{j=1}^{3}a_{j}[(\zpar/\zpar x_{j})\zte
dx_{j}+ (\zpar/\zpar{\tilde x}_{j})\zte d{\tilde x}_{j}]$

$+\zsu_{j=1}^{2}[(\zpar/\zpar y_{2j})\zte dy_{2j-1}- (\zpar/\zpar
y_{2j-1})\zte dy_{2j}$

$+(\zpar/\zpar {\tilde y}_{2j-1})\zte
d{\tilde y}_{2j}- (\zpar/\zpar {\tilde y}_{2j})\zte d{\tilde
y}_{2j-1}]$

$+(\zpar/\zpar y_{1})\zte dy_{3}+(\zpar/\zpar y_{2})\zte
dy_{4}+ (\zpar/\zpar {\tilde y}_{3})\zte d{\tilde
y}_{1}+(\zpar/\zpar {\tilde y}_{4})\zte d{\tilde y}_{2}$

$+[(y_{3}g_{1}-y_{4}g_{2})(\zpar/\zpar y_{1})+
(y_{3}g_{2}+y_{4}g_{1})(\zpar/\zpar y_{2})$

$-({\tilde
y}_{1}g_{1}+{\tilde y}_{2}g_{2})(\zpar/\zpar {\tilde
y}_{3})+({\tilde y}_{1}g_{2}-{\tilde y}_{2}g_{1})(\zpar/\zpar
{\tilde y}_{4})]\zte dx_{1}$

$+\zsu_{j=1}^{3}(\zpar/\zpar {\tilde
x}_{j})\zte\zb_{j}$

\noindent where $\zb_{1},\zb_{2},\zb_{3}$ are functional
combinations of
$dx_{1},dx_{2},dx_{3}$, $dy_{3},dy_{4}$,
$d{\tilde y}_{1},d{\tilde y}_{2}$.

Recall that in our case ${\mathcal G}_{0}$ is a $2$-dimensional
foliation, isotropic and symplecticly complete for $\zw$ and
$\zw_{1}$, spanned by the $\zw$-hamiltonians of $\za_{1}\zci
J^{-1},\za_{2}\zci J^{-1}$ or by the $\zw_{1}$-hamiltonians of
$\za_{1},\za_{2}$, when $\za_{1},\za_{2},\za_{1}\zci
J^{-1},\za_{2}\zci J^{-1}$ are regarded as $1$-forms on
$T^{*}{\mathbb R}^{7}$ in the obvious way.

Therefore by projection the Poisson structures ${\zL}_{\zw}$ and
${\zL}_{\zw_{1}}$ give rise to a bihamiltonian structure
$(\zL,\zL_{1})$ on the global quotient $M=T^{*}{\mathbb
R}^{7}/{\mathcal G}_{0}$ (proposition 1.4).

Since the $\zw$-hamiltonians of $\za_{1}\zci J^{-1},\za_{2}\zci
J^{-1}$ equal $-a_{1}^{-1}\zpar/\zpar{\tilde
x}_{1}+a_{2}^{-1}\zpar/\zpar{\tilde x}_{2}$,
$-a_{2}^{-1}x_{2}\zpar/\zpar{\tilde
x}_{2}+a_{3}^{-1}\zpar/\zpar{\tilde x}_{3}$, the submanifold of
$T^{*}{\mathbb R}^{7}$ defined by ${\tilde x}_{2}={\tilde
x}_{3}=0$ is transverse to ${\mathcal G}_{0}$, which allows us to
identify it with $M$ endowed with coordinates
$(x_{1},x_{2},x_{3},y_{1},...,y_{4}, {\tilde x}_{1},{\tilde
y}_{1}, ...,{\tilde y}_{4})$, whereas $\zL,\zL_{1}$ are given by
the restriction to $M$ of $\za_{1}\zci J^{-1},\za_{2}\zci
J^{-1},\zw$ and $\za_{1},\za_{2},\zw_{1}$ respectively.

In general (see the proof of proposition 1.4) $\zL+t\zL_{1}$ is
defined by the restriction to $M$ of $\za_{1}\zci
(S^{*}+tI)^{-1}=\za_{1}\zci (J+tI)^{-1}$, $\za_{2}\zci
(S^{*}+tI)^{-1}=\za_{1}\zci (J+tI)^{-1}$ and
$\zw((I+t(S^{*})^{-1})^{-1},\quad)$. Therefore
the rank of $(\zL,\zL_{1})$ equals 10, the primary axis of
$(\zL,\zL_{1})$ is the foliation $dx_{1}=dx_{2}=dx_{3}=0$ and the
secondary one the foliation spanned by $\zpar/\zpar{\tilde
x}_{1}$; in particular the dimension of the symplectic factor
is 8 everywhere and 4 that of the Kronecker factor.
Thus the global quotient of $M$ by the secondary axis is
identified, in a natural way, to the submanifold $P'$ of
$T^{*}{\mathbb R}^{7}$ defined by the equations ${\tilde
x}_{1}={\tilde x}_{2}={\tilde x}_{3}=0$ endowed with coordinates
$(x,y,{\tilde y})$, while the foliation ${\mathcal A}$ of the
Veronese flag on $P'$ induced by $(\zL,\zL_{1})$ is given by
$dx_{1}=dx_{2}=dx_{3}=0$, and the Veronese web is defined in
variables $x=(x_{1},x_{2},x_{3})$ by $J,\za_{1},\za_{2}$.

On the other hand, as $\zw_{1}=\zw(S^{*},\quad)$ and $S^{*}$
projects on $P'$ in the $(1,1)$-tensor field

\noindent $G=\zsu_{j=1}^{3}a_{j}(\zpar/\zpar x_{j})\zte dx_{j}$

$+ \zsu_{j=1}^{2}[(\zpar/\zpar y_{2j})\zte dy_{2j-1}- (\zpar/\zpar
y_{2j-1})\zte dy_{2j}$

$+(\zpar/\zpar {\tilde y}_{2j-1})\zte
d{\tilde y}_{2j}- (\zpar/\zpar {\tilde y}_{2j})\zte d{\tilde
y}_{2j-1}]$

$+(\zpar/\zpar y_{1})\zte dy_{3}+(\zpar/\zpar y_{2})\zte
dy_{4}+ (\zpar/\zpar {\tilde y}_{3})\zte d{\tilde
y}_{1}+(\zpar/\zpar {\tilde y}_{4})\zte d{\tilde y}_{2}$

$+[(y_{3}g_{1}-y_{4}g_{2})(\zpar/\zpar y_{1})+
(y_{3}g_{2}+y_{4}g_{1})(\zpar/\zpar y_{2})$

$-({\tilde y}_{1}g_{1}+{\tilde y}_{2}g_{2})(\zpar/\zpar {\tilde
y}_{3})+({\tilde y}_{1}g_{2}-{\tilde y}_{2}g_{1})(\zpar/\zpar
{\tilde y}_{4})]\zte dx_{1}$,

\noindent this last one is a prolongation of the partial
$(1,1)$-tensor field $\zlma:{\mathcal F}\zfl TP'$, which projects
in $J$.

Since $\zlma_{\zbv\mathcal A}=G_{\zbv\mathcal A}$ is $0$-deformable
because it is written with constant coefficients, the algebraic
model of the symplectic factor of $(\zL,\zL_{1})$, which is
completely determined by $\zlma_{\zbv\mathcal A}$, does not depend
on the point considered. In particular its characteristic polynomial
equals $(t^{2}+1)^{4}$ and the hypothesis of theorem 7.1 on the
coefficients of this polynomial automatically holds.

Note that the algebraic model of the Veronese web in variables $x$
does not depend on the point as in dimension three and codimension
two there is only one model. Thus the algebraic model of the
Kronecker factor is independent of the point, $(\zL,\zL_{1})$
defines a $G$-structure and $M$ is the regular open set
of $(\zL,\zL_{1})$.

Assume that, around some point $q$ of $M$, the bihamiltonian
structure $(\zL,\zL_{1})$ decomposes into a product
Kronecker-symplectic. Then considering the local quotient by the
secondary axis on each factor separately implies the existence
about of some point $p\zpe P'$ of a $(1,1)$-tensor field ${\widetilde
G}$, which prolongs $\zlma$ and projects in $J$, whose Nijenhuis
torsion vanishes. In other words, around $p$ there exist two
vector fields $Z_{1},Z_{2}\zpe{\mathcal A}$ such that $N_{\widetilde
G}=0$ where ${\widetilde G}=G+Z_{1}\zte\za_{1}+Z_{2}\zte\za_{2}$.

Now set $y_{5}={\tilde y}_{3}$, $y_{6}=-{\tilde y}_{4}$,
$y_{7}={\tilde y}_{1}$, $y_{8}=-{\tilde y}_{2}$ and consider
complex variables $z_{1}=y_{1}+\zima y_{2}$, $u_{1}=y_{3}+\zima y_{4}$,
$z_{2}=y_{5}+\zima y_{6}$ and $u_{2}=y_{7}+\zima y_{8}$. Then
$G=J+H+Z\zte dx_{1}$ where
$Z=u_{1}f\zpar/\zpar z_{1}-u_{2}f\zpar/\zpar z_{2}$ and
$f=g_{1}+\zima g_{2}$.

By theorem 8.1 one may
choose function $f$ in such a way that the Nijenhuis torsion of ${\widetilde G}$
never vanishes about any point, which implies that $(\zL,\zL_{1})$ does
not decompose into a product Kronecker-symplectic around any point.

In short, {\it one has constructed a counter-example to theorem 7.1
in the $C^{\zinf}$ case}.
\bigskip

{\bf Appendix: A splitting property for $(1,1)$-tensor fields}

Nowadays it is well known, and belongs to the mathematical folklore,
that a $(1,1)$-tensor  fields whose Nijenhuis torsion vanishes locally follows
the decomposition of its characteristic polynomial (see \cite{FN}).
Nevertheless, and for making our text more self-contained, we will
prove this result here. More exactly:
\bigskip

{\bf Proposition A.1.} {\it Consider a $(1,1)$-tensor fields $G$ on
a $n$-manifold $M$. Let $\zf$ be its characteristic polynomial.
Assume that:

\noindent (1) $N_{G}=0$,

\noindent (2) $\zf=\zf_{1}\zpu\zf_{2}$ where $\zf_{1},\zf_{2}$ are monic
polynomials, of respective degrees $n_{1}$ and $n_{2}$, relatively prime
at each point.

Then, around every point, $(M,G)$ decomposes into a product
$(M_{1},G_{1})\zpor(M_{2},G_{2})$, where $dimM_{1}=n_{1}$,
$dimM_{2}=n_{2}$, $N_{G_{1}}=N_{G_{2}}=0$, $\zf_{1}$ is the characteristic
polynomial of $G_{1}$ (more exactly $\zf_{1}$ is the pull-back of
the characteristic polynomial of $G_{1}$
by the first projection) and $\zf_{2}$ that of $G_{2}$.}
\bigskip

Let us prove proposition A.1. Set $H_{1}=\zf_{2}(G)$ and
$H_{2}=\zf_{1}(G)$. By algebraic reasons $KerH_{1}=ImH_{2}$,
$KerH_{2}=ImH_{1}$, $ImH_{1}$ and $ImH_{2}$ are vector sub-bundles
of dimension $n_{1}$ and $n_{2}$ respectively and $TM=ImH_{1}\zdi
ImH_{2}$. Moreover $ImH_{1}$ and $ImH_{2}$ are $G$-invariant,
$\zf_{1}(H_{1})=\zf_{2}(H_{2})=0$, and
$H_{1},\zf_{2}(H_{1}):ImH_{1}\zfl ImH_{1}$,
$H_{2},\zf_{1}(H_{2}):ImH_{2}\zfl ImH_{2}$ are isomorphisms.

Since $N_{G}=0$ one has
$(L_{G^{k}X}(G^{r}))Y=(G^{k}L_{X}(G^{r}))Y$ for any vector fields
$X,Y$ and natural numbers $k,r$. Recall that
if ${\widetilde H}$ is a $(1,1)$-tensor
field then $L_{fZ}{\widetilde
H}=fL_{Z}{\widetilde H}+({\widetilde H}Z)\zte df-Z\zte
(df\zci{\widetilde H})$. Now a straightforward calculation shows:
\bigskip

{\bf Lemma A.1.} {\it Consider functions $h_{0},...,h_{s}$ and set
$H=\zsu_{k=0}^{s}h_{k}G^{k}$. Then
$N_{H}(X,Y)=\zsu_{j=0}^{n-1}[\za_{j}(X)G^{j}Y-\za_{j}(Y)G^{j}X]$
where each $\za_{j}$ is a $1$-form functional combination of
$dh_{k}\zci G^{r}$, $k=0,...,s$, $r=0,...,n-1$.

In particular $N_{H}=0$ if $h_{0},...,h_{s}$ are constant.}
\bigskip

By definition of Nijenhuis torsion
$[H_{1}X,H_{1}Y]-N_{H_{1}}(X,Y)$ is a section of $ImH_{1}$.
Therefore given vector fields $X,Y\zpe ImH_{1}$, since
$G^{k}(ImH_{1})\zco ImH_{1}$, from lemma A.1 follows that
$[H_{1}X,H_{1}Y]\zpe ImH_{1}$. But the vector fields $H_{1}Z$ such
that $Z\zpe ImH_{1}$ span $ImH_{1}$, so $ImH_{1}$ is involutive;
in turn and by a similar reason $ImH_{2}$ is involutive too.

In other words, locally, $M$ can be regarded as a product
$M_{1}\zpor M_{2}$ associated to the decomposition of the tangent
bundle $TM=ImH_{1}\zdi ImH_{2}$; moreover $G(TM_{1}\zpor
\{0\})\zco TM_{1}\zpor \{0\}$ and $G(\{0\}\zpor TM_{2})\zco
\{0\}\zpor TM_{2}$. Thus there exist two $(1,1)$-tensor field
$G_{1}:TM_{1}\zfl TM_{1}$, perhaps depending on $M_{2}$, and
$G_{2}:TM_{2}\zfl TM_{2}$, perhaps depending on $M_{1}$, such that
$G=G_{1}+G_{2}$ when $G_{1},G_{2}$ are considered on $TM$ in the
natural way (that is $G_{1}(\{0\}\zpor TM_{2})=0$ and
$G_{2}(TM_{1}\zpor \{0\})=0$). The proof will be finished if we
are able to show that $G_{1},G_{2}$ do not depend on $M_{2}$ and
$M_{1}$ respectively, since in this case $N_{G}=0$ obviously implies
$N_{G_{1}}=N_{G_{2}}=0$.

We start dealing with the case where there exist a symplectic form $\zw$
and a closed $2$-form $\zw_{1}$ such that $\zw_{1}=\zw(G,\quad)$;
recall that $\zw(G,\quad)=\zw(\quad, G)$. Then
$\zw(ImH_{1},ImH_{2})=\zw(Im(\zf_{2}(G)),Im(\zf_{1}(G)))=
\zw(Im(\zf_{1}(G)\zci \zf_{2}(G)),TM)=0$; in a analogous way one
has $\zw_{1}(ImH_{1},ImH_{2})=0$. Now consider coordinates
$(x,y)=(x_{1},...,x_{n_{1}},y_{1},...,y_{n_{2}})$ on $M$ such that
$\zpar/\zpar x_{1},...,\zpar/\zpar x_{n_{1}}$ span $ImH_{1}$ and
$\zpar/\zpar y_{1},...,\zpar/\zpar y_{n_{2}}$ span $ImH_{2}$. Then
$\zw=\zw'+\zw''$ and $\zw_{1}=\zw'_{1}+\zw''_{1}$ where
$\zw'=\zsu_{1\zmei i<j\zmei n_{1}}f_{ij}(x)dx_{i}\zex dx_{j}$,
$\zw''=\zsu_{1\zmei i<j\zmei n_{2}}g_{ij}(y)dy_{i}\zex dy_{j}$,
$\zw'_{1}=\zsu_{1\zmei i<j\zmei n_{1}}{\tilde f}_{ij}(x)dx_{i}\zex
dx_{j}$ and $\zw''_{1}=\zsu_{1\zmei i<j\zmei n_{2}}{\tilde
g}_{ij}(y)dy_{i}\zex dy_{j}$, because $d\zw=d\zw_{1}=0$ and
$\zw(\zpar/\zpar x_{k},\zpar/\zpar y_{r})=\zw_{1}(\zpar/\zpar
x_{k},\zpar/\zpar y_{r})=0$, $k=1,...,n_{1}$, $r=1,...,n_{2}$.
Thus $\zw'_{1}=\zw'(G_{1},\quad)$ in coordinates
$(x_{1},...,x_{n_{1}})$ regarded on $M_{1}$ and
$\zw''_{1}=\zw''(G_{2},\quad)$ in coordinates
$(y_{1},...,y_{n_{2}})$ on $M_{2}$; whereby $G_{1}$ only depends
on $(x_{1},...,x_{n_{1}})$ and $G_{2}$ on $(y_{1},...,y_{n_{2}})$,
which proves proposition A.1 in this case.

In the general case consider the prolongation $G^{*}$ of $G$ to
$T^{*}M$ (see sub-section 1.2) whose characteristic polynomial
equals $\zf^{2}$, or more exactly the pull-back of $\zf^{2}$ by the
canonical projection $\zp:T^{*}M\zfl M$. Now $N_{G^{*}}=0$,
$\zf^{2}=\zf^{2}_{1}\zpu\zf^{2}_{2}$ and, since on $T^{*}M$ there
exist $\zw$ and $\zw_{1}$ as before, $G^{*}$ decomposes into
a sum $G^{*}=G^{*}_{1}+G^{*}_{2}$ in such a way that
$ImG^{*}_{1}=Im\zf^{2}_{2}(G^{*})$,
$KerG^{*}_{1}=Im\zf^{2}_{1}(G^{*})$,
$ImG^{*}_{2}=Im\zf^{2}_{1}(G^{*})$ and
$KerG^{*}_{2}=Im\zf^{2}_{2}(G^{*})$. Moreover
$N_{G^{*}_{1}}=N_{G^{*}_{2}}=0$ as $N_{G^{*}}=0$.

Again, consider coordinates
$(x,y)=(x_{1},...,x_{n_{1}},y_{1},...,y_{n_{2}})$ on $M$
such that $\zpar/\zpar x_{1},...,\zpar/\zpar x_{n_{1}}$ span $ImH_{1}$
and $\zpar/\zpar y_{1},...,\zpar/\zpar y_{n_{2}}$ span $ImH_{2}$.
Identify $M$ to the zero section $S_{0}$ of $T^{*}M$. If
$(x,y,{\tilde x},{\tilde y})$ are the associated coordinates on $T^{*}M$,
in which the zero section is given by ${\tilde x}=0$, ${\tilde y}=0$,
from the formula of the prolongation given in sub-section 1.2 easily
follows that $G^{*}(TS_{0})\zco TS_{0}$, $G^{*}_{1}(TS_{0})\zco TS_{0}$
and $G^{*}_{2}(TS_{0})\zco TS_{0}$. Besides
$(Im\zf^{2}_{2}(G^{*}))\zin TS_{0}=Im\zf_{2}(G)$,
$(Im\zf^{2}_{1}(G^{*}))\zin TS_{0}=Im\zf_{1}(G)$,
$G^{*}_{\zbv S_{0}}=G$, ${G^{*}_{1}}_{\zbv S_{0}}=G_{1}$ and
${G^{*}_{2}}_{\zbv S_{0}}=G_{2}$ [it is just an algebraic verification
at each point of $S_{0}$]. Thus $N_{G_{1}}=N_{G_{2}}=0$ on $M$.
In particular from $N_{G_{1}}(\zpar/\zpar y_{r},\quad)=0$ follows
$L_{(\zpar/\zpar y_{r})}G_{1}=0$, $r=1,...,n_{2}$, that is $G_{1}$ does
not depend on $M_{2}$. Analogously one shows that $G_{2}$ does
not depend on $M_{1}$. Therefore {\it the proof of proposition A.1 is finished}.

\end{document}